\documentclass[11pt]{amsart}
\usepackage{amsmath, amssymb, amsthm, xcolor, mathtools}
\usepackage[margin=1in]{geometry}
\usepackage{hyperref}
\usepackage{cleveref}
\usepackage{comment}
\usepackage{float}

\newtheorem{theorem}{Theorem}[section]
\newtheorem{lemma}[theorem]{Lemma}

\theoremstyle{definition}
\newtheorem{definition}[theorem]{Definition}

\theoremstyle{remark}
\newtheorem{remark}[theorem]{Remark}
\newtheorem{remarks}[theorem]{Remarks}

\numberwithin{equation}{section}

\renewcommand\d{\,\mathrm{d}}
\newcommand{\E}{\mathbf{E}}
\newcommand{\Z}{\mathbb{Z}}
\newcommand{\R}{\mathbb{R}}
\newcommand{\N}{\mathbb{N}}
\newcommand{\Prob}{\mathbf{P}}
\newcommand{\Cov}{\mathbf{Cov}}
\newcommand{\Var}{\mathbf{Var}}
\newcommand{\eps}{\varepsilon}
\newcommand{\llll}{{\ll\!\!\ll}}
\newcommand{\asum}{\sideset{}{^{\ast}}\sum}
\let\oldpmod\pmod
\renewcommand{\pmod}[1]{\hspace{-0.12cm}\oldpmod {#1}}

\title[Quantitative correlations and prime factors]{Quantitative correlations and some problems on prime factors of consecutive integers}
\author[Terence Tao]{Terence Tao}
\address{Department of Mathematics, University of California, Los Angeles, CA 90095-1555, USA}
\email{tao@math.ucla.edu}

\author[Joni Ter\"{a}v\"{a}inen]{Joni Ter\"{a}v\"{a}inen}
\address{Department of Pure Mathematics and Mathematical Statistics, University of Cambridge, Cambridge CB3 0WB, UK}
\email{joni.p.teravainen@gmail.com}

\subjclass[2020]{11N25 (primary), 11N36, 11K65 (secondary)}

\begin{document}

\begin{abstract} We consider several old problems involving the number of prime divisors function $\omega(n)$, as well as the related functions $\Omega(n)$ and $\tau(n)$. Firstly, we show that there are infinitely many positive integers $n$ such that $\omega(n+k) \leq \Omega(n+k) \ll k$ for all positive integers $k$, establishing a conjecture of Erd\H{o}s and Straus. Secondly, we show that the series $\sum_{n=1}^{\infty} \omega(n)/2^n$ is irrational, settling a conjecture of Erd\H{o}s. Thirdly, we prove an asymptotic formula conjectured by Erd\H{o}s, Pomerance and S\'ark\"ozy for the number of $n\leq x$ satisfying $\omega(n)=\omega(n+1)$, for almost all $x$, with similar results for $\Omega$ and $\tau$. Common to the resolution of all these problems is the use of the probabilistic method. For the first problem, this is combined with computations involving a high-dimensional sieve of Maynard-type. For the second and third problems, we instead make use of a general quantitative estimate for two-point correlations of multiplicative functions with a small power of logarithm saving that may be of independent interest. This correlation estimate is derived by using recent work of Pilatte.
\end{abstract}

\maketitle

\section{Introduction}

In this paper we establish several conjectures involving the number of prime divisors function (not counting multiplicity)
\begin{equation}\label{prime-count}
\omega(n) \coloneqq \sum_{p} 1_{p\mid n},
\end{equation}
as well as the multiplicity-counting variant
\begin{equation}\label{prime-count-mult}
\Omega(n) \coloneqq \sum_p \nu_p(n),
\end{equation}
where 
 $$\nu_p(n) \coloneqq \sum_{j=1}^\infty 1_{p^j \mid n}$$
 is the number of times $p$ divides $n$. Here and in the sequel all sums and products over $p$ (or $p'$, $p_1$, etc.) are understood to be over primes.  The divisor function
\begin{equation}\label{tau-def}
\tau(n) \coloneqq \sum_d 1_{d \mid n} 
\end{equation}
is of course closely related to the above two functions, for instance obeying the identity
\begin{equation}\label{taulog}
\frac{\log \tau(n)}{\log 2} = \sum_p \frac{\log (1+\nu_p(n))}{\log 2}
\end{equation}
which in turn implies the inequalities
\begin{equation}\label{tau-bound}
 \omega(n) \leq \frac{\log \tau(n)}{\log 2} \leq \Omega(n),
\end{equation}
with equality holding in both cases if and only if $n$ is squarefree.

\subsection{Heuristic distribution}\label{heuristic-sec}

If $\mathbf{n}$ is drawn uniformly at random from the natural numbers $n \leq x$ for some large $x$, then the famous Hardy--Ramanujan and Erd\H{o}s--Kac theorems tell us that $\omega(\mathbf{n})$ and $\Omega(\mathbf{n})$ (and hence also $\log \tau(\mathbf{n})/\log 2)$) behave like a Gaussian random variable with mean and variance\footnote{See \Cref{notation-sec} for our conventions on asymptotic notation, as well as notation for iterated logarithms such as $\log_2 x = \log\log x$.} $\sim \log_2 x$.  The mean and variance of the first two variables were computed more precisely in~\cite{diaconis}, which established the bounds
\begin{align*}
\E\ \omega(\mathbf{n}) &= \log_2 x + B_1 + O\left(\frac{1}{\log x}\right) \\
\E\ \Omega(\mathbf{n}) &= \log_2 x + B_2 + O\left(\frac{1}{\log x}\right) \\
\Var\ \omega(\mathbf{n}) &= \log_2 x + B_3 + O\left(\frac{\log_2 x}{\log x}\right) \\
\Var\ \Omega(\mathbf{n}) &= \log_2 x + B_4 + O\left(\frac{\log_2 x}{\log x}\right) 
\end{align*}
where
$$ B_1 \coloneqq \gamma + \sum_p \log \left(1-\frac{1}{p}\right) + \frac{1}{p} =0.26149\dots \text{(\href{https://oeis.org/A077761}{OEIS A077761})}$$
is the Meissel--Mertens constant, and $B_2,B_3,B_4$ are the related constants
\begin{align*}
    B_2 &\coloneqq B_1 + \sum_p \frac{1}{p(p-1)} = 1.03465\dots \text{(\href{https://oeis.org/A083342}{OEIS A083342})} \\
    B_3 &\coloneqq B_1 - \frac{\pi^2}{6} - \sum_p \frac{1}{p^2} = -1.83568\dots \\
    B_4 &\coloneqq B_1 - \frac{\pi^2}{6} + \sum_p \frac{2p-1}{p(p-1)^2} = 0.76478\dots \text{(\href{https://oeis.org/A091589}{OEIS A091589})}. 
\end{align*}
The constants $B_1,B_2,B_3,B_4$ are lower order terms in the asymptotic limit $x \to \infty$, but have a noticeable impact on the numerics for small $x$, e.g., $x \leq 10^7$, due to the very slow growth of $\log_2 x$.

Consecutive values $\omega(\mathbf{n}), \omega(\mathbf{n}+1)$ would \emph{almost} be expected heuristically to behave independently of each other, but there are some negative correlations at small primes $p$, due to the fact that at most one of $\mathbf{n}$ and $\mathbf{n}+1$ can be divisible by $p$.  A routine calculation shows that
$$ -\Cov( 1_{p|\mathbf{n}}, 1_{p|\mathbf{n}+1}) \sim \frac{1}{p^2}$$
for any fixed prime $p$, and hence one would heuristically predict
$$ -\Cov( \omega(\mathbf{n}), \omega(\mathbf{n}+1) ) \sim \sum_p \frac{1}{p^2} = 0.45224\dots \text{(\href{https://oeis.org/A085548}{OEIS A085548})};$$
a similar calculation leads to the heuristic
$$ -\Cov( \Omega(\mathbf{n}), \Omega(\mathbf{n}+1) ) \sim \sum_p \frac{1}{(p-1)^2} = 1.37506\dots \text{(\href{https://oeis.org/A086242}{OEIS A086242})}.$$
The situation with $\tau$ is more complicated, due to the fact that $\log \tau(\mathbf{n})/\log 2$ is not always an integer, so that the raw computation of variance or covariance does not give a full picture of the asymptotic distribution; we will return to this point later.

Only partial progress has been made towards making the above heuristics rigorous.  For instance, in the recent work~\cite{charamaras} it was shown that any bounded observables $f(\Omega(\mathbf{n}))$, $g(\Omega(\mathbf{n}+1))$ are asymptotically decoupled if $\mathbf{n}$ is drawn from a logarithmically weighted distribution rather than from a uniform one.  Such recent results are powered by progress on (logarithmically averaged) versions of the Elliott conjecture~\cite{Tao},~\cite{TaoTeravainenGeneral},~\cite{TaoTeravainenAlmost},~\cite{Klurman},~\cite{teravainen-binary},~\cite{Helfgott-Radziwill},~\cite{Pilatte}.

\subsection{A sieve-theoretic result}

Our first result, proven in \Cref{consecutive-sec}, settles a conjecture of Erd\H{o}s and Straus~\cite{erd},~\cite[Problem \#248]{bloom}, also stated by Erd\H{o}s and Graham~\cite[p.61]{erdos-graham}:

\begin{theorem}[Erd\H{o}s \#248]\label{erd-prog}
 There exists an absolute constant $C>0$ such that for infinitely many positive integers $n$, for all positive integers $k$ we have
\begin{equation}\label{kn}
\omega(n+k) \leq \Omega(n+k) \leq C k.
\end{equation}
In particular, by~\eqref{tau-bound} we have $\tau(n+k) \leq 2^{Ck}$ for such $n$ for all positive integers $k$.
\end{theorem}

If $k$ was restricted to a bounded range then this claim would easily follow from the fundamental lemma of sieve theory; the difficulty thus lies in the unbounded nature of $k$.  Indeed, in~\cite[p. 61]{erdos-graham}, Erd\H{o}s and Graham write ``We just know too little about sieves to handle such a question (``we'' here means not just us but the collective wisdom (?) of our poor struggling human race)''.  However, it turns out that a variant of the Maynard sieve~\cite{maynard}, which was introduced several decades after that paper, can establish this result  (via an approach similar to that in~\cite{pintz},~\cite{polymath}).  

The strategy of proof can be easily described as follows.  For simplicity we begin with the task of showing~\eqref{kn} for $\omega$ and some $n$ in a large interval $[x,2x]$.  In view of the Hardy--Ramanujan law that\ $\omega(n) \approx \log_2 x$ for most $n \in [x,2x]$, the main difficulty is to establish~\eqref{kn} in the range $k \leq \frac{1}{C_0} \log_2 x$ for a large constant $C_0$.  Here we perform a Maynard-type sieve (in fact, a product of smooth one-dimensional Selberg sieves) on $[x,2x]$, sieving $n+k$ to mostly avoid being divisible by primes less than $x^{1/100^k}$ for $1 \leq k \leq \frac{1}{C_0} \log_2 x$.  As is now well known, one can compute bounded order correlations (e.g., pair correlations or fourth order correlations) of events such as $p\mid \mathbf{n}+k$ for $\mathbf{n}$ sampled using this sieve to quite high accuracy for all primes less than, say, $x^{1/100}$, and one can then establish~\eqref{kn} with positive probability on this sieve by standard applications of the moment method (we will rely primarily on the second and fourth moments).  The main technical issue is that the dimension $\lfloor \frac{1}{C_0} \log_2 x\rfloor$ of the sieve is growing slowly with $x$, in contrast with most other applications of Maynard-type sieves in which the sieve dimension is bounded; but the dimension growth turns out to be sufficiently slow that a suitably careful computation of main terms and error terms still lets one close the argument.  To pass from $\omega$ to $\Omega$, one needs to account for divisibility by squares and higher powers of $p$, but the fact that the sieve can be built out of divisibility criteria with respect to \emph{squarefree} moduli makes these additional terms relatively easy to handle.

\begin{remark} In~\cite{erdos-unconventional} (see also~\cite[Problem \#679]{bloom}), the variant question was posed whether, for any $\eps>0$, there existed infinitely many $n$ such that $\omega(n-k) < (1+\eps) \log k / \log_2 k$ for all $k < n$ that are sufficiently large depending on $\eps$.  This appears to be beyond the capability of our methods to show, as is the variant in~\cite{erdos-mathmag},~\cite{erdos-unconventional},~\cite[Problem \#647]{bloom} asking whether there is any $n >24$ such that $\tau(n-k) \leq k+2$ for all $1 \leq k \leq n$.  

In~\cite{erdos-mathmag} (see also~\cite[Problem \#413]{bloom},~\cite{erdos-unconventional},~\cite{erdos-graham},~\cite[Problem B8]{guy}, and \href{https://oeis.org/A005236}{OEIS A005236}), it was asked whether there are infinitely many $n$ which were \emph{barriers} for $\omega$ in the sense that $\omega(n-k) \leq k$ for all $1 \leq k \leq n$.  It may be that the methods of this paper can be pushed (particularly if one is more careful about how bounds depend on the sieve dimension) to establish\footnote{In fact, after the release of this preprint, bounds of the form $\omega(n \pm k) \leq \Omega(n \pm k) \ll \log k$ were established for all $2 \leq k \leq n$, by refining the methods of this paper; see \cite{lau}.} a result of this form for \emph{sufficiently large} $k$; however, to handle small values of $k$, such as $k=1,2$, it seems that the problem is of comparable difficulty to the prime tuples conjecture (or the Sophie Germain conjecture that there are infinitely many primes $p$ with $2p+1$ also prime), and we make no progress on it here.
\end{remark}

\subsection{An irrationality result}

The motivation in~\cite{erdos-graham} for asking the above question was to resolve the following irrationality problem, which we also establish in \Cref{irrational-sec}.

\begin{theorem}[Erd\H{o}s \#69]\label{thm:irrational}
The series
\begin{equation}\label{omega-sum}
\sum_{n=1}^{\infty}\frac{\omega(n)}{2^n} = \sum_p \frac{1}{2^p-1} = 0.5169428\dots \text{\rm (\href{https://oeis.org/A262153}{OEIS A262153})}
\end{equation}
is irrational.
\end{theorem}

 Recently, \Cref{thm:irrational} was shown conditionally on a form of the prime tuples conjecture in~\cite{pratt}; our result makes the irrationality unconditional.
This answers a question of Erd\H{o}s (\cite{erdos-irrat},~\cite[p. 61]{erdos-graham},~\cite[p. 102]{erdos-irrat-2},~\cite[Problem \#69]{bloom}); in~\cite[p.61]{erdos-graham} it is stated that ``it is very annoying that~\eqref{omega-sum} cannot be proven irrational'', in comparison with similar sums such as the Erd\H{o}s--Borwein constant 
$$\sum_{n=1}^{\infty}\frac{\tau(n)}{2^n} = \sum_{n=1}^\infty \frac{1}{2^n-1} = 1.606695\dots \text{(\href{https://oeis.org/A065442}{OEIS A065442})}$$ 
which has long been known to be irrational~\cite{erdos-divisor}. This result also establishes a single special case of the conjecture (\cite{erdos-irrat},~\cite[p. 62]{erdos-graham},~\cite[p. 105]{erdos-irrat-2},~\cite[Problem \#257]{bloom}) that $\sum_{n \in A} \frac{1}{2^n-1}$ is irrational for any infinite set $A \subset \N$, namely the case when $A$ is the set of primes.  (The aforementioned result from~\cite{erdos-divisor} corresponds to the case $A=\N$.)  It seems likely that the arguments in this paper
can also treat other sets $A$ that are sufficiently similar to the primes, but we do not pursue this question here.  The method can also be modified to establish the irrationality of $\sum_{n=1}^{\infty}\frac{\omega(n)}{b^n}$ for any integer base $b \geq 2$; in fact the case $b>2$ is somewhat easier due to the faster convergence of certain coefficients that arise in the analysis.  A similar argument can also establish the irrationality of
$$\sum_{n=1}^{\infty}\frac{\Omega(n)}{2^n} = \sum_{j=1}^\infty \sum_p \frac{1}{2^{p^j}-1} = 0.5895033\dots$$
(corresponding to the case of~\cite[Problem \#257]{bloom} where $A$ is the set of prime powers); we leave the details of such modifications to the interested reader.

To establish \Cref{thm:irrational}, upper bounds on consecutive values $\omega(n+1), \omega(n+2),\dots, \omega(n+H)$ of $\omega$, such as~\eqref{kn}, are insufficient; some control on the distribution of various linear combinations of such values (e.g., $\sum_{h=1}^H \frac{\omega(n+h)}{2^h}$) is needed. Furthermore, the contribution of large prime factors (e.g., $p \geq n^{1/100}$) is no longer negligible and requires more delicate analysis.  If a sufficiently quantitative form of the prime tuples conjecture is available, one can obtain quite good control on such values for selected $n$; this is the strategy taken in~\cite{pratt}.  In the absence of the prime tuples conjecture, known distributional results such as those in~\cite{mangerel},~\cite{charamaras} turned out to be too weak for our purposes.  However, after applying some manipulations of Gowers uniformity norm type to various expressions relating to~\eqref{omega-sum}, we were able to arrive at a certain linear combination of consecutive values of $\omega$ that could be adequately controlled by a combination of Erd\H{o}s--Kac-type results, as well as pairwise correlation estimates for multiplicative functions; but the latter requires a version of the recent strong quantitative results of Elliott-type due to Pilatte~\cite{Pilatte} in order to keep the error terms under control.  Indeed, the main result of~\cite{Pilatte} gives the Chowla-type bound
$$ \frac{1}{N} \sum_{n \leq N} \lambda(n) \lambda(n+1) \ll \log^{-c} N$$
for almost all scales $N$ (where ``almost all'' is in the sense of logarithmic density), where $c>0$ is a small absolute constant, improving upon previous results towards the Chowla conjecture in~\cite{Helfgott-Radziwill},~\cite{TaoTeravainenAlmost},~\cite{TaoTeravainenGeneral},~\cite{Tao}.  By using Pilatte's method in \Cref{corr-sec}, we will be able to more generally obtain an estimate of the form
$$ \frac{W}{N}\sum_{N < n \leq 2N} (g_1(n+h_1)-\delta_N) g_2(n+h_2) 1_{n \equiv b \pmod W} \ll {\mathcal L}^{-c}
$$
for all scales $N$ outside of a small exceptional set and various $1$-bounded multiplicative functions $g_1,g_2$, where $1 \leq {\mathcal L} \ll \log N$ is a parameter, $h_1,h_2$ are distinct integer shifts, $\delta_N$ is some approximation for a mean value of $g_1$ on $[N,2N]$, $c>0$ is a small absolute constant, $W, b$ are natural numbers bounded by ${\mathcal L}^c$, and one assumes that $g_1$ is either real-valued and ``equidistributed'' (has roughly the same mean on all residue classes of small modulus) or complex-valued and ``non-pretentious'' (not close to any twisted Dirichlet character $\chi(n) n^{it}$ with $\chi$  having small modulus and $t$ not too large), with the parameter ${\mathcal L}$ used to quantify imprecise terms such as ``roughly'' or ``close''.  See \Cref{thm:correlation} for a precise (but somewhat technical) statement.  

\begin{remark} As observed\footnote{\url{https://www.erdosproblems.com/forum/thread/258}} by Przemek Chojecki using GPT 5.4 Thinking, \Cref{erd-prog} also easily resolves another irrationality problem of Erd\H{o}s and Graham \cite[p. 62]{erdos-graham}, \cite[p. 80]{erdos-irrat-2}, \cite[Problem \#258]{bloom}, namely that the quantity
$$ \sum_n \frac{\tau(n)}{a_1 \dots a_n}$$
is irrational whenever $a_1,a_2,\dots$ are natural numbers going to infinity; this had been previously established in the case that the $a_n$ were non-decreasing in \cite{erdos-straus}.  We give the proof as follows.  If for contradiction this expression was an integer multiple of $\frac{1}{q}$ for some natural number $q$, then by multiplying by $a_1 \dots a_n$ we conclude that for any $n$, the quantity
$$ \sum_{k=1}^\infty \frac{\tau(n+k)}{a_{n+1} \dots a_{n+k}}$$
is a positive multiple of $1/q$ and hence at least $1/q$.  However, from \Cref{erd-prog} and the hypothesis that the $a_n$ go to infinity, one can this quantity arbitrarily small by choosing $n$ appropriately.  This gives the desired contradiction.  We remark that this argument also works more generally if we replace $\tau(n)$ with any positive quantity bounded by $\exp(O(\Omega(n))$.
\end{remark}

\subsection{Consecutive equal values of arithmetic functions}

It is the equidistributed case of \Cref{thm:correlation} that will be relevant for proving \Cref{thm:irrational}.  By using instead the non-pretentious case of this theorem, we can obtain a further application, namely an asymptotic formula for the quantities
\begin{align*}
\Prob(\omega(\mathbf{n}) = \omega(\mathbf{n}+1)) &= \frac{1}{x} |\{ n \leq x: \omega(n)=\omega(n+1)\}| \\
\Prob(\Omega(\mathbf{n}) = \Omega(\mathbf{n}+1)) &= \frac{1}{x}|\{ n \leq x: \Omega(n)=\Omega(n+1)\}| \\
\Prob(\tau(\mathbf{n}) = \tau(\mathbf{n}+1)) &= \frac{1}{x}|\{ n \leq x: \tau(n)=\tau(n+1)\}|
\end{align*}
(which are the densities for \href{https://oeis.org/A006049}{OEIS A006049}, \href{https://oeis.org/A045920}{OEIS A045920}, and \href{https://oeis.org/A005237}{OEIS A005237} respectively), where $\mathbf{n}$ is drawn uniformly from the natural numbers $n \leq x$ for some large $x$.
These three quantities (which are expected to have similar behavior) have received some attention in previous literature.   As far back as~\cite{erdos-chowla}, it was shown that all three densities are $o(1)$. The lower bound
\begin{align}\label{eq:taulower}\Prob(\tau(\mathbf{n}) = \tau(\mathbf{n}+1)) \gg \log^{-7} x\end{align}
was shown by Heath-Brown~\cite{heath-brown-divisor} (who also announced the same result with $\tau$ replaced by $\Omega$), in particular answering a question of Erd\H{o}s and Mirsky~\cite{erdos-mirsky},~\cite[Problem \# 946]{bloom} regarding whether there were infinitely many solutions to $\tau(n)=\tau(n+1)$. See also the earlier work of Spiro~\cite{spiro} on the equation $\tau(n)=\tau(n+h)$ for $h=5040$ (later generalized to arbitrary $h$ in \cite{pinner}).  Hildebrand~\cite{hildebrand} improved the lower bound~\eqref{eq:taulower} to
$$\Prob(\tau(\mathbf{n}) = \tau(\mathbf{n}+1)) \gg \log^{-3}_2 x.$$
The bounds
$$\Prob(\omega(\mathbf{n}) = \omega(\mathbf{n}+1)) \leq \Prob(\omega(\mathbf{n}) = \omega(\mathbf{n}+1)+O(1)) \asymp \log^{-1/2}_2 x$$
were established in~\cite{eps1},~\cite{eps}, with a remark that similar results hold\footnote{In the former case one replaces the additive error $+O(1)$ with a multiplicative one $\times 2^{O(1)}$.} for $\tau$ and $\Omega$.  As noted in these papers, it was observed by Hall and Tenenbaum that the Barban--Vinogradov theorem gives the slightly weaker upper bound
$$\Prob(\omega(\mathbf{n}) = \omega(\mathbf{n}+1)) \ll \log^{-1/2}_2 x \log_3 x \log^{-1}_4 x.$$
We also note that a short proof of the infinitude of solutions to $\omega(n)=\omega(n+1)$ can be found in~\cite{schlage}.

In view of the heuristics from \Cref{heuristic-sec}, one expects to have
\begin{align*}
\Var( \omega(\mathbf{n}+1)-\omega(\mathbf{n})) &\approx 2 (\log_2 x + B_5)  \\
\Var( \Omega(\mathbf{n}+1)-\Omega(\mathbf{n})) &\approx 2 (\log_2 x + B_6),
\end{align*}
where
\begin{align*}
    B_5 &\coloneqq B_3 + \sum_p \frac{1}{p^2} = -1.3834\dots,\\
    B_6 &\coloneqq B_4 + \sum_p \frac{1}{(p-1)^2} = 2.1398\dots\\
\end{align*}
which, when combined with Gaussian heuristics coming from local limit theorems, lead to the predictions
\begin{align}
\Prob(\omega(\mathbf{n}) = \omega(\mathbf{n}+1))&\approx \frac{1}{2\sqrt{\pi (\log_2 x + B_5)}} \label{b5} \\
\Prob(\Omega(\mathbf{n}) = \Omega(\mathbf{n}+1))&\approx \frac{1}{2\sqrt{\pi (\log_2 x + B_6)}}. \label{b6} 
\end{align}

In order to understand the event $\tau(\mathbf{n}) = \tau(\mathbf{n}+1)$, it turns out to be necessary to first understand the larger (and simpler) event that $\tau(\mathbf{n}+1) / \tau(\mathbf{n})$ is a power of two, which by~\eqref{taulog} is equivalent to the condition
\begin{equation}\label{psum}
\sum_p \frac{\log (1+\nu_p(\mathbf{n}+1)) - \log (1+\nu_p(\mathbf{n}))}{\log 2}  = 0 \pmod{1}.
\end{equation}
One can compute the probability of this event asymptotically, via the following definitions:

\begin{definition}[Definition of $c_\tau$]\label{ctau-def}  For each prime $p$, let $(\mathbf{a}_{0,p}, \mathbf{a}_{1,p})$ be a random variable drawn from the set
$$ \{ (0,0)\} \cup \{ (0,j): j \in \N \} \cup \{ (j,0): j \in \N\}$$
with probability distribution
\begin{align*}
    \Prob( (\mathbf{a}_{0,p}, \mathbf{a}_{1,p}) = (0,0) ) &= 1 - \frac{2}{p} \\
    \Prob( (\mathbf{a}_{0,p}, \mathbf{a}_{1,p}) = (0,j) ) &= \left(1-\frac{1}{p}\right) \frac{1}{p^j} \\
    \Prob( (\mathbf{a}_{0,p}, \mathbf{a}_{1,p}) = (j,0) ) &= \left(1-\frac{1}{p}\right) \frac{1}{p^j}    
\end{align*}
for all $j \in \N$ (the probabilities sum to $1$ by the geometric series formula), with these pairs jointly independent in $p$.  We then set $\mathbf{c}_p \in \R/\Z$ to be the random variable
\begin{align}\label{eq:cpdef}
\mathbf{c}_p \coloneqq \frac{\log(1+\mathbf{a}_{1,p}) - \log(1+\mathbf{a}_{0,p})}{\log 2} \pmod{1};    
\end{align} 
in particular the $\mathbf{c}_p$ are jointly independent in $p$, and on considering the contribution of the cases when
$(\mathbf{a}_{0,p}, \mathbf{a}_{1,p}) = (0,0), (0,1), (1,0)$ one sees that
\begin{equation}\label{cp0}
 \Prob( \mathbf{c}_p = 0) \geq 1 - \frac{2}{p^2}.
 \end{equation}
In particular, by the Borel--Cantelli lemma it is almost surely true that all but finitely many of the $\mathbf{c}_p$ vanish.  We then set
$$ c_\tau \coloneqq \Prob\left( \sum_p \mathbf{c}_p = 0\right).$$
\end{definition}

\begin{remark} The number-theoretic interpretation of these quantities is as follows.  For each prime $p$, the pair
$(\nu_p(\mathbf{n}), \nu_p(\mathbf{n}+1))$ converges in distribution to $(\mathbf{a}_{0,p}, \mathbf{a}_{1,p})$ as $x \to \infty$; for instance, the fact that $\mathbf{a}_{0,p}, \mathbf{a}_{1,p}$ cannot simultaneously be positive reflects the coprime nature of $\mathbf{n}$ and $\mathbf{n}+1$.  The random variable $\mathbf{c}_p$ is therefore the distributional limit of the contribution
 of the prime $p$ to the expression~\eqref{psum}.  By an application of dominated convergence (using~\eqref{cp0}), the quantity $c_\tau$ is thus the asymptotic probability that $\tau(\mathbf{n}+1)/\tau(\mathbf{n})$ is a power of two.  
 
 Empirically, we have $c_\tau \approx 0.4888$; see \Cref{ctau-fig}.  A lower bound for $c_\tau$ arises from computing the probability that $\tau(\mathbf{n}+1), \tau(\mathbf{n})$ are both powers of two (\href{https://oeis.org/A372690}{OEIS A372690}); this contribution can easily computed to be
$$ c_\tau^{(1)} \coloneqq \prod_p \left(1 - \frac{2}{p^2} + \sum_{j=2}^\infty \frac{2}{p^{2^j-1}} - \frac{2}{p^{2^j}} \right) = 0.44446\dots \text{(\href{https://oeis.org/A388069}{OEIS A388069})}.$$
The next most significant contribution arises from the event that $\tau(\mathbf{n}+1), \tau(\mathbf{n})$ are both three times powers of two; this contribution has the more complicated exact expression
$$ c_\tau^{(3)} \coloneqq c_\tau^{(1)} \sum_{p_0 \neq p_1} \prod_{i=0}^1 \frac{\sum_{j=0}^\infty \left(\frac{1}{p_i^{3 \times 2^j-1}} - \frac{1}{p_i^{3 \times 2^j}}\right)}{1 - \frac{2}{p_i^2} + \sum_{j=2}^\infty \left(\frac{2}{p_i^{2^j-1}} - \frac{2}{p_i^{2^j}}\right)} \approx 0.04358.$$
Together, these two terms sum to about $0.48804$, capturing most of the observed value of $c_\tau$.  As for upper bounds, one can control the probability that $\nu_p(\tau(\mathbf{n}+1))=\nu_p(\tau(\mathbf{n}))$ for all small odd $p$ (e.g., $p=3,5,7$), which empirically provides a very good upper bound and is relatively easy to calculate asymptotically as the return time of a low-dimensional random walk; again, see \Cref{ctau-fig}.
\end{remark}

\begin{figure}[H]
    \centering
\centerline{\includegraphics[width=0.7\textwidth]{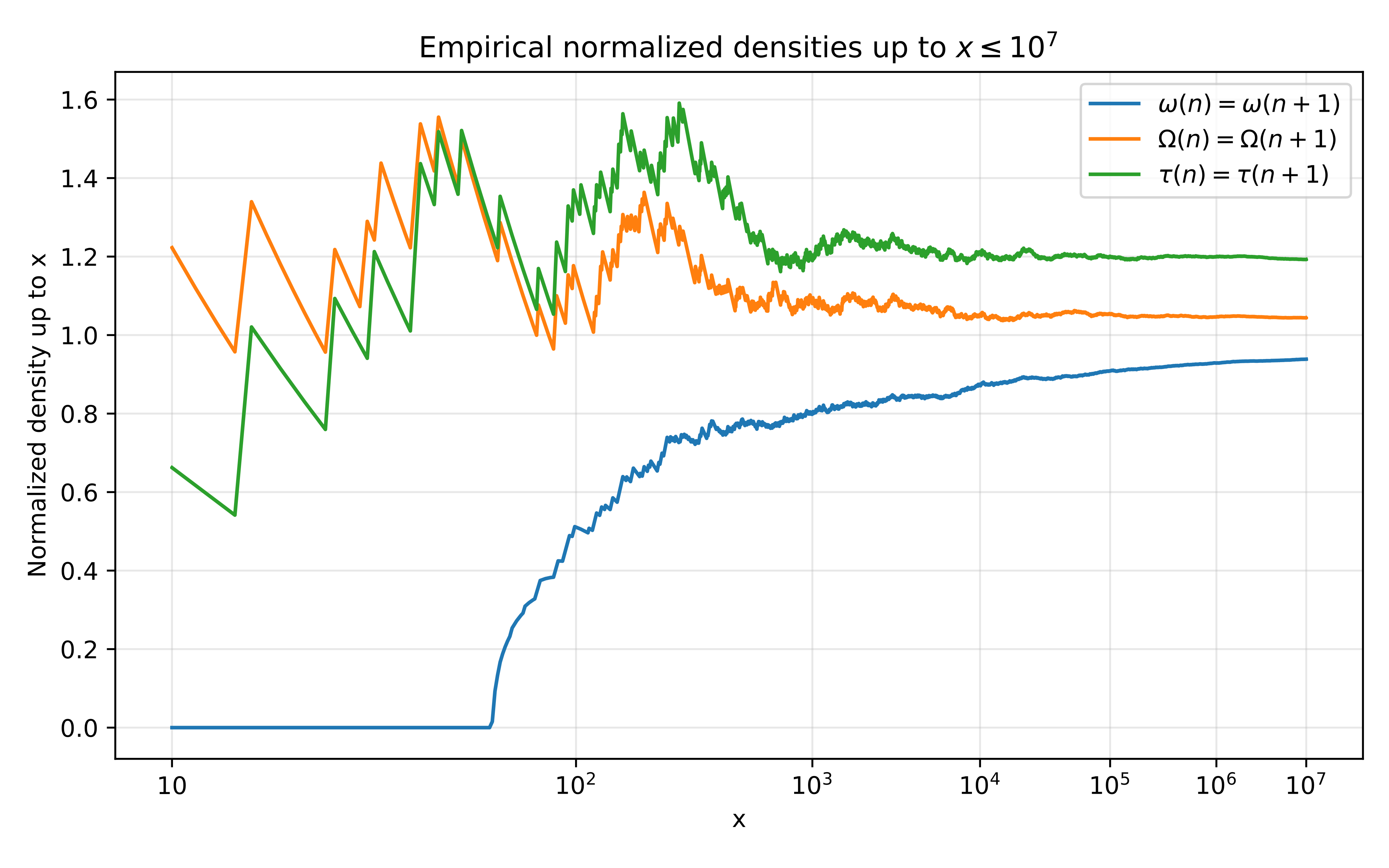}}
    \caption{The probabilities $\Prob(f(\mathbf{n}) = f(\mathbf{n}+1))$ for $f=\omega,\Omega,\tau$ and $10 \leq n \leq 10^7$, normalized by multiplying against $2\sqrt{\pi(\log_2 x + B_5)}$, $2\sqrt{\pi(\log_2 x + B_6)}$, and $2\sqrt{\pi \log_2 x}/c_\tau$ respectively, with the $x$-axis plotted using $\log_2 x$.  (Due to the negativity of $B_5$, the normalized $f=\omega$ probability vanishes for small $x$.)  \Cref{consecutive-omega-thm} asserts that these densities converge to $1$ asymptotically, although the convergence is very slow numerically, suggesting the presence of additional correction terms that decay as some negative power of $\log_2 x$, as well as some negative correction $B_7$ in the $f=\tau$ case.  }
    \label{fig:omega}
\end{figure}

\begin{figure}[H]
    \centering
\centerline{\includegraphics[width=\textwidth]{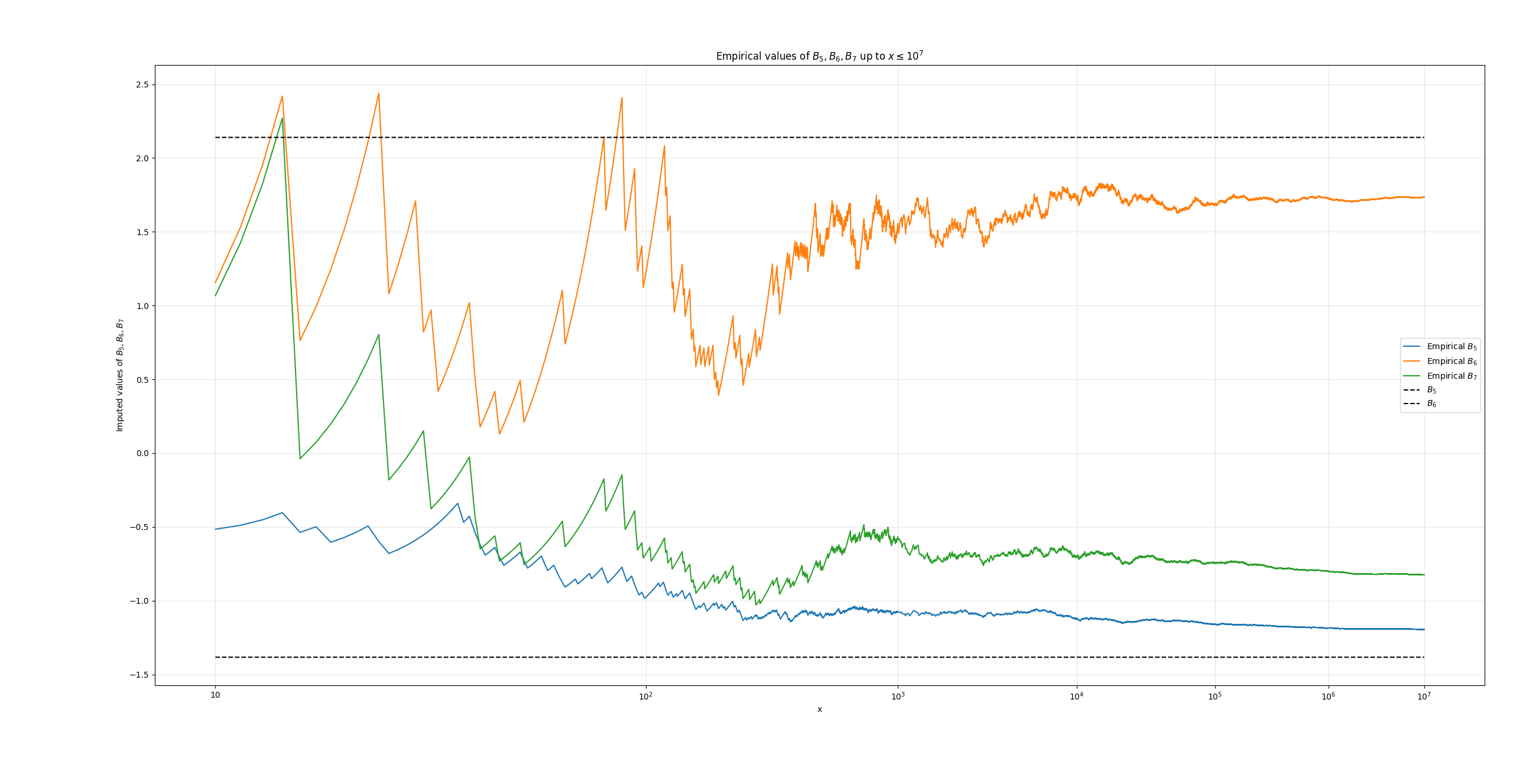}}
    \caption{Empirical values of $B_5, B_6, B_7$ imputed from treating the heuristic approximations~\eqref{b5},~\eqref{b6},~\eqref{b7} as exact, compared against the predicted values of $B_5$ and $B_6$.  These numerics tentatively suggest that $B_7 \approx -1$ to one significant figure. }
    \label{fig:imputed}
\end{figure}

\begin{figure}[H]
    \centering
\centerline{\includegraphics[width=\textwidth]{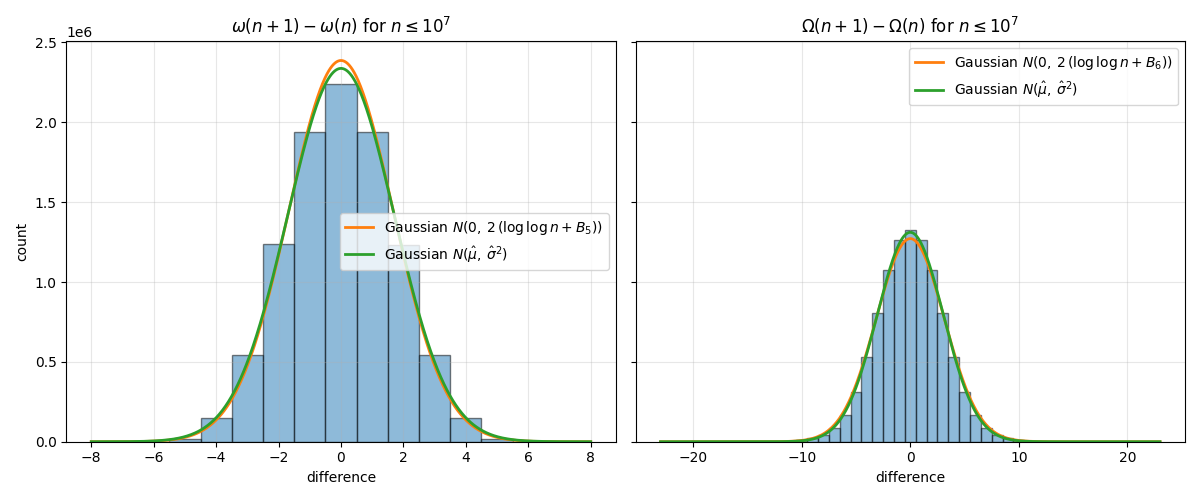}}
    \caption{Histograms of $\omega(n+1)-\omega(n)$ and $\Omega(n+1)-\Omega(n)$ for $n \leq 10^7$, compared to the Gaussian with the empirical mean $\hat \mu$ and variance $\hat \sigma^2$ for these data sets, as well as the predicted Gaussian behavior using~\eqref{b5},~\eqref{b6}. At this scale the effects of the lower order terms $B_5, B_6$ are significant. }
    \label{fig:omega-hist}
\end{figure}

\begin{figure}[H]
    \centering
\centerline{\includegraphics[width=0.7\textwidth]{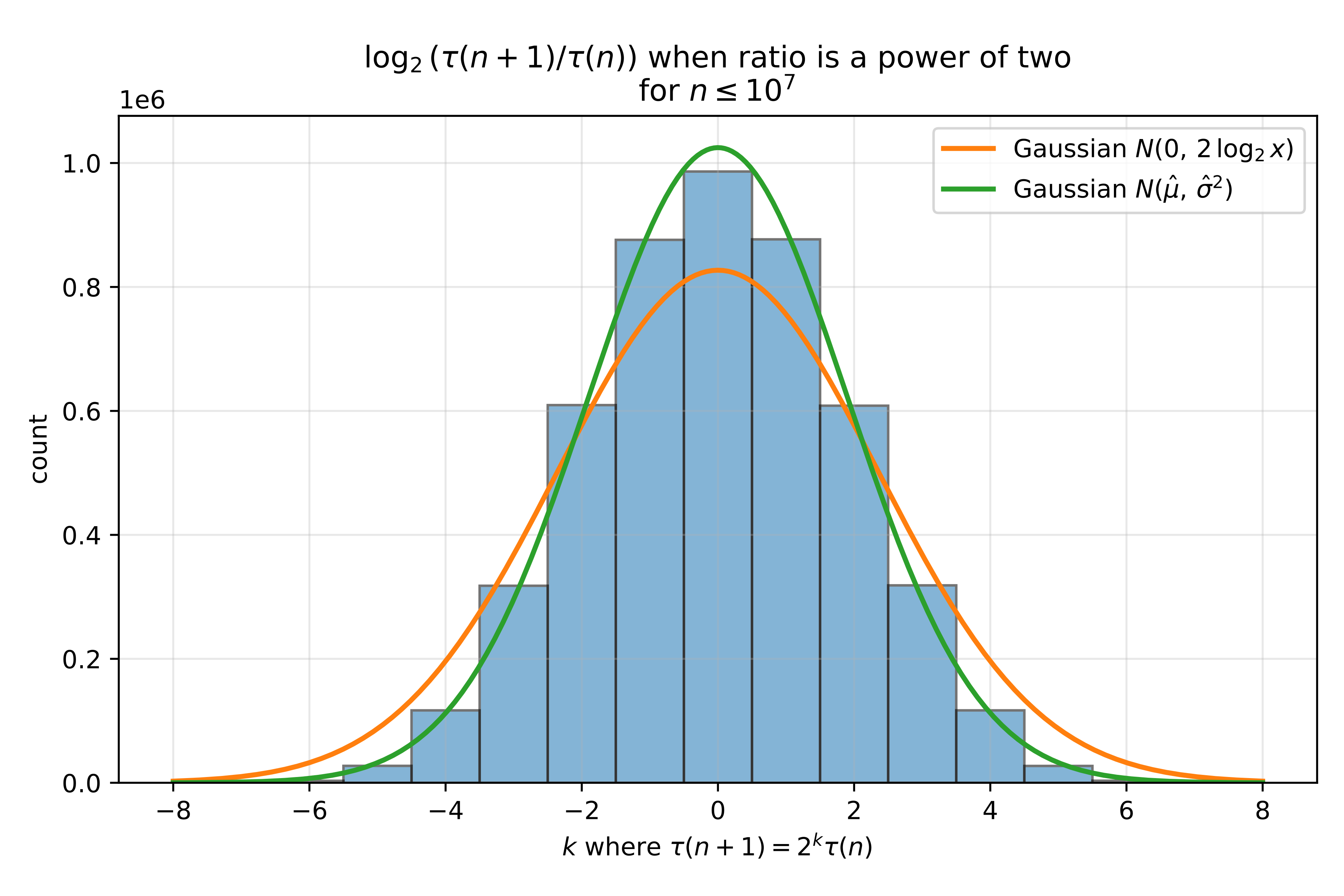}}
    \caption{Histogram of $\frac{\log \tau(n+1)-\log \tau(n)}{\log 2}$ for those $n \leq 10^7$ with $\tau(n+1)/\tau(n)$ a power of two, compared to the Gaussian with the empirical mean $\hat \mu$ and variance $\hat \sigma^2$ for this data set, as well as the predicted Gaussian behavior setting $B_7=0$.  There is a significant deviation, suggesting that the undetermined constant $B_7$ is negative.}
    \label{fig:tau-hist}
\end{figure}

It is natural to conjecture that, conditioning to the event that $\tau(\mathbf{n}+1) / \tau(\mathbf{n})$ is a power of two, that the quantity $\frac{\log \tau(\mathbf{n}+1) - \log \tau(\mathbf{n})}{\log 2}$ (which is now an integer) behaves like a Gaussian of mean zero and variance $\approx 2(\log_2 x+B_7)$ for some constant $B_7$ which could theoretically be computed, although the formula appears to be somewhat complicated, and we do not attempt to calculate it here (though numerically it seems to be roughly $-1$; see \Cref{fig:imputed}).  Thus, in analogy with~\eqref{b5},~\eqref{b6}, we have the additional heuristic prediction
 \begin{align}
\Prob(\tau(\mathbf{n}) = \tau(\mathbf{n}+1))& \approx \frac{c_\tau}{2\sqrt{\pi(\log_2 x+B_7)}}. \label{b7}
 \end{align}
Numerically, these predictions appear to be reasonably accurate (though they hint at additional lower order corrections); see \Cref{fig:omega}.  The Gaussian heuristic also appears to be strongly supported by numerics; see \Cref{fig:omega-hist}, \Cref{fig:tau-hist}.  Slightly less precise versions of these heuristics were explicitly stated as conjectures by Erd\H{o}s, Pomerance and S\'ark\"ozy~\cite{epsI} (for $\omega$) and by Hildebrand~\cite{hildebrand} (for $\tau$), with the latter noting that heuristic support for this conjecture was also found by Bateman and Spiro.

In \Cref{hildebrand-sec}, we will show the following asymptotic formulae consistent with the above predictions and numerics (though with insufficient accuracy to discern the lower order terms $B_5, B_6, B_7$), and in particular confirming the conjectures in~\cite{epsI} and~\cite{hildebrand} at almost all scales.  

\begin{theorem}[Problems of Erd\H{o}s--Pomerance--S\'ark\"ozy and Hildebrand]\label{consecutive-omega-thm}
 There is a set $\mathcal{X}\subset \mathbb{N}$ with
 $$ \frac{1}{\log X}\sum_{n \in [1,X]\setminus \mathcal{X}}\frac{1}{n} \ll \log^{-c}_2 X$$
 for some absolute constant $c>0$ and all sufficiently large $X$ (so in particular, ${\mathcal X}$ has logarithmic density $1$) such that for $x\in \mathcal{X}$ and $\mathbf{n}$ drawn uniformly at random from $[1,x]\cap \mathbb{N}$ we have the asymptotic formulae
 \begin{align*}
\Prob(\omega(\mathbf{n}) = \omega(\mathbf{n}+1))&= \frac{1+O(\log^{-c}_2 x)}{2\sqrt{\pi \log_2 x}} \\
\Prob(\Omega(\mathbf{n}) = \Omega(\mathbf{n}+1))&= \frac{1+O(\log^{-c}_2 x)}{2\sqrt{\pi \log_2 x}}  \\
\Prob(\tau(\mathbf{n}) = \tau(\mathbf{n}+1))&= \frac{c_\tau+O(\log^{-c}_2 x)}{2\sqrt{\pi\log_2 x}} 
 \end{align*}
 where $c>0$ is an absolute constant.
\end{theorem}

The proof of \Cref{consecutive-omega-thm} starts by using the circle method to express (for instance) $\Prob (\Omega(\mathbf{n})=\Omega(\mathbf{n}+1))$ as a ``low-frequency'' term that is dominated by small prime factors of $\mathbf{n}$ and $\mathbf{n}+1$, and can therefore be handled by known estimates on smooth numbers; and a ``high-frequency'' term which requires one to understand correlations of $e(\alpha \Omega(\mathbf{n}))$ and $e(\alpha \Omega(\mathbf{n}+1))$ for various frequencies $\alpha$ that are not too close to an integer.  The latter can be handled by the ``non-pretentious'' case of \Cref{thm:correlation}, with the parameter ${\mathcal L}$ depending on how close $\alpha$ is to an integer.

\subsection{Other applications}

The correlation estimate given in Theorem~\ref{thm:correlation} has also other applications; for example, it can be used to give an asymptotic formula for the number of pairs of consecutive smooth numbers. 

\begin{theorem}[Consecutive smooth numbers]\label{thm:smooth} There is a set $\mathcal{X}\subset \mathbb{N}$ with 
\begin{align}\label{eq:exceptional}
 \frac{1}{\log X}\sum_{n\in [1,X]\setminus \mathcal{X}}\frac{1}{n}\ll \log^{-c}X   
\end{align}
for some small absolute constant $c>0$ and all sufficiently large $X$ (so in particular, $\mathcal{X}$ has logarithmic density $1$) such that for any $x\in \mathcal{X}$ and any $u,v\geq 1$ we have
\begin{align}\label{eq:smoothcorr}
\frac{1}{x}\sum_{n\leq x}1_{n\textnormal{ is } x^{1/u} \textnormal{ smooth}}\, 1_{n+1\textnormal{ is } x^{1/v} \textnormal{ smooth}}=\rho(u)\rho(v)+O(\log^{-c}x),
\end{align}
where $\rho$ is the Dickman--de Bruijn function.
\end{theorem}

This strengthens~\cite[Theorem 1.14]{teravainen-binary} and~\cite[Remark 3.3]{TaoTeravainenAlmost} by giving an asymptotic formula in the range  $u,v\in [1,c(\log_2 x)/(\log_3 x)]$ for some constant $c>0$, whereas in the aforementioned works $u,v$ were fixed.

\textbf{Acknowledgements:}  TT was supported by the James and Carol Collins Chair, the Mathematical Analysis \& Application Research Fund, and by NSF grants DMS-2347850, and is particularly grateful to recent donors to the Research Fund. 

JT was supported by funding from the European Union's Horizon Europe research and innovation programme under ERC grant agreement no. 101162746. 

We thank C\'edric Pilatte for helpful discussions and suggestions, especially with regards to the modifications to his arguments required to establish Theorem \ref{thm:correlation}, Amiram Eldar for a correction and contributing a relevant OEIS sequence, and C\'ecile Dartyge for an additional correction.

ChatGPT was used to generate initial code for the graphs displayed in this paper, but the text in this paper is entirely human-generated.

\subsection{Notation}\label{notation-sec}

We use $1_E$ to denote the indicator of an event $E$, equal to $1$ when $E$ is true and $0$ otherwise.  If $S$ is a finite set, we use $|S|$ to denote its cardinality.

We let $x$ be an asymptotic parameter tending to infinity (possibly excluding a small exceptional set).  We write $X \ll Y$, $X \gg Y$, or $X = O(Y)$ to denote the estimate $|X| \leq CY$ for some constant $C$, and $X = o(Y)$ to denote the bound $|X| \leq c(x) Y$ for some $c(x)$ that goes to zero as $x \to \infty$.  We write $X \asymp Y$ for $X \ll Y \ll X$, $X \lll Y$ for $X = o(Y)$, and $X \llll Y$ for $X \ll Y^{o(1)}$.  

In later sections we shall also introduce custom variants $Y = \tilde o(Z)$  of this notation. We caution the reader that this definition will have slightly different meanings in each of the Sections~\ref{consecutive-sec},~\ref{corr-sec},~\ref{hildebrand-sec} (as this will be convenient for our calculations); see~\eqref{eq:tildedef}~\eqref{eq:tildedef2}~\eqref{eq:tildedef3}, respectively. Lastly, in Section~\ref{consecutive-sec}, we also use the variant $Y \approx Z$ introduced in~\eqref{eq:approx}.

We use the iterated logarithm notation
$$ \log_2 x \coloneqq \log\log x; \quad \log_3 x \coloneqq \log\log\log x; \quad \log_4 x \coloneqq \log\log\log\log x.$$
The argument of these iterated logarithms will always be assumed large enough that these expressions exist and are positive.  We make $\log$ take precedence over all other operators except for parentheses; thus for instance
$$ \log x \log^{-1}_2 x = (\log x) (\log \log x)^{-1}.$$
We write $e(x) \coloneqq e^{2\pi i x}$.

In addition to the arithmetic functions $\omega, \Omega, \tau$ already introduced, we shall also use the M\"obius function $\mu$ and the Euler totient function $\varphi$.

Given two multiplicative functions $f, g \colon \N \to {\mathbb C}$ and a scale $X \geq 1$, we define the \emph{pretentious distance}
$$ \mathbb{D}(f,g; X) \coloneqq \left( \sum_{p \leq X} \frac{1 - \mathrm{Re}(f(p) \overline{g(p)})}{p} \right)^{1/2}$$
and then for a modulus $Q\in \mathbb{N}$ we define the quantities
\begin{equation}\label{mdef}
M(g; X) \coloneqq \inf_{|t| \leq X} \mathbb{D}(g, n \mapsto n^{it}; X)^2
\end{equation}
and
\begin{equation}\label{mdef2}
M(g; X,Q) \coloneqq \inf_{|t| \leq X; q \leq Q; \chi \pmod{q}} \mathbb{D}(g, n \mapsto \chi(n) n^{it}; X)^2
\end{equation}
where $\chi$ ranges over Dirichlet characters of modulus $q \leq Q$.  We refer the reader to~\cite{pretentious} for a discussion of the role of pretentiousness in the theory of bounded multiplicative functions.

We use boldface notation for random variables, with $\Prob, \E, \Var, \Cov$ denoting probability, expectation, variance, and covariance respectively (with respect to a probability measure that will always be clear from context). Thus for instance $\Cov(\mathbf{X}, \mathbf{Y}) = \E(\mathbf{X} \mathbf{Y}) - (\E \mathbf{X}) (\E \mathbf{Y})$ and $\Var(\mathbf{X}) = \Cov(\mathbf{X},\mathbf{X})$.

We write $[N]$ for $[1,N]\cap \mathbb{N}$, let $d \mid n$ denote the claim that $d$ divides $n$, let $(a,b)$ denote the greatest common divisor of $a$ and $b$, and let $[a_1,\dots,a_m]$ denote the least common multiple of $a_1,\dots,a_m$. We write $\mathbb{P}$ for the set of all primes.   

A function $f \colon G \to \mathbb{C}$ is $1$-bounded if $|f(x)| \leq 1$ for all $x \in G$.
For $k\in \mathbb{N}$, a finite abelian group $G$ and a function $f\colon G\to \mathbb{C}$, we define the Gowers norm $\|f\|_{U^k(G)}$ by
\begin{align}\label{def:gowers}
\|f\|_{U^k(G)}\coloneqq \left(\mathbb{E}_{x,h_1,\ldots, h_k\in G}\prod_{\eps \in \{0,1\}^k}\mathcal{C}^{|\eps|}f(x+\eps_1h_1+\cdots+\eps_k h_k))\right)^{1/2^k},
\end{align}
where $|\eps|=|\eps_1|+\cdots +|\eps_k|$ for a vector $\eps=(\eps_1,\ldots, \eps_k)$ and $\mathcal{C}$ is the complex conjugation operator $z\mapsto \overline{z}$.

If $z>1$ is a parameter, we let $\N_{\leq z}$ and $\N_{> z}$ be the semigroups generated by primes $p \leq z$ and primes $p>z$ respectively.  Each natural number $n$ then factors uniquely as a product $n = n_{\leq z} n_{> z}$ with $n_{\leq z} \in \N_{\leq z}$ and $n_{> z} \in \N_{> z}$. If $f \colon \N \to {\mathbb C}$ is an arithmetic function, we write $f_{\leq z}(n)$ and $f_{>z}(n)$ for $f(n_{\leq z})$ and $f(n_{>z})$ respectively.  Thus for instance we have $f = f_{\leq z} f_{>z}$ if $f$ is multiplicative, and $f = f_{\leq z} + f_{>z}$ if $f$ is additive.  

Similarly, if $w, z>1$ are parameters, we let $\N_{w < \cdot \leq z}$ denote the semigroup generated by primes $w < p \leq z$, and factor $n = n_{\leq w} n_{w < \cdot \leq z} n_{> z}$ and define $f_{w < \cdot \leq z}$ as above, with the convention that $\N_{w < \cdot \leq z}$ is trivial and $n_{w < \cdot \leq z} = 1$ if $w \geq z$. Thus for instance
$$ \Omega_{w < \cdot \leq z}(n) = \sum_{w < p \leq z} \nu_p(n).$$
More generally, we have
\begin{equation}\label{fmult}
 f = f_{\leq w} f_{w < \cdot \leq z} f_{>z}
 \end{equation}
for a multiplicative function $f$, and
\begin{equation}\label{fadd}
 f = f_{\leq w} + f_{w < \cdot \leq z} +f_{>z}
 \end{equation}
 for an additive $f$.  Similarly if we partition the primes into more than three groups.

We record the following crude bound on smooth numbers: if $2 \leq z \leq x$, then we have
\begin{equation}\label{nsmooth}
| \N_{\leq z} \cap [x]| \ll x \exp\left( - c\frac{\log x}{\log z} \right)
\end{equation}
for some absolute constant $c>0$; see~\cite[(1.9)]{debruijn}.  Much sharper bounds than this are known (see, e.g.,~\cite{granville-smooth}), often involving the Dickman--de Bruijn function $\rho$, but we will not need such bounds here.

\section{A problem on prime factors of consecutive integers}\label{consecutive-sec}

In this section we establish \Cref{erd-prog}.
Let $x$ be large.
Suppose we can find a random variable $\mathbf{n}$ taking values in $[x,2x]$ with the property that for every $k \geq 1$, one has
\begin{equation}\label{omega-prob}
 \Prob (\Omega(\mathbf{n}+k) > Ck) \leq \frac{1}{10k^{3/2}}
 \end{equation}
(say) for some large constant $C$. Then, since $\sum_{k=1}^{\infty} \frac{1}{10k^{3/2}} <1$,  there exists $n \in [x,2x]$ such that
$$ \Omega(n+k) \leq C k$$
for all $k$; sending $x$ to infinity, we obtain the claim.

It remains to produce a random variable obeying~\eqref{omega-prob}.
We introduce the shift scales $1 < K < K_+$ by the formulae
\begin{align*}
K &\coloneqq \frac{1}{C_0} \log_2 x\\
 K_+ &\coloneqq  \log x.
\end{align*}
for a large constant $C_0$, and divide the shifts $k \geq 1$ into three classes:
\begin{itemize}
    \item \emph{near shifts} $1 \leq k \leq K$;
    \item \emph{far shifts} $K < k \leq K_+$; and
    \item \emph{very far shifts} $k > K_+$.
\end{itemize}
For each near or far shift $1 \leq k \leq K_+$ we introduce the scale
\begin{equation}\label{rk-def}
 R_k \coloneqq \begin{cases}
 x^{1/100^{k}},\quad &1\leq k\leq K,\\
 w,\quad &K < k \leq K_+
 \end{cases}
\end{equation}
where
\begin{align}\label{eq:w}
w\coloneqq \log^{1/2}x.     
\end{align}
Thus for any near shift $1\leq k\leq K$ (and for $x$ large enough) we have
\begin{equation}\label{logrk}
\log^{0.9} x \leq \log R_k \leq \log x.
\end{equation}

From the crude bound $\Omega(n) \ll \log n$ we have $\Omega(n+k) \ll k$ for all very far shifts $k > K_+$, so~\eqref{omega-prob} is trivial in this range.  Hence we can restrict attention to the near and far shifts $1 \leq k \leq K_+$.  It will be convenient to adopt the notation  that
\begin{align}\label{eq:tildedef}
 Y=\tilde{o}(Z)\,\, \textnormal{ means } \,\, Y=O(Z \log^{-c} x)\,\,  \textnormal{ for some absolute constant }\,\,  c>0   
\end{align}
(with $c$ independent of $C_0$). We also write 
\begin{align}\label{eq:approx}
Y \approx Z \,\,\textnormal{ for }\,\,Y = (1+\tilde o(1)) Z.   
\end{align}
A key point is that any factor of the form $O(1)^K$ can be harmlessly absorbed into a $\tilde o(\cdot)$ error term if $C_0$ is large enough. Roughly speaking, this will allow us to safely utilize $K$-dimensional sieves in our analysis.

\begin{figure}
    \centering
\centerline{\includegraphics{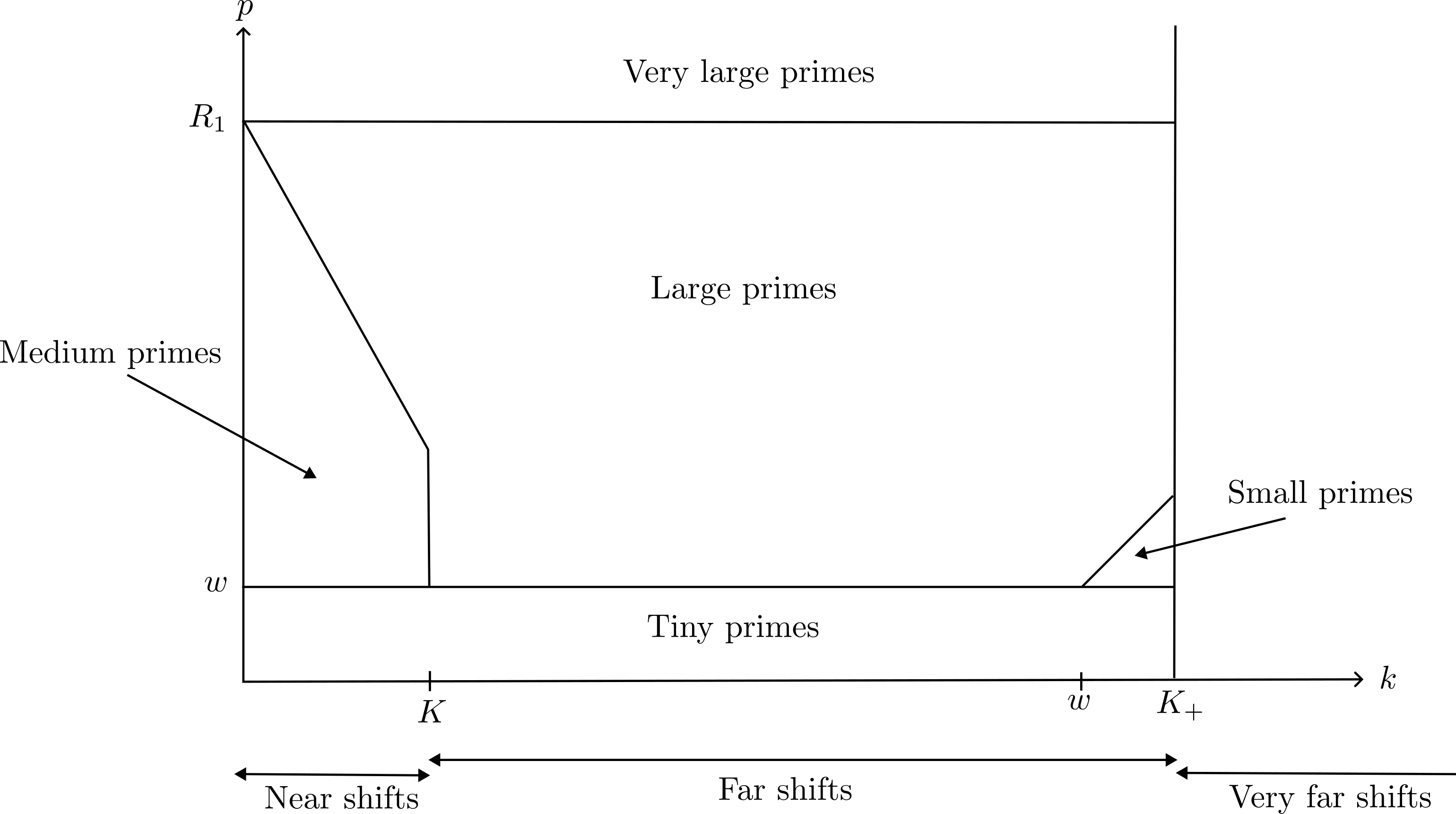}}
    \caption{A schematic diagram (not to scale) of the various ranges of primes for a given shift.  We do not perform any subdivision of primes when $k$ is a very far shift.  We perform an exact sieve at tiny primes and a Selberg type sieve at medium primes, in order to ensure $\omega(\mathbf{n}+k) \ll_{C_0} k$ with high probability for all shifts $k$.}
    \label{fig:shift}
\end{figure}

Recall the definition of $w$ in~\eqref{eq:w}.
For a given near or far shift $1 \leq k \leq K_+$, it is convenient to divide the primes into the following five classes:
\begin{itemize}
    \item \emph{tiny primes} $p \leq w$.
    \item \emph{small primes} $w < p \leq k$.  (This class will be empty if $k \leq w$.)
    \item \emph{medium primes} $w < p \leq R_k$.  (This class will be empty if $k$ is a far shift.)
    \item \emph{large primes} $\max(k,R_k) < p \leq R_1$.
    \item \emph{very large primes} $p > R_1$.
\end{itemize}
See \Cref{fig:shift}.  It will be the medium and large primes which require the most attention; the very large and tiny primes are easy to handle, and the small primes, while posing some technical difficulties with regards to sieve-theoretic estimates, will not make a major contribution to $\omega(\mathbf{n}+k)$.

Roughly speaking, the strategy is to apply a sieve weight to $n$ so that the large primes dominate the prime factorization of $n+k$ for both near and far shifts $k$.  More precisely, suppose we can construct a random variable $\mathbf{n}$ obeying the following axioms:
\begin{itemize}
    \item[(i)] Almost surely, $\mathbf{n}$ takes values in $[x,2x]$ and is divisible by $p^2$ for every tiny prime $p\leq w$.
    \item[(ii)]  If $1 \leq k \leq K_+$ is a near or far shift, $1 \leq j \leq 4$ and $\max(k,R_k) < p_1 < \dots < p_j \leq R_1$ are large primes, then
    $$ \Prob( p_1 \cdots p_j \mid  \mathbf{n}+k) \approx \frac{1}{p_1 \dots p_j}.$$
    \item[(iii)]  If $1 \leq k \leq K$ is a near shift, $1 \leq j \leq 2$, and $w < p_1 < \dots < p_j \leq R_k$ are medium primes, then
    $$ \Prob( p_1 \dots p_j \mid  \mathbf{n}+k) \ll \prod_{i=1}^j \frac{1}{p_i}\left(\frac{\log p_i}{\log R_k}\right)^2.$$
    \item[(iv)]  If $1 \leq k \leq K_+$ is a near or far shift, $1 \leq j \leq 2$, $a_1,\dots,a_j \geq 2$, and $p_1 < \dots < p_j$ are (small, medium, or large) primes with $p_i^{a_i} \leq R_1$ for all $1 \leq i \leq j$, then
    $$ \Prob( p_1^{a_1} \cdots p_j^{a_j} \mid  \mathbf{n}+k) = \frac{\Prob( p_1 \cdots p_j \mid  \mathbf{n}+k)}{p_1^{a_1-1} \dots p_j^{a_j-1}}  + O(x^{-0.1}).$$
    If the $p_i$ are assumed to be tiny primes rather than small, medium, or large primes, we instead have the variant bound
    $$ \Prob( p_1^{a_1} \cdots p_j^{a_j} \mid  \mathbf{n}+k) = \frac{\Prob( p_1^2 \cdots p_j^2 \mid  \mathbf{n}+k)}{p_1^{a_1-2} \dots p_j^{a_j-2}}  + O(x^{-0.1}).$$
    
\end{itemize}

\subsection{Proof that axioms (i)--(iv) imply the theorem.}
We first show how the axioms (i), (ii), (iii), (iv) imply~\eqref{omega-prob} for any near or far shift $1 \leq k \leq K_+$ (and a constant $C$ depending on $C_0$).  We split up the primes into the tiny and small primes $p \leq \max(k,w)$, the medium primes $w < p \leq R_k$, the large primes $\max(k,R_k) < p \leq R_1$, and the very large primes $p > R_1$.  
Since $\mathbf{n}+k \leq 3x$, the quantity $\mathbf{n}+k$ has at most $O(1)$ very large prime factors, each of which divides $\mathbf{n}+k$ at most $O(1)$ times, thus 
$$ \Omega_{>R_1}(\mathbf{n}+k) = O(1).$$
For all other primes $p \leq R_1$, we have 
$$\nu_p(\mathbf{n}+k) \ll 1_{p|\mathbf{n}+k} + \sum_{\substack{j \geq 2\\ p^j \leq R_1}} 1_{p^j|\mathbf{n}+k}$$
and thus by~\eqref{fadd}
\begin{align*}
    \Omega(\mathbf{n}+k) &= \Omega_{\leq R_1}(\mathbf{n}+k) + O(1) \\
    &\ll 
    \omega_{\leq R_1}(\mathbf{n}+k) + \sum_{j=2}^\infty \sum_{p \leq R_1^{1/j}} 1_{p^j \mid \mathbf{n}+k} + 1\\
    \\
    &=
    \omega_{\leq \max(k,w)}(\mathbf{n}+k) + \omega_{\max(k,w) < \cdot \leq R_k}(\mathbf{n}+k) + \omega_{R_k < \cdot \leq R_1}(\mathbf{n}+k) \\
    &\quad + \sum_{j=2}^\infty \sum_{p \leq R_1^{1/j}} 1_{p^j \mid \mathbf{n}+k} + 1.
\end{align*} 
From axiom (i), we have $1_{p\mid \mathbf{n}+k} = 1_{p\mid k}$ for tiny primes $p \leq w$, thus we have good control at both tiny and small primes:
$$ \omega_{\leq \max(k,w)}(\mathbf{n}+k) = \omega_{\leq w}(k) + \sum_{w < p \leq k} 1_{p\mid \mathbf{n}+k} \ll k.$$
We split
$$ \omega_{R_k < \cdot \leq R_1}(\mathbf{n}+k) = \sum_{R_k < p \leq R_1} \frac{1}{p} +
\sum_{R_k < p \leq R_1} \left( 1_{p|\mathbf{n}+k} - \frac{1}{p} \right)$$
and observe from Mertens' theorem that
\begin{equation}\label{wk}
 \sum_{R_k < p \leq R_1} \frac{1}{p} \ll \log \frac{\log R_1}{\log R_k} \ll C_0 k
 \end{equation}
(treating the case of near shifts $k \leq K$ and far shifts $K < k \leq K_+$ separately). Hence by the union bound
 it will suffice to show the medium prime bound
\begin{align*}
\Prob\left(\omega_{w < p \leq R_k}(\mathbf{n}+k)\geq C_0k\right)\ll \frac{1}{C_0^2k^2}    
\end{align*}
and the large prime bound
\begin{align*}
\Prob\left(\sum_{\max(k, R_k) < p \leq R_1} \left(1_{p\mid \mathbf{n}+k}-\frac{1}{p}\right)\geq C_0k\right)\ll \frac{1}{C_0^2k^2}        
\end{align*}
and the prime power bound
$$\Prob\left( \sum_{j=2}^\infty \sum_{p \leq R_1^{1/j}} 1_{p^j \mid \mathbf{n}+k} \geq C_0 k \right)\ll \frac{1}{C_0^2k^{3/2}}.$$
Applying  Chebyshev's inequality, it will suffice to show the second moment bound
\begin{equation}\label{second}
 \E \left|\sum_{w < p \leq R_k} 1_{p\mid \mathbf{n}+k}\right|^2 \ll 1,
 \end{equation}
 for medium primes, the fourth moment bound
\begin{equation}\label{fourth}
\E \left| \sum_{\max(k, R_k) < p \leq R_1} \left(1_{p\mid \mathbf{n}+k}-\frac{1}{p}\right)\right|^4 \ll C_0^2 k^2
\end{equation}
for large primes, and the additional second moment bound
\begin{equation}\label{second-alt}
\E \left| \sum_{j=2}^\infty \sum_{p \leq R_1^{1/j}} 1_{p^j \mid \mathbf{n}+k} \right|^2 \ll k^{1/2}
\end{equation}
for tiny, small, medium, and large primes.

We first show the second moment bound~\eqref{second} for medium primes.  We can assume that we have a near shift $k \leq K$, since otherwise there are no medium primes.  We expand the square and use axiom (iii) to bound the left-hand side by
$$ \ll \sum_{w < p \leq R_k} \frac{\log^2 p}{p \log^2 R_k} + \sum_{w < p < p' \leq R_k} \frac{\log^2 p \log^2 p'}{pp' \log^4 R_k}$$
and the claim then follows from Mertens' theorem.

We now turn to the fourth moment bound~\eqref{fourth} for large primes.  Let $\boldsymbol{\xi}_p$ be independent Bernoulli indicators for $\max(k,R_k) < p \leq R_1$ with mean $\E \boldsymbol{\xi}_p = \frac{1}{p}$.  From axiom (ii), we have
$$ \E \prod_{i=1}^j 1_{p_i\mid \mathbf{n}+k} = \E \prod_{i=1}^j \boldsymbol{\xi}_{p_i} + \tilde o\left( \frac{1}{[p_1,\dots,p_j]}\right)$$
for $1 \leq j \leq 4$ and large primes $\max(k, R_k) < p_1, \dots, p_j \leq R_1$ (not necessarily distinct).  This implies that
$$ \E \left(\prod_{i=1}^j \left(1_{p_i\mid \mathbf{n}+k}-\frac{1}{p_i}\right)\right) = \E \prod_{i=1}^j \left(\boldsymbol{\xi}_{p_i}-\frac{1}{p_i}\right) + \tilde o\left( \frac{1}{[p_1,\dots,p_j]}\right)$$
under the same hypotheses.  We can therefore write the left-hand side of~\eqref{fourth} as
\begin{align}\label{eq:p1p2p3p4}
\E \left| \sum_{\max(k, R_k) < p \leq R_1} \boldsymbol{\xi}_p-\frac{1}{p}\right|^4 + \tilde o\left( \sum_{\max(k, R_k) < p_1,p_2,p_3,p_4 \leq R_1} \frac{1}{[p_1,p_2,p_3,p_4]} \right).
\end{align}
By several applications of Mertens' theorem (dividing into $O(1)$ cases depending on what collisions occur between the $p_1,p_2,p_3,p_4$) we have
$$ \sum_{\max(k, R_k) < p_1,p_2,p_3,p_4 \leq R_1} \frac{1}{[p_1,p_2,p_3,p_4]} \leq \sum_{p_1,p_2,p_3,p_4 \leq R_1} \frac{1}{[p_1,p_2,p_3,p_4]} \ll \log^{O(1)}_2 x,$$
and so the error term in~\eqref{eq:p1p2p3p4} is $\tilde o(\log^{O(1)}_2 x)\ll 1 \ll k^2$; so we reduce to showing that
$$\E \left| \sum_{\max(k, R_k) < p \leq R_1} \boldsymbol{\xi}_p-\frac{1}{p}\right|^4 \ll C_0^2k^2.$$
But as the $\boldsymbol{\xi}_p - \frac{1}{p}$ are independent with mean zero, we can bound the left-hand side by
$$ \ll  \sum_{\max(k, R_k) < p,p' \leq R_1} \left(\E \left|\boldsymbol{\xi}_p-\frac{1}{p}\right|^2\right) \left(\E \left|\boldsymbol{\xi}_{p'}-\frac{1}{p'}\right|^2\right) 
+  \sum_{\max(k, R_k) < p \leq R_1}\E \left|\boldsymbol{\xi}_p-\frac{1}{p}\right|^4.$$
Since $\E |\boldsymbol{\xi}_p-\frac{1}{p}|^j \ll \frac{1}{p}$ for $j=2,4$, the claim now follows from~\eqref{wk}.

Finally, we turn to the additional second moment bound~\eqref{second-alt}.  By the triangle inequality, we can consider the tiny, small, medium, and large primes separately.  We first consider the contribution of large primes $\max(k, R_k) < p \leq R_1$. Expanding out the square and collecting terms involving the same large prime $p$, we can bound the left-hand side by
$$ \ll \sum_{\max(k, R_k) < p \leq R_1} \sum_{\substack{j \geq 2\\ p^j \leq R_1}} j \Prob(p^j \mid \mathbf{n}+k) +
\sum_{\max(k, R_k) < p_1 < p_2 \leq R_1} \sum_{\substack{j_1,j_2 \geq 2\\ p_1^{j_1}, p_2^{j_2} \leq R_1}} \Prob(p_1^{j_1} p_2^{j_2} \mid \mathbf{n}+k).$$
Applying axiom (iv) and $\sum_{j\geq 2}\frac{j}{p^{j-1}}\ll \frac{1}{p}$, one can bound this by
$$ \ll \sum_{\max(k, R_k) < p \leq R_1} \frac{\Prob(p \mid \mathbf{n}+k)}{p} 
+ \sum_{\max(k, R_k) < p_1 < p_2 \leq R_1} \frac{\Prob(p_1 p_2 \mid \mathbf{n}+k)}{p_1 p_2}  + R_1^2 x^{-0.1}.$$
Applying axiom (ii), this is bounded by
$$ \ll \sum_{\max(k, R_k) < p \leq R_1} \frac{1}{p^2}
+ \sum_{\max(k, R_k) < p_1 < p_2 \leq R_1} \frac{1}{p_1^2 p_2^2} + R_1^2 x^{-0.1}$$
which is certainly bounded, and thus is acceptable with room to spare.

A similar argument (using axiom (iii) instead of axiom (ii), and discarding the $\log p_i/\log R_k$ factors) treats the medium primes 
$w < p \leq R_k$, again with some room to spare.  For small primes $w < p \leq k$, we use the above arguments to arrive at an upper bound of the form
$$ \ll \sum_{w < p \leq k} \frac{ \Prob(p \mid \mathbf{n}+k)}{p}
+ \sum_{w < p_1 < p_2 \leq k} \frac{\Prob(p_1 p_2 \mid \mathbf{n}+k)}{p_1 p_2}  + k^2 x^{-0.1}.$$
Crudely bounding $\Prob(p \mid \mathbf{n}+k), \Prob(p_1 p_2 \mid \mathbf{n}+k) \leq 1$, we see from Mertens' theorem that this contribution is $O(\log_2^2 k)$ which is acceptable.

Finally we need to treat the contribution of tiny primes $p \leq w$.  Using the above arguments, but with the second case of axiom (iv) instead of the first, we arrive at the upper bound
$$ \ll \sum_{p \leq w}  \Prob(p^2 \mid \mathbf{n}+k)
+ \sum_{p_1 < p_2 \leq w} \Prob(p_1^2 p^2_2 \mid \mathbf{n}+k)  + w^2 x^{-0.1}.$$
But from axiom (i) we have
$$ \Prob(p^2 \mid \mathbf{n}+k) = 1_{p^2|k}; \quad  \Prob(p_1^2 p^2_2 \mid \mathbf{n}+k) = 1_{p_1^2 p_2^2|k}$$
so we can upper bound the first two summands by $\tau(k)$ and $\tau(k)^2$ respectively.  The claim now follows from the divisor bound.

\subsection{Setup and verification of axiom (i).}
It remains to construct $\mathbf{n}$ obeying axioms (i), (ii), (iii), (iv).  Here we will use a variant of the Maynard sieve~\cite{maynard}, and more specifically a product of $K$ smoothed Selberg sieves of Goldston--Pintz--Y{\i}ld{\i}r{\i}m-type, and perform calculations similar to those in~\cite{pintz},~\cite[Proposition 4.2]{polymath}, which roughly speaking would correspond to the case in which $K$ was replaced by a bounded quantity.

We first introduce the multiplier
$$ W \coloneqq \prod_{p \leq w} p^2$$
arising from tiny primes, noting from the prime number theorem and choice of $w$ that
$$ W \ll \exp(O(\log^{1/2} x)).$$
Next, let $\eta$ be smooth function supported on $[-1,1]$ that equals $1$ at the origin and is an element of the Gevrey class of order $2$, that is to say one has the derivative estimates
$$ \eta^{(m)}(u) \ll O(1)^m (m!)^{2}$$
for all $m \geq 0$ and $u \in \R$; the existence of such a function is standard and can be found for instance in~\cite{chung}.  By the Gevrey--Paley--Wiener theorem (see, e.g.,~\cite{gelfand}), we have a Fourier expansion
$$
\eta(u) = \int_\R \hat \eta(t) e^{-itu}\d t
$$
where the Fourier coefficients $\hat \eta(t)$ have the half-exponential decay
\begin{equation}\label{hat-decay}
\hat \eta(t) \ll \exp( - c |t|^{1/2} )
\end{equation}
for some absolute constant $c>0$.  If we then define the modified function
$$ \tilde \eta(u) \coloneqq e^{-u} \eta(u)$$
then $\tilde \eta$ is also supported on $[-1,1]$, and equals $1$ at the origin, and has the modified Fourier expansion
\begin{equation}\label{etau}
\tilde\eta(u) = \int_\R \hat \eta(t) e^{-(1+it)u}\d t
\end{equation}
We then introduce the sieve weight
$$ \nu(n) \coloneqq 1_{[x,2x]}(n) 1_{W\mid n} \prod_{1 \leq k \leq K} \left( \sum_{\substack{d \in \N_{>w}\\ d\mid n+k}} \mu(d) \tilde \eta\left( \frac{\log d}{\log R_k} \right) \right)^2,$$
which roughly speaking is sieving out medium primes while enforcing congruence conditions at tiny primes.
Clearly $\nu$ is non-negative; we will show later that it does not vanish entirely for $x$ large enough.  We let $\mathbf{n}$ be drawn using the probability density function
$$ \frac{\nu(n)}{\sum_{n'} \nu(n')}.$$
In particular, axiom (i) is obeyed.

\subsection{Verification of axioms (ii)--(iii): initial steps.} Now we verify axioms (ii) and (iii); initially we will do this in a unified fashion, and split into cases later.
Fix a near or far shift $1\leq k_* \leq K_+$, and let $d_* = p_1 \cdots p_j$ for some $0 \leq j \leq 4$ and some medium or large primes $\max(w,k_*) < p_1 < \dots < p_j \leq R_1$.  Then we have
$$ \Prob(d_*\mid \mathbf{n}+k_*) = \frac{P_{k_*}(d_*)}{P_{k_*}(1)},$$
where
\begin{equation}\label{quant}
P_{k_*}(d_*) = \sum_n \nu(n) 1_{d_*\mid n+k_*}.
\end{equation}
We now compute this quantity $P_{k_*}(d_*)$.  We can expand this as
\begin{align}\label{eq:sievesum}\begin{split}
 \sum_{d_1,\dots,d_K,d'_1,\dots,d'_K \in \N_{>w}} &\left(\prod_{k=1}^K \mu(d_k) \mu(d'_k) \tilde \eta\left(\frac{\log d_k}{\log R_k}\right) \tilde \eta\left(\frac{\log d'_k}{\log R_k}\right)\right) \\
&\quad \sum_{\substack{n \in [x,2x]\\ W\mid n}} 1_{d_*\mid n+k_*} \prod_{k=1}^K 1_{[d_k,d'_k]\mid n+k}.
\end{split}
\end{align}
The summand vanishes unless the $[d_k,d'_k]$ are pairwise coprime, and also coprime to $d_*$ if $k \neq k_*$ (since all prime factors of $d_*$ exceed both $k_*$ and $k$, since $k \leq K \leq w$).  The product of all the $d_k$ and $d'_k$, together with $d_*$ and $W$, is at most
\begin{equation}\label{xo}
 \prod_{k=1}^K R_k^2 \cdot R_1^4 \cdot W \leq x^{0.1}
\end{equation}
(say).  Thus in these cases, the inner sum evaluates to
$$ \frac{x}{W} \left(\prod_{\substack{k \leq K\\ k \neq k_*}} \frac{1}{[d_k,d'_k]}\right) \frac{1}{[d_{k_*}, d'_{k_*}, d_*]} + O( x^{0.1} )$$
with the convention that $d_k=d'_k = 1$ for $k > K$.    The contribution of the $O(x^{0.1})$ error to~\eqref{eq:sievesum} is at most
$$ \ll (\prod_{k=1}^K R_k^2) \cdot O(1)^K\cdot x^{0.1} \ll x^{0.2}$$
which in practice will be negligible.  The main term is then
$$ \frac{x}{W} \asum_{d_1,\dots,d_K,d'_1,\dots,d'_K \in \N_{>w}} \frac{\prod_{k=1}^K \mu(d_k) \mu(d'_k) \tilde \eta\left(\frac{\log d_k}{\log R_k}\right) \tilde \eta\left(\frac{\log d'_k}{\log R_k}\right)}{[d_{k_*}, d'_{k_*}, d_*] \prod_{\substack{k \leq K\\ k \neq k_*}} [d_k,d'_k]} 
$$
where the asterisk in the sum means that we restrict the $[d_k,d'_k]$ to be pairwise coprime, and also coprime to $d_*$ if $k \neq k_*$.  Using~\eqref{etau}, we conclude that
\begin{equation}\label{xw}
P_{k_*}(d_*)= \frac{x}{W} \int_\R \dots \int_\R F(\vec t, {\vec t}\,')\ \prod_{k=1}^K \hat \eta(t_k) \hat \eta(t'_k) \d t_k \d t'_k + O(x^{0.2}),
 \end{equation}
where $\vec t \coloneqq (t_1,\dots,t_K)$, ${\vec t}\,' \coloneqq (t'_1,\dots,t'_K)$, and
$$ F(\vec t, {\vec t}\,') \coloneqq 
\asum_{d_1,\dots,d_K,d'_1,\dots,d'_K \in \N_{>w}}\frac{\prod_{k=1}^K \frac{\mu(d_k) \mu(d'_k)}{d_k^{\frac{1+it_k}{\log R_k}} (d'_k)^{\frac{1+it'_k}{\log R_k}}}}{ [d_{k_*}, d'_{k_*}, d_*] \prod_{\substack{k \leq K\\ k \neq k_*}} [d_k,d'_k]}.
$$
We can factorize this latter expression as an Euler product
$$ F(\vec t, {\vec t}\,') = \prod_{p>w} E_{k_*,d_*,p}(\vec t, {\vec t}\,')$$
where the Euler factor $E_{k_*,d_*,p}(\vec t, {\vec t}\,')$ is defined by
$$E_{k_*,d_*,p}(\vec t, {\vec t}\,')= 1 - \sum_{k=1}^K \left(\frac{1}{p^{1+\frac{1+it_k}{\log R_k}}} + \frac{1}{p^{1+\frac{1+it'_k}{\log R_k}}} - \frac{1}{p^{1+\frac{2+it_k+it'_k}{\log R_k}}}\right)$$
when $p \nmid  d_*$, by
\begin{equation}\label{mix}
\begin{split}
    &
E_{k_*,d_*,p}(\vec t, {\vec t}\,')=\frac{1}{p} - \frac{1}{p^{1+\frac{1+it_{k_*}}{\log R_{k_*}}}} - \frac{1}{p^{1+\frac{1+it'_{k_*}}{\log R_{k_*}}}} + \frac{1}{p^{1+\frac{2+it_{k_*}+it'_{k_*}}{\log R_k}}}\\
&= \frac{1}{p} \left(1 - p^{-\frac{1+it_{k_*}}{\log R_{k_*}}}\right) \left(1 - p^{-\frac{1+it'_{k_*}}{\log R_{k_*}}}\right)
\end{split}
\end{equation}
when $p\mid d_*$ and $k_* \leq K$ is a near shift, and by
\begin{equation}\label{pw}
E_{k_*,d_*,p}(\vec t, {\vec t}\,')=\frac{1}{p}
\end{equation}
when $p\mid d_*$ and $K < k_*$ is a far shift (in which case $d_{k_*}=d_{k_*}'=1$).

We can crudely bound $E_{k_*,d_*,p}(\vec t, {\vec t}\,')=1 + O( K / p^{1 + \frac{1}{\log x}})$ when $p \nmid d_*$ and $E_{k_*,d_*,p}(\vec t, {\vec t}\,')=O(1/p)$ when $p \mid  d_*$.  Using the Euler product $\prod_{p}(1-p^{-1-1/\log R})^{-1}=\zeta(1+1/\log R)\ll \log R$ and recalling that $d_{*}$ is squarefree and has $O(1)$ prime factors, this gives the crude upper bound
$$ F(\vec t, {\vec t}\,') \ll \frac{\log^{O(K)} x}{d_*} \ll \frac{\exp(O(\log^2_2 x))}{d_*}.$$
If one has $|t_k| \geq \sqrt{\log x}$ or $|t_{k'}| \geq \sqrt{\log x}$ for some $1 \leq k \leq K$, one can use the decay property~\eqref{hat-decay} of $\eta$ to control this contribution to~\eqref{xw} by
$$ \frac{x}{W} \exp(-c \log^{1/4} x) O(1)^K \frac{\exp(O(\log^2_2 x))}{d_*} \ll \frac{x \exp(-c' \log^{1/4} x)}{W d_*}$$
for some constants $c,c' > 0$.  We conclude that
\begin{equation}\label{kds}
P_{k_*}(d_*) =  \frac{x}{W} \int_{|t_j|, |t'_j| \leq \sqrt{\log x}\, \forall j\in [K]} F(\vec t, {\vec t}\,')\ \prod_{k=1}^K \hat \eta(t_k) \hat \eta(t'_k) \d t_k \d t'_k + O\left( \frac{x \exp(-c' \log^{1/4} x)}{W d_*} \right)
\end{equation}
(the error term here can absorb the previous $O(x^{0.2})$ error).  

Using the Euler product $\zeta(s) = \prod_p (1-1/p^s)^{-1}$ for $\textnormal{Re}(s)>1$, applied with $s$ equal to $1+\frac{1+it_k}{\log R_k}, 1+\frac{1+it_k'}{\log R_k}, 1+\frac{2+it_k+it_k'}{\log R_k}$, we can factor $F(\vec t, {\vec t}\,')$ as 
$$ \left(\frac{W}{\varphi(W)}\right)^K \prod_{k=1}^K \frac{\zeta(1 + \frac{2+it_{k}+it'_{k}}{\log R_k})}{\zeta(1 + \frac{1+it_k}{\log R_k}) \zeta(1 + \frac{1+it'_k}{\log R_k})} \prod_p \tilde E_{k_*,d_*,p}(\vec t, {\vec t}\,'),$$
where the normalized Euler factor $\tilde E_{k_*,d_*,p}(\vec t, {\vec t}\,')$ is defined by
$$ \tilde E_{k_*,d_*,p}(\vec t, {\vec t}\,')=\prod_{k=1}^K \frac{\left(1 - \frac{1}{p^{1+\frac{2+it_k+it'_k}{\log R_k}}}\right) (1-\frac{1}{p})}{\left(1 - \frac{1}{p^{1+\frac{1+it_k}{\log R_k}}}\right) \left(1 - \frac{1}{p^{1+\frac{1+it'_k}{\log R_k}}}\right)}$$
for tiny primes $p \leq w$, and by
$$ \tilde E_{k_*,d_*,p}(\vec t, {\vec t}\,')=E_{k_*,d_*,p}(\vec t, {\vec t}\,') \prod_{k=1}^K \frac{\left(1 - \frac{1}{p^{1+\frac{2+it_k+it'_k}{\log R_k}}}\right)}{\left(1 - \frac{1}{p^{1+\frac{1+it_k}{\log R_k}}}\right) \left(1 - \frac{1}{p^{1+\frac{1+it'_k}{\log R_k}}}\right)}$$
for non-tiny primes $p > w$.
From the restriction $|t_k|, |t'_k| \leq \sqrt{\log x}$ and~\eqref{logrk}, we have
$$ \frac{1+it_k}{\log R_k}, \frac{1+it'_k}{\log R_k} \ll \log^{-0.4} x$$
for near shifts $1 \leq k \leq K$.  By Taylor expansion, this gives
$$ \tilde E_{k_*,d_*,p}(\vec t, {\vec t}\,') = 1 + \tilde o\left( \frac{1}{p} \right)$$
for tiny primes $p \leq w$,
$$ \tilde E_{k_*,d_*,p}(\vec t, {\vec t}\,') = 1 + O\left( \frac{K}{p^2} \right)$$
for non-tiny primes $p > w$ with $p \nmid d_*$, and
$$ \tilde E_{k_*,d_*,p}(\vec t, {\vec t}\,') \approx E_{k_*,d_*,p}(\vec t, {\vec t}\,') $$
for $p \mid  d_*$ (which implies $p > w$, i.e., $p$ is not tiny).  Multiplying all this together using Mertens' theorem (and the hypothesis that $d_*$ has $O(1)$ prime factors), we have
$$ \prod_p \tilde E_{k_*,d_*,p}(\vec t, {\vec t}\,') \approx \prod_{p\mid d_*} E_{k_*,d_*,p}(\vec t, {\vec t}\,').$$
Also, using $\zeta(s) = \frac{1}{s-1} (1 + O(|s-1|))$ for $|s-1|\leq 1/10$, we have
$$ \prod_{k=1}^K \frac{\zeta(1 + \frac{2+it_{k}+it'_{k}}{\log R_k})}{\zeta(1 + \frac{1+it_k}{\log R_k}) \zeta(1 + \frac{1+it'_k}{\log R_k})} 
\approx \prod_{k=1}^K \frac{(1+it_k) (1+it'_k)}{(2+it_k+it'_k) \log R_k}.$$
We thus have
\begin{equation}\label{fax}
 F(\vec t, {\vec t}\,') \approx \left(\frac{W}{\varphi(W)}\right)^K \left(\prod_{k=1}^K \frac{(1+it_k) (1+it'_k)}{(2+it_k+it'_k) \log R_k}\right)
\prod_{p\mid d_*} E_{k_*,d_*,p}(\vec t, {\vec t}\,').
\end{equation}

\subsection{Verification of axiom (ii) for far shifts \texorpdfstring{$k$}{k}.}
We now specialize to the task of proving axiom (ii) in the case $K < k_*$ is a far shift (with $k_{*}$ taking the place of $k$ in that axiom).  For this case, we see from~\eqref{pw} that we have
$$ \prod_{p\mid d_*} E_{k_*,d_*,p}(\vec t, {\vec t}\,') = \frac{1}{d_*}.$$
From~\eqref{hat-decay} we have
\begin{equation}\label{ta}
\int_{|t|, |t'| \leq \sqrt{\log x}} \frac{(1+it) (1+it')}{(2+it+it')} \hat \eta(t) \hat \eta(t')\d t \d t'  = c_0 + O(\exp(-c\log^{1/4} x)),
\end{equation}
where
\begin{align*}
c_0 &\coloneqq \int_\R\int_\R \frac{(1+it) (1+it')}{(2+it+it')} \hat \eta(t) \hat \eta(t')\d t \d t' \\
&= \int_0^\infty \int_\R \int_\R (1+it) (1+it') e^{-(2+it+it')u} \hat \eta(t) \hat \eta(t')\d t \d t' \\
&= \int_0^\infty \left(\int_\R (1+it) e^{-(1+it)u} \hat \eta(t)\d t\right)^2 \d u\\
&= \int_0^\infty \left(\frac{\mathrm{d}}{\mathrm{d}u} \int_\R e^{-(1+it)u} \hat \eta(t)\d t\right)^2\d u \\
&= \int_0^\infty \tilde \eta'(u)^2\d u \\
&> 0.
\end{align*} 
Using~\eqref{fax}, the quantity~\eqref{kds} can now be written as
$$ (1+\tilde o(1)) \frac{x}{Wd_* \prod_{k=1}^K \log R_k} \left(\frac{W}{\varphi(W)}\right)^K c_0^K + O\left(\frac{O(1)^K x \exp(-c\log^{1/4} x)}{W d_*}\right)$$
and we conclude in this case that
$$ P_{k_*}(d_*) \approx \frac{x}{Wd_* \prod_{k=1}^K \log R_k} \left(c_0 \frac{W}{\varphi(W)}\right)^K$$
(noting that the error term involving $\exp(-c\log^{1/4} x)$ may be easily absorbed in the $\approx$ error); the same argument gives
\begin{equation}\label{kstar} P_{k_*}(1) \approx \frac{x}{W\prod_{k=1}^K \log R_k} \left(c_0 \frac{W}{\varphi(W)}\right)^K
\end{equation}
and thus
$$ \Prob( d_* \mid  \mathbf{n}+k_*) \approx \frac{1}{d_*}.$$
This gives axiom (ii) in the $K < k_*$ case.  For future reference we observe that~\eqref{kstar} also holds for $k_* \leq K$, since $P_{k_*}(1)$ does not depend on $k_*$.

\subsection{Verification of axiom (iii).}
Now we establish axiom (iii) (with $k_*$ in place of $k$). In this case $k_* \leq K$ is a near shift and $w < p_1,\dots,p_j \leq R_{k_*}$, and by~\eqref{kstar} the objective is to show that
$$
P_{k_*}(d_*) \ll \frac{x\left(c_0 \frac{W}{\varphi(W)}\right)^K \prod_{i=1}^j \left(\frac{\log p_i}{\log R_{k_*}}\right)^2}{Wd_*\prod_{k=1}^K \log R_k} .$$
The error term in~\eqref{kds} is easily seen to be acceptable. Thus it suffices to show that
\begin{equation}\label{iii}
\int_{|t_j|, |t'_j| \leq \sqrt{\log x}\, \forall j\in [K]} F(\vec t, {\vec t}\,')\ \prod_{k=1}^K \hat \eta(t_k) \hat \eta(t'_k) \d t_k \d t'_k
\ll \frac{\left(c_0 \frac{W}{\varphi(W)}\right)^K \prod_{i=1}^j \left(\frac{\log p_i}{\log R_{k_*}}\right)^2}{d_*\prod_{k=1}^K \log R_k} .
\end{equation}
If $p$ divides $d_*$ (and is thus one of the $p_1,\dots,p_j$), we see from~\eqref{mix} and Taylor expansion that
\begin{equation}\label{edp}
 E_{k_*,d_*,p}(\vec t, {\vec t}\,') \ll \frac{(1+|t_{k_*}|) (1+|t'_{k_*}|) \log^2 p}{p\log^2 R_{k_*}}
 \end{equation}
and thus
$$ \prod_{p\mid d_*} E_{k_*,d_*,p}(\vec t, {\vec t}\,') \ll \frac{(1+|t_{k_*}|)^j (1+|t'_{k_*}|)^j}{d_*} \left(\prod_{i=1}^j \frac{\log p_i}{\log R_{k_*}}\right)^2.$$
The contribution of the error term implicit in~\eqref{fax} to~\eqref{iii} can then be bounded (using the rapid decrease of $\hat \eta$) by
\begin{align*}
&\ll \tilde o\left( \frac{1}{d_*}O(1)^K \left(\frac{W}{\varphi(W)}\right)^K \int_\R \dots \int_\R \prod_{k=1}^K \frac{\d t_k \d t'_k}{(1+|t_k|)^2 (1+|t'_k|)^2\log R_k} \left(\prod_{i=1}^j \frac{\log p_i}{\log R_{k_*}}\right)^2 \right) \\
&\ll \tilde o \left( \frac{1}{d_*}O(1)^K \left(\frac{W}{\varphi(W)}\right)^K \prod_{k=1}^K\frac{1}{\log R_k} \left(\prod_{i=1}^j \frac{\log p_i}{\log R_{k_*}}\right)^2 \right),
\end{align*}
which is acceptable for $C_0$ large enough (since then the factor $O(1)^K (W/\varphi(W))^K$ can be absorbed into the $\tilde{o}(\cdot)$ error.).  It remains to control contribution of the main term of~\eqref{fax} to~\eqref{iii}. Thus we need to show
\begin{align*}   
&\int_{|t_j|, |t'_j| \leq \sqrt{\log x}\, \forall j\in [K]} \left(\prod_{k=1}^K \frac{(1+it_k) (1+it'_k)}{(2+it_k+it'_k)}\right)
\prod_{p\mid d_*} E_{k_*,d_*,p}(\vec t, {\vec t}\,')\ \prod_{k=1}^K \hat \eta(t_k) \hat \eta(t'_k) \d t_k \d t'_k\\
&\quad \ll  \frac{1}{d_*} c_0^K \prod_{i=1}^j \left(\frac{\log p_i}{\log R_{k_*}}\right)^2.
\end{align*}
Using~\eqref{ta}, we may cancel out all integrals except the ones involving $t_{k_*}$ and $t'_{k_*}$, and reduce to showing that
$$\int_{|t_{k_*}|, |t'_{k_*}| \leq \sqrt{\log x}} \frac{(1+it_{k_*}) (1+it'_{k_*})}{(2+it_{k_*}+it'_{k_*})}
\prod_{p\mid d_*} E_{k_*,d_*,p}(t_{k_*}, t'_{k_*}) \hat \eta(t_{k_*}) \hat \eta(t'_{k_*})\d t_{k_*} \d t'_{k_*}
\ll  \frac{1}{d_*} c_0 \left(\prod_{i=1}^j \frac{\log p_i}{\log R_{k_*}}\right)^2$$
where by abuse of notation we have dropped the (irrelevant) variables $t_k, t'_k$ for $k \neq k_*$ from $E_{k_*,d_*,p}$.  But by~\eqref{edp} and the rapid decay of $\hat \eta$, we can bound the left-hand side by
$$\ll  \int_\R \int_\R \frac{1}{d_*} \left(\prod_{i=1}^j \frac{\log p_i}{\log R_{k_*}}\right)^2\ \frac{\d t_{k_*} \d t'_{k_*}}{(1+|t_{k_*}|)^2 (1+|t'_{k_*}|)^2},$$
which gives the desired bound (noting that $c_0$ is a positive constant independent of $x$).

\subsection{Verification of axiom (ii) (with \texorpdfstring{$k_{*}$}{k*} in place of \texorpdfstring{$k$}{k}) for near shifts \texorpdfstring{$k$}{k}.}
It remains to verify axiom (ii) in the case $k_* \leq K$ is a near shift, so now $p_1,\dots,p_j > R_{k_*}$.   By~\eqref{kstar},~\eqref{kds}, it suffices to show that
$$
\int_{|t_j|, |t'_j| \leq \sqrt{\log x}\, \forall j\in [K]} F(\vec t, {\vec t}\,')\ \prod_{k=1}^K \hat \eta(t_k) \hat \eta(t'_k) \d t_k \d t'_k 
\approx \frac{\left(c_0 \frac{W}{\varphi(W)}\right)^K}{d_* \prod_{k=1}^K \log R_k} .$$
We apply~\eqref{fax}.  The contribution of the error term can be seen to be acceptable as before, so it suffices to show that
$$ \int_{|t_j|, |t'_j| \leq \sqrt{\log x}\, \forall j\in [K]} \left(\prod_{p\mid d_*} E_{k_*,d_*,p}(\vec t, {\vec t}\,')\right) \prod_{k=1}^K \frac{(1+it_k) (1+it'_k)\hat \eta(t_k) \hat \eta(t'_k)\ \d t_k \d t'_k}{(2+it_k+it'_k)}
\approx \frac{c_0^K}{d_*}.$$
We can again use~\eqref{ta} to cancel out all integrals except the ones involving $t_{k_*}$ and $t'_{k_*}$, and reduce to showing that
$$ \int_{|t_{k_*}|, |t'_{k_*}| \leq \sqrt{\log x}} \left(\prod_{p\mid d_*} E_{k_*,d_*,p}(t_{k_*}, t'_{k_*})\right)\ \frac{(1+it_{k_*}) (1+it'_{k_*})\hat \eta(t_{k_*}) \hat \eta(t'_{k_*})\d t_{k_*} \d t'_{k_*}}{2+it_{k_*}+it'_{k_*}}
\approx \frac{c_0}{d_*};$$
note that $E_{k_*,d_*,p}$ depends only $t_{k_*},t_{k_*}$, so we have abbreviated the notation by writing only those variables.
Using~\eqref{hat-decay} and the crude bound $E_{k_*,d_*,p}(t_{k_*}, t'_{k_*}) = O(1/p)$ we can easily delete the cutoff to $\sqrt{\log x}$, and reduce to establishing
\begin{equation}\label{ut}
 \int_\R \int_\R \left(\prod_{p\mid d_*} E_{k_*,d_*,p}(t_{k_*}, t'_{k_*})\right)\ \frac{(1+it_{k_*}) (1+it'_{k_*})\hat \eta(t_{k_*}) \hat \eta(t'_{k_*})\d t_{k_*} \d t'_{k_*}}{2+it_{k_*}+it'_{k_*}}
\approx \frac{c_0}{d_*}.
\end{equation}
From~\eqref{mix} we have
$$ \prod_{p\mid d_*} E_{k_*,d_*,p}(t_{k_*}, t'_{k_*}) = \frac{1}{d_*} \sum_{\epsilon \in \{0,1\}^j} \sum_{\epsilon' \in \{0,1\}^j} (-1)^{|\epsilon|+|\epsilon'|}
\prod_{i=1}^j p_i^{-\frac{\epsilon_i (1+t_{k_*})}{\log R_{k_*}}-\frac{\epsilon'_i (1+t'_{k_*})}{\log R_{k_*}} }.$$
Writing
$$ \frac{1}{2+it_{k_*}+it'_{k_*}} = \int_0^\infty e^{-(1+it_{k_*})u- (1+it'_{k_*})u}\d u$$
we can then write the left-hand side of~\eqref{ut} as
$$
\frac{1}{d_*} \int_0^\infty \left(\sum_{\epsilon \in \{0,1\}^j} (-1)^{|\epsilon|} \int_\R \left(\prod_{i=1}^j p_i^{-\frac{\epsilon_i(1+it)}{\log R_{k_*}}}\right) e^{-(1+it)u} (1+it) \hat \eta(t)\ \d t\right)^2\d u,$$
which by~\eqref{etau} simplifies to
$$
\frac{1}{d_*} \int_0^\infty \left(\sum_{\epsilon \in \{0,1\}^j} (-1)^{|\epsilon|} \tilde \eta'\left( u + \sum_{i=1}^j \epsilon_i \frac{\log p_i}{\log R_{k_*}} \right)\right)^2\d u.$$
But as all the large primes $p_i$ exceed $R_{k_*}$, and $\tilde \eta$ is supported on $[-1,1]$, the summand vanishes unless $\epsilon=0$, in which case the expression simplifies exactly to $\frac{c_0}{d_*}$, as required.

\subsection{Verification of axiom (iv).}  We can now establish axiom (iv) (again replacing $k$ by $k_*$).  Suppose first that the $p_1,\dots,p_j$ are small, medium or large primes. If we introduce the weight
$$ f(n) \coloneqq 1_{p_1^{a_1} \dots p_j^{a_j} \mid n+k_*} - \frac{1_{p_1 \dots p_j \mid n+k_*}}{p_1^{a_1-1} \dots p_j^{a_j-1}}  $$
then it will suffice to show that
$$ \E f(\mathbf{n}) \ll x^{-0.1}$$
which by~\eqref{kstar} is equivalent to
$$ \sum_n f(n) \nu(n) \ll \frac{x^{0.9}}{W\prod_{k=1}^K \log R_k} \left(c_0 \frac{W}{\varphi(W)}\right)^K.$$
The left-hand side expands out as
\begin{align}\label{eq:sievesum-f}\begin{split}
 \sum_{d_1,\dots,d_K,d'_1,\dots,d'_K \in \N_{>w}} &\left(\prod_{k=1}^K \mu(d_k) \mu(d'_k) \tilde \eta\left(\frac{\log d_k}{\log R_k}\right) \tilde \eta\left(\frac{\log d'_k}{\log R_k}\right)\right) \\
&\quad \sum_{\substack{n \in [x,2x]\\ W\mid n}} f(n) \prod_{k=1}^K 1_{[d_k,d'_k]\mid n+k}.
\end{split}
\end{align}
The product $1_{W \mid n} \prod_{k=1}^K 1_{[d_k,d'_k]\mid n+k}$ either vanishes completely, or restricts $n$ to a residue class modulo $[d_1,\dots,d_K,d'_1,\dots,d'_K] W$, which is squarefree except at tiny primes.  The function $f$ is $1$-bounded, periodic of period at most $R_1^j$, and will have mean zero on this residue class (here it is important that none of the $p_1,\dots,p_j$ are tiny).  From this we see that
$$ \sum_{\substack{n \in [x,2x]\\ W\mid n}} f(n) \prod_{k=1}^K 1_{[d_k,d'_k]\mid n+k} \ll R_1^j.$$
By the triangle inequality, the expression~\eqref{eq:sievesum-f} is then at most
$$ \ll R_1^2 \cdots R_K^2 O(1)^K R_1^j$$
which can easily be seen to be acceptable (recall that $j \leq 2$).

If instead the primes $p_1,\dots,p_j$ are tiny, we argue similarly using
$$ f(n) \coloneqq 1_{p_1^{a_1} \dots p_j^{a_j} \mid n+k_*} - \frac{1_{p_1^2 \dots p_j^2 \mid n+k_*}}{p_1^{a_1-2} \dots p_j^{a_j-2}}$$
again noting that this function has mean zero on residue classes modulo $[d_1,\dots,d_K,d'_1,\dots,d'_K] W$.

We have now shown that axioms (i)--(iv) can be satisfied simultaneously, which was enough to conclude the proof of Theorem~\ref{erd-prog}.

\begin{remark} A more careful analysis reveals that the set of $n \leq x$ for which the conclusion of \Cref{erd-prog} holds can be lower bounded by $x \exp(-c\log_2^2 x)$ for some absolute constant $c>0$ and sufficiently large $x$.  This is consistent with the (nearly) independent Gaussian heuristics discussed in the introduction, and it is plausible that one could also use further sieve-theoretic analysis to obtain a matching upper bound (with a different value of $c$).  We will not pursue these questions here.
\end{remark}

\section{A quantitative correlation estimate}\label{corr-sec}

The objective of this section is to establish the following bound.

\begin{theorem}[A quantitative correlation estimate]\label{thm:correlation} Let $X \geq 2$.  Suppose that $g_1,g_2\colon \N\to \mathbb{C}$ are $1$-bounded multiplicative functions.  Let $1 \leq {\mathcal L} \leq \log X$ and $\delta_N$ be real numbers for all $X^{0.4} \leq N \leq X$, obeying one of the following two axioms:
\begin{itemize}
    \item[(i)] (Equidistributed case) $g_1$ is real-valued and for all $X^{0.4} \leq N \leq X$ and $a,q\in \mathbb{N}$ one has the equidistribution bound
    \begin{equation}\label{equi}
    \sum_{\substack{N<n\leq 2N\\ n\equiv a\pmod q}}g_1(n) = \frac{N}{q}\delta_N + O(N {\mathcal L}^{-1}).
    \end{equation}
    Furthermore, we assume the technical condition
\begin{equation}\label{technical-1}
g_1(p) = 1
\end{equation}
whenever $\exp(\log^{1/11} X) \leq p \leq \exp(\log^{1/10} X)$.      
    \item[(ii)] (Non-pretentious case) One has $\delta_N = 0$ for all $X^{0.4} \leq N \leq X$, as well as the non-pretentious inequality
    \begin{equation}\label{nonpret}
    \exp(M( g_1; X^2, \log^{1/125} X )) \gg {\mathcal L}.
    \end{equation}
\end{itemize}
Let $c>0$ be a sufficiently small absolute constant. Then, there exists
a set $\mathcal{E}\subset [\sqrt{X},X]$ satisfying the logarithmic density bound
\begin{align*}
\frac{1}{\log X}\int_{\mathcal{E}}\frac{\d t}{t}\ll {\mathcal L}^{-c}
\end{align*}
such that for any $W \in [{\mathcal L}^c]$ and integers $b,h_1,h_2 = O({\mathcal L}^c)$ with $h_1\neq h_2$, we have
\begin{align}\label{eq:correlation_sum}
\frac{W}{N}\sum_{N < n \leq 2N} (g_1(n+h_1)-\delta_N) g_2(n+h_2) 1_{n \equiv b \pmod W} \ll {\mathcal L}^{-c}
\end{align}
for all $N\in [\sqrt{X},X]\setminus \mathcal{E}$.
\end{theorem}

\begin{remarks}\begin{itemize}
\item  A qualitative version of this theorem (with $o(1)$ type error terms), using logarithmic averaging instead of natural averaging outside of an exceptional set, was established in~\cite[Theorem 1.3]{Tao} (in the non-pretentious case) and~\cite[Theorem 1.4]{teravainen-binary} (in the equidistributed case).  The ability to replace logarithmic averaging with natural averaging outside of an exceptional set was first observed (in the non-pretentious case) in~\cite{TaoTeravainenAlmost}.  The ${\mathcal L}^{-c}$ type error terms were first obtained by Pilatte~\cite{Pilatte} (with ${\mathcal L} = \log X$ and $W=1$) in the model case of the Liouville function $\lambda$ (which is both equidistributed and non-pretentious), which improved on ~\cite{Helfgott-Radziwill} and~\cite{Tao}; the arguments here will rely heavily upon those of Pilatte.  The technical condition~\eqref{technical-1} at small primes can most likely be weakened or even dropped with some additional work, but in any event this condition is satisfied in our applications, so we retain it for this paper.

\item Note that in the equidistributed case (i) we permit $a$ to share a common factor with $q$; thus it asserts that $g_1$ is quantitatively equidistributed in both primitive and non-primitive residue classes of moduli of size $O({\mathcal L})$.  As noted in~\cite[Remark 1.3]{teravainen-binary}, this type of definition includes all sufficiently non-pretentious $1$-bounded multiplicative functions. Note also that by setting $q=1$ in~\eqref{equi} we have
    \begin{equation}\label{deltan}
     \delta_N = \frac{1}{N} \sum_{N < n \leq 2N} g_1(n) + O({\mathcal L}^{-1}) = O(1).
    \end{equation}

\item It is also worth noting that Theorem~\ref{thm:correlation} contains as a special case the following uniform bound for Chowla-type averages: For some constant $c>0$, for any $a_1,a_2,b_1,b_2\in [1,(\log N)^{c}]\cap \mathbb{N}$ with $a_1b_2\neq a_2b_1$ we have
\begin{align*}
 \frac{1}{N}\sum_{n\leq N}\lambda(a_1n+b_1)\lambda(a_2n+b_2)\ll (\log N)^{-c}   
\end{align*}
for all $N\in [2,\infty)\setminus \mathcal{E}$ with $\mathcal{E}$ satisfying $\frac{1}{\log X}\int_{[1,X]\cap \mathcal{E}}\frac{\d t}{t}\ll (\log X)^{-c}$ for all $X\geq 2$. This follows from Theorem~\ref{thm:correlation}(ii) with $W=a_1a_2$ since the Liouville function $\lambda$ satisfies~\eqref{nonpret} by~\cite[equation (1.12)]{MRT} and since we can write $$\lambda(a_1n+b_1)\lambda(a_2n+b_2)=\lambda(a_1)\lambda(a_2)\lambda(a_1a_2n+a_2b_1)\lambda(a_1a_2n+a_1b_2).$$
\end{itemize}
\end{remarks}

\subsection{A decoupling inequality}

A key tool in the proof of \Cref{thm:correlation} will be Pilatte's decoupling inequality.

\begin{theorem}[Pilatte's decoupling inequality]\label{pilatte-thm}  Let $\eps_1>0$ be a sufficiently small absolute constant.  Let $H > 0$ be sufficiently large in terms of $\eps_1$.  Let $J \in [\eps_1^2 \log_2 H]$ and set $H_0 \coloneqq \exp(\log^{1-\eps_1} H)$ and $C \coloneqq \exp(\eps_1 (\log_2 H) / (2J))$.  For $i\in [J]$, let ${\mathcal P}_i$ be the set of all primes $p$ with
$$ C^{2i-2} < \frac{\log p}{\log H_0} < C^{2i-1}$$
(so in particular $H_0 < p < H$),
and let $V\geq 2$ be such that 
\begin{equation}\label{mert}
\sum_{p \in {\mathcal P}_i} \frac{1}{p} \leq V
\end{equation}
for all $i\in [J]$.  Then, for integer $N \geq \exp(\log^{2.01} H)$, any interval $I$ of length $N$, any residue classes $b_p\pmod{p}$ for each prime $p$, any $1$-bounded functions $g_1,g_2 \colon\Z \to {\mathbb C}$, and any indices $1\leq i_1<\cdots <i_{J'}\leq J$ with $J'\in [J]$, one has
\begin{align*}
    &\sum_{(p_1,\dots,p_{J'}) \in {\mathcal P}_{i_1} \times \cdots \times {\mathcal P}_{i_{J'}}}\left| \sum_{n \in I} \left(1_{n\equiv b_{p_1}\pmod{p_1}}-\frac{1}{p_1}\right) \cdots \left(1_{n \equiv b_{p_{J'}} \pmod{p_{J'}}}-\frac{1}{p_{J'}}\right) g_1(n) g_2(n+\sigma p_1 \cdots p_{J'})\right|\\
&\quad\quad \ll V^{0.51{J'}} N
\end{align*} 
for $\sigma \in \{-1,+1\}$.
\end{theorem}

As noted in~\cite[Remark 2.5]{Pilatte}, the trivial triangle inequality bound here is basically of the form $V^{J'} N$, so the point is the gain of $V^{-0.49{J'}}$.  

\begin{proof}  We begin with some initial reductions. By splitting both $g_1$ and $g_2$ into their real and imaginary parts and applying the triangle inequality, we may assume that $g_1,g_2$ are real-valued. By the Chinese remainder theorem, there exists some (large) positive integer $B$ such that 
$B \equiv b_p \pmod{p}$ for all primes $p < H$. Hence, by shifting the interval $I$ (and $g_1, g_2$) by $B$ we may assume that $B$ vanishes, thus we may replace the constraints $n \equiv b_{p_j} \pmod{p_j}$ by $p_j \mid n$ without loss of generality. By interchanging $g_1,g_2$ if necessary and shifting $n$ by $p_1 \cdots p_{J'}$ (causing a negligible error by the triangle inequality) we may assume $\sigma=+1$. 

Now the remaining claim is 
\begin{align}\label{eq:corr1}
    \sum_{d\in \mathcal{D}}c_d\left( \sum_{n \in I} \left(1_{p_1\mid n}-\frac{1}{p_1}\right) \cdots \left(1_{p_{J'}\mid n}-\frac{1}{p_{J'}}\right) g_1(n) g_2(n+d)\right) \ll V^{0.51{J'}} N
\end{align} 
for any signs $c_d\in \{-1,+1\}$, where $\mathcal{D}=\{p_1\cdots p_{J'}\colon p_k\in \mathcal{P}_{i_k} \textnormal{ for all } k\in [{J'}]\}$. 
In the model case where $g_1=g_2=\lambda$, $I = I_N \coloneqq [N+1,2N]$ and ${J'}=J$, this claim is essentially~\cite[Theorem 2.4]{Pilatte} (or more precisely, the key input $S_2 \ll e^{O(J)}V^{J/2}N$ in the proof of that theorem in~\cite[end of \S 3]{Pilatte}), after incorporating~\cite[Remark 2.6]{Pilatte}.  Thus, we need to adapt the arguments so that
\begin{itemize}
    \item[(i)] The functions $g_1,g_2$ are not required to be equal to each other;
    \item[(ii)]  The functions $g_1, g_2$ are not required to be equal to the Liouville function $\lambda$; 
    \item[(iii)]  The length $N$ interval $I$ is not required to be equal to $I_N$; and
    \item[(iv)]  The number ${J'}$ of sets $\mathcal{P}_{i_k}$ of primes is not required to be $J$.
\end{itemize}
Item (ii) is straightforward, as  the deduction of the decoupling inequality from the key spectral radius bound~\cite[Proposition 3.4]{Pilatte} does not use any property of $\lambda$ other than $1$-boundedness (as long as $g_1,g_2$ are equal). Item (iv) is also straightforward to implement, as the only property required of the sets $\mathcal{P}_i$ of primes in the proof (besides \cite[Lemma 2.3(b)]{Pilatte}) is that they are disjoint subsets of $(H_0,H)$ with $1\ll \sum_{p\in \mathcal{P}_i}1/p\leq V$.

For item (iii), the fact that $I=I_N$ is used in precisely two places in \cite{Pilatte}, both involving estimates from \cite{norton}.  Firstly, in~\cite[end of Section 3]{Pilatte} (and replacing $J$ with $J'$), the bound
$$ \sum_{n \in I_N} \left(\frac{\omega_{\mathcal P}(n)/V+J'}{J'}\right)^{2J'} \leq e^{O(J')} N$$
is established, where ${\mathcal P} \coloneqq \bigcup_{k \in [J']} {\mathcal P}_{i_k}$.  But on expanding out the powers of $\omega_{\mathcal P}$, and noting that $H^{2J'} \ll N^{1/2}$, say, one sees that
$$  \sum_{n \in I} \left(\frac{\omega_{\mathcal P}(n)/V+J'}{J'}\right)^{2J'}  = \sum_{n \in I_N} \left(\frac{\omega_{\mathcal P}(n)/V+J'}{J'}\right)^{2J'}+O(N^{1/2+o(1)})$$
for any interval $I$ of length $N$, and using $(a/J')^{J'}\leq e^{a}$ for any $a\geq 0$, we get 
\begin{align*}
\sum_{n \in I_N} \left(\frac{\omega_{\mathcal P}(n)/V+J'}{J'}\right)^{2J'}\leq e^{2J'}\sum_{n\in I_N}e^{2\omega_{\mathcal{P}}(n)/V}\ll e^{O(J')}N    
\end{align*}
by~\cite[Lemma (3.10)]{norton}, just as in~\cite{Pilatte}, so item (iii) is not a difficulty here.  Secondly, in \cite[Equation (23)]{Pilatte}, the bound
$$ \sum_{n \in I_N} 1_{\omega_{\mathcal P}(n) > 2JV} \ll e^{-JV/3} N$$
is obtained.  Setting $m \coloneqq \lfloor J'V \rfloor$ and applying \cite[Lemma (3.10)]{norton} one has
$$
\sum_{n \in I_N} \frac{(\omega_{\mathcal P}(n)/2)^m}{m!} \leq \sum_{n \in I_N} e^{\omega_{\mathcal P}(n)/2}
\ll e^{(e^{1/2}-1)J'V} N,
$$
whence
\begin{align*}
\sum_{n \in I_N} 1_{\omega_{\mathcal P}(n) > 2J'V} 
&\ll \frac{m!}{(J'V)^m} e^{(e^{1/2}-1)J'V} N \\
&\ll \left( \frac{m}{eJ'V} \right)^{m}m^{1/2} e^{(e^{1/2}-1)V} N \\
&\ll e^{-J'V+(e^{1/2}-1)J'V+J'V/1000} N
\end{align*}
which gives the claim since $e^{1/2}-2-1/1000< -1/3$.  But again we have
$$\sum_{n \in I} \frac{(\log 2 \omega_{\mathcal P}(n))^m}{m!}= \sum_{n \in I_N} \frac{(\log 2 \omega_{\mathcal P}(n))^m}{m!}+O(N^{1/2+o(1)}),$$
so there is no difficulty with item (iii) in this step either.

Finally, we turn to (i), which requires a little more care, as the arguments in~\cite{Pilatte} rely heavily on the fact that a certain adjacency matrix is Hermitian (i.e., that one is working with undirected graphs rather than directed graphs), and allowing $g_1 \neq g_2$ will break this symmetry.  However, this can be addressed by the standard technique\footnote{A closely related technique is to convert a non-Hermitian matrix $A$ into the larger Hermitian matrix $\begin{pmatrix} 0 & A \\ A^\dagger & 0 \end{pmatrix}$, which has the same operator norm as $A$.} of ``doubling up'' a directed graph to create an undirected graph.

We introduce a ``double cover'' $\tilde I \coloneqq I \times \{-1,+1\}$ of $I$.  If we define the weighted graph $G=(\tilde I,w_0)$ with weights
\begin{align*}
w_0((m,\sigma),(n,\sigma')) \coloneqq \begin{cases} c_d\prod_{p\mid d}\left(1_{p\mid n}-\frac{1}{p}\right),\quad &n-m=d\sigma \in \mathcal{D}, \quad \sigma \neq \sigma'\\
0,& \textnormal{otherwise},
\end{cases}    
\end{align*}
for $(m,\sigma), (n,\sigma') \in \tilde I$, then we observe the symmetry
$$ w_0((m,\sigma),(n,\sigma')) = w_0((n,\sigma'),(m,\sigma)).$$
Now, writing $A=(w_0(m,n))_{m,n\in \tilde I}$ for the weighted adjacency matrix of $G_0$ and $\boldsymbol{g_1}, \boldsymbol{g_2}$ for the column vectors
$$\boldsymbol{g_1} \coloneqq (g_1(n) 1_{\sigma=+1})_{(n,\sigma)\in \tilde I}; \quad \boldsymbol{g_2} \coloneqq (g_2(n) 1_{\sigma=-1})_{(n,\sigma)\in \tilde I},$$
we have
\begin{align*}
\boldsymbol{g_1}^{\intercal} A \boldsymbol{g_2}=\sum_{n\in I}\sum_{\substack{d\in \mathcal{D}\\ n+d\in I}}c_d\prod_{p\mid d}\left(1_{p\mid n}-\frac{1}{p}\right) g_1(n)g_2(n+d). 
\end{align*}
Note that, thanks to the introduction of the double cover, we have $g_2(n+d)$ instead of $g_2(n+\sigma d)$ here.
Since $\mathcal{D}\subset [1,H]$, removing the constraint $n+ d\in I$ causes an error term of $O(H^2)$, which is negligible, since $V\geq (\log C)/2>1$ (say) by Mertens' theorem.  It will suffice to prove an analogue of \cite[Proposition 3.4]{Pilatte}, namely that there exists $X \subset I$ with $|I_N \backslash X| \ll e^{-J' V/3} N$ such that all eigenvalues of $A|_{X \times \{-1,+1\}}$ are bounded in absolute value by $e^{O(J')} V^{J'/2}$.  By using the Ihara--Bass formula as in \cite[\S 4]{Pilatte}, one can reduce to establishing an analogue of \cite[Proposition 4.7]{Pilatte}.  Namely, there exists a set $Y \subset \Z$ with $|I \backslash Y| \ll H_0^{-1/3} N$ such that all but $N/H^3$ eigenvalues of $M_{Y \times \{-1,1\}}$ have absolute value at most $e^{O(J')} V^{J'/2}$, where the non-backtracking operator $M_{Y \times \{-1,1\}}$ is defined in exact analogy with the operator $M_Y$ defined in \cite[Definitions 4.4, 4.6]{Pilatte}.  The arguments in \cite[\S 5]{Pilatte} then reduce to establishing a high trace bound for this operator (see \cite[Proposition 5.3]{Pilatte}).  The crucial bound in \cite[Lemma 5.4]{Pilatte} continues to hold in our setting; the coefficients $c_d$ that now appear in our analogue of $M_{Y \times \{-1,1\}}$ are absorbed in the absolute value signs on the right-hand side (cf. \cite[Remark 2.6]{Pilatte}), and the passage to the double cover restricts the signs of the shifts $d_1,\dots,d_{2K}$ in \cite[(28)]{Pilatte} to now be alternating, but this only serves to improve the bound.
\end{proof}

\subsection{Initial technical reductions}

We begin with some initial technical reductions to the task of proving Theorem~\ref{thm:correlation}.  We can assume that ${\mathcal L}$ is larger than any absolute constant, as the claim is trivial otherwise.  By shifting $n, h_2, b$ by $h_1$ (and adjusting $c$ slightly), we may normalize $h_1=0$, thus $h_2 = \sigma h$ for some $h \in [{\mathcal L}^c]$ and $\sigma = \pm 1$.  If we can prove the claim for any single quadruple $(W,b,\sigma,h)$, then we can obtain it for all $(W,b,\sigma,h)$ (after replacing $c$ by, say, $c/4$), since there are only $O({\mathcal L}^{3c})$ possible choices of $(W,b,\sigma,h)$ and one can concatenate all the exceptional sets by the union bound.

In a similar spirit as~\eqref{eq:tildedef} in \Cref{consecutive-sec}, we adopt the notation 
\begin{align}\label{eq:tildedef2}
Y = \tilde o(Z)\,\, \textnormal{ if }\,\, Y =O(Z {\mathcal L}^{-c'})    
\end{align}
for some constant $c'>0$ that is independent of the $c$ parameter; in particular any factor of the form $O({\mathcal L}^{O(c)})$ can be absorbed into the $\tilde o$ notation if $c$ is small enough, and any quantity of the form $O(\log^{-c'} X)$ for some $c'>0$ independent of $c$ will qualify to be of the form $\tilde o(1)$.  By Markov's inequality, it will now suffice to establish the estimate
$$
\frac{1}{\log X} \int_{\sqrt{X}}^X \left| \frac{W}{N}\sum_{N < n \leq 2N} (g_1(n)-\delta_N) g_2(n+\sigma h) 1_{n \equiv b \pmod W} \right| \frac{\d N}{N} = \tilde o(1)
$$
for $c$ sufficiently small.  It will be convenient to treat all the $b$ at once, so we will in fact prove
\begin{equation}\label{fid}
T\coloneqq \frac{1}{\log X} \int_{\sqrt{X}}^X \sum_{b \in [W]} \left| \frac{W}{N}\sum_{N < n \leq 2N} (g_1(n)-\delta_N)g_2(n+\sigma h) 1_{n \equiv b \pmod W} \right| \frac{\d N}{N} = \tilde o(1).
\end{equation}

In preparation for applying Pilatte's decoupling inequality, we let $\eps_1\in (0,1/100)$ be a small absolute constant (not depending on $c$), and introduce the following scales:
\begin{align*}
    H &\coloneqq \exp( \log^{1/10} X ) \\
    H_0 &\coloneqq \exp(\log^{1-\eps_1} H) = \exp(\log^{(1-\eps_1)/10} X) \\
    J &\coloneqq \lfloor c \log {\mathcal L} \rfloor \\
    V &\coloneqq \frac{\eps_1 \log_2 H}{J}\\
    C &\coloneqq \exp(V/2) = \exp(\eps_1 \log_2 H / (2J))
\end{align*}
and define ${\mathcal P}_i$ as in Theorem~\ref{pilatte-thm}.  
For $c$ small enough and $X$ large enough, we can check that $J \in [\eps_1^2 \log_2 H]$, and from Mertens' theorem one has~\eqref{mert}, as well as the matching lower bound
\begin{equation}\label{mert-2}
\sum_{p \in {\mathcal P}_i} \frac{1}{p} \geq \frac{1}{3} V
\end{equation}
(say). Note that all the primes in the ${\mathcal P}_i$ lie in $(H_0,H)$ and all such primes are coprime to $W$. Also note that all the $\mathcal{P}_j$ are disjoint and that $p_1\cdots p_{J}\leq p_J^{J}\leq H$ for any primes $p_j\in \mathcal{P}_j$. Observe that any factor of the form $O(1)^J$ is of the shape
$O({\mathcal L}^{O(c)})$ and so can be absorbed into the $\tilde o$ notation if $c$ is small enough.

Let $i_1<\cdots <i_{J'}$ be all the indices $i\in [J]$ for which
\begin{align}\label{eq:Vlower}
\sum_{\substack{p\in \mathcal{P}_i\\ |g_1(p)|,|g_2(p)|\geq 1/2}}\frac{1}{p}\geq \frac{V}{10}.    
\end{align}
Note that by the triangle inequality and Hal\'asz's theorem for nonnegative multiplicative functions (see, e.g.,~\cite[Equation (2)]{hall-tenenbaum}), we have the bound
\begin{align*}
\left|\sum_{N<n\leq 2N}(g_1(n)-\delta_N)g_2(n+\sigma h)\right|&\ll \min_{j\in \{1,2\}} \sum_{N<n\leq 2N}|g_j(n)|+O(\mathcal{L}^c/N)\\
&\ll \min_{j\in \{1,2\}}N\exp\left(-\sum_{p\leq 2N}\frac{1-|g_j(p)|}{p}\right)+O(\mathcal{L}^c/N)\\
&\ll N\exp\left(-\frac{1}{2}\cdot \frac{9}{20}\cdot \frac{V}{3}(J-{J'})\right)+O(\mathcal{L}^c/N)\\
&\ll N\mathcal{L}^{-c},
\end{align*}
unless ${J'}\geq (1-20c)J$, say. Hence, we may assume that ${J'}$ obeys this lower bound for some small enough constant $c>0$.

For a given $N \in [\sqrt{X},X]$, consider the quantity
\begin{align*} S_{N,\emptyset} \coloneqq 
\sum_{\substack{(p_1,\dots,p_{J'}) \in {\mathcal P}_{i_1} \times \cdots \times {\mathcal P}_{i_{J'}}\\ |g_1(p_m)|,|g_2(p_m)|\geq 1/2\,\forall \, m\in [{J'}]}}(g(p_1)\cdots g(p_{J'}))^{-1}
\sum_{b \in [W]} &\Bigg| \frac{W}{N}\sum_{N < n \leq 2N} (g_1(n)-\delta_N) g_2(n+\sigma h p_1 \cdots p_{J'})\\
& 1_{n \equiv b \pmod W} 1_{p_1 \cdots p_{J'} \mid n} \Bigg|.
\end{align*}
Making the change of variables $n = p_1 \cdots p_{J'} n'$ (and permuting the moduli $b$), we can write 
\begin{align}\label{ppp}\begin{split}
S_{N,\emptyset}= \sum_{\substack{(p_1,\dots,p_{J'}) \in {\mathcal P}_{i_1} \times \cdots \times {\mathcal P}_{i_{J'}}\\ |g_1(p_m)|,|g_2(p_m)|\geq 1/2\,\forall \, m\in [{J'}]}} \frac{1}{p_1 \cdots p_{J'}}
\sum_{b' \in [W]} &\Bigg| \frac{W}{N/p_1 \cdots p_{J'}}\sum_{\frac{N}{p_1 \cdots p_{J'}} < n' \leq \frac{2N}{p_1 \cdots p_{J'}}} (g_1(n')-\delta_N) g_2(n'+\sigma h)\\
&1_{n' \equiv b' \pmod W} \Bigg|+O({J'}V^{{J'}-1}/H_0),
\end{split}
\end{align}
where we used the fact that the map $b\mapsto rb$ permutes the residue classes modulo $W$ for any integer $r$ coprime to $W$ and the error term arose from those $n$ for which $n$ or $n+\sigma hp_1\cdots p_J$ is divisible by $p_j^2$ for some $j\in [{J'}]$. Since $V^J\ll \exp((\log \log N)^{1+o(1)})$, the $O(\cdot)$ error term above is $\tilde o(1)$.

We claim the Lipschitz-type bound
\begin{equation}\label{lipp}
 \delta_{N/P} = \delta_N + \tilde o(1)
 \end{equation}
whenever $1 \leq P \leq H$ and $\sqrt{X} \leq N \leq X$.  This is trivial in the non-pretentious case (ii) since $\delta_{N/P}, \delta_N$ both vanish there.  In the equidistributed case (i), we observe from stability properties of mean values of real $1$-bounded multiplicative functions (see, e.g.,~\cite[(13)]{mr-annals} combined with telescopic summation) that we have
\begin{equation}\label{lipo}
\frac{1}{N} \sum_{N < n \leq 2N} g(n) - \frac{P}{N} \sum_{N/P < n \leq 2N/P} g(n) \ll \left(\frac{\log 2P}{\log N}\right)^{1/4} = \tilde o(1),
\end{equation}
and the claim follows from~\eqref{deltan}, noting that $X^{0.4} \leq N/P \leq N \leq X$.

From~\eqref{lipp}, we can replace $\delta_N$ by $\delta_{N/(p_1 \cdots p_{J'})}$ in~\eqref{ppp} with acceptable error.
Observe from the change of variables $N'=N/(p_1\cdots p_{J'})$ that for any $(p_1,\ldots, p_{J'})\in \mathcal{P}_{i_1}\times \cdots \times \mathcal{P}_{i_{J'}}$, the integral
$$ \frac{1}{\log X} \int_{\sqrt{X}}^X \sum_{b' \in [W]} \left| \frac{W}{N/p_1 \cdots p_{J'}}\sum_{\frac{N}{p_1 \cdots p_{J'}} < n' \leq \frac{2N}{p_1 \cdots p_{J'}}} (g_1(n')-\delta_{N/(p_1 \cdots p_{J'})}) g_2(n'+\sigma h) 1_{n' \equiv b' \pmod W} \right| \frac{\d N}{N}$$
is equal to $T+\tilde o(1)$.  Thus, in view of~\eqref{ppp}, to prove~\eqref{fid}, it will suffice to establish the bound
$$ \frac{1}{\log X} \int_{\sqrt{X}}^X S_{N,\emptyset}\frac{\d N}{N} \ll \tilde o\left( \sum_{\substack{(p_1,\ldots,p_{J'}) \in {\mathcal P}_{i_1} \times \cdots \times {\mathcal P}_{i_{J'}}\\ |g_1(p_m)|,|g_2(p_m)|\geq 1/2\,\forall \, m\in [{J'}]}} \frac{1}{p_1 \cdots p_{J'}} \right).$$
By~\eqref{eq:Vlower} and the fact that the $\tilde{o}(\cdot)$ notation can absorb a factor of $10^J$, it suffices to establish the pointwise estimate
\begin{equation}\label{pointsn}
S_{N,\emptyset} \ll \tilde o(V^{{J'}})
\end{equation}
for all $N \in [\sqrt{X},X]$.

By the triangle inequality, we can write
$$ S_{N,\emptyset} \leq 2^{J'}\left(S_N + \sum_{I \subset \{1,\dots,{J'}\}: I \neq \emptyset} S_{N,I}\right)$$
where
\begin{align*}
&S_N \coloneqq\\
&\sum_{(p_1,\dots,p_{J'}) \in {\mathcal P}_{i_1} \times \cdots \times {\mathcal P}_{i_{J'}}} 
\sum_{b \in [W]} \left| \frac{W}{N}\sum_{N < n \leq 2N} (g_1(n)-\delta_N) g_2(n+\sigma h p_1 \cdots p_{J'}) 1_{n \equiv b \pmod W} \prod_{j=1}^{J'} \left( 1_{p_j \mid n} - \frac{1}{p_j} \right) \right|
\end{align*}
and
\begin{align*}
&S_{N,I} \coloneqq\\
&\sum_{(p_1,\dots,p_{J'}) \in {\mathcal P}_{i_1} \times \cdots \times {\mathcal P}_{i_{J'}}} \frac{1}{\prod_{j \in I} p_j} 
\sum_{b \in [W]} \left| \frac{W}{N}\sum_{N < n \leq 2N} (g_1(n)-\delta_N) g_2(n+\sigma h p_1 \cdots p_{J'}) 1_{n \equiv b \pmod W} 1_{\prod_{j \notin I} p_j \mid n} \right|.
\end{align*}
We claim that
\begin{align}\label{eq:SNgoal}
S_N \ll \tilde o(V^{J'}).    
\end{align} 
By splitting the sum in the definition of $S_N$ into residue classes modulo $h$, it suffices to show that
\begin{align*}
S_N &\coloneqq \sum_{(p_1,\dots,p_{J'}) \in {\mathcal P}_{i_1} \times \cdots \times {\mathcal P}_{i_{J'}}} 
\sum_{b \in [W]}\sum_{a\in [h]} \left| \frac{W}{N}\sum_{N/h < n' \leq 2N/h} (g_1(hn'+a)-\delta_N) g_2(h(n'+\sigma p_1 \cdots p_{J'})+a)\right.\\
&\quad \quad \quad \quad \quad \quad \quad \quad\quad \quad \quad\quad \quad \quad\left. 1_{hn'+a \equiv b \pmod W} \prod_{j=1}^J \left( 1_{hn'+a\equiv b\pmod{p_j}} - \frac{1}{p_j} \right) \right|\ll \tilde o(V^{J'}).
\end{align*}
Now,~\eqref{eq:SNgoal} follows from Theorem~\ref{pilatte-thm} applied to the functions $n\mapsto (g_1(hn+a)-\delta_N)1_{hn+a\equiv b\pmod W}$, $n\mapsto g_2(hn+a)$ for each $a\in [h], b\in [W]$, and using that $V^{-0.49{J'}}=\tilde o(1)$ for ${J'}\geq (1-20c)J$.

 Since the $\tilde o(\cdot)$ notation can absorb factors of $2^J$, it thus suffices to show that
\begin{align}\label{eq:SNI}
S_{N,I} \ll \tilde o(V^{|I|})    
\end{align} 
for all non-empty $I \subset \{1,\dots,{J'}\}$.

\subsection{Reduction to a short exponential sum estimate}

Fix $I$.  By summing over the primes $p_j\in \mathcal{P}_{i_j}$, $j \notin I$ using~\eqref{mert}, proving~\eqref{eq:SNI} reduces to showing that for each choice of $p_j\in \mathcal{P}_{i_j}$, $j \notin I$ we have
$$
\sum_{p_i \in {\mathcal P}_{i_k} \forall k\in I} \frac{1}{\prod_{j \in I} p_j} 
\sum_{b \in [W]} \left| \frac{W}{N}\sum_{N < n \leq 2N} (g_1(n)-\delta_N) g_2(n+\sigma h p_1 \cdots p_{J'}) 1_{n \equiv b \pmod W} 1_{\prod_{j \notin I} p_j \mid n} \right|
\ll \tilde o\left(\frac{V^{|I|}}{\prod_{j \notin I} p_j}\right).
$$
Making the change of variables $n = P n'$ with $P \coloneqq \prod_{j \notin I} p_j$, and permuting the indices $b \in [W]$, this is
$$
\sum_{p_i \in {\mathcal P}_{i_k} \forall k\in I} \frac{1}{\prod_{j \in I} p_j} 
\sum_{b' \in [W]} \left| \frac{W}{N'}\sum_{N' < n' \leq 2N'} (g_1(P n')-\delta_N) g_2(P (n' +\sigma h \prod_{j \in I} p_j)) 1_{n' \equiv b' \pmod W} \right|
\ll \tilde o\left(V^{|I|}\right),
$$
where $N' \coloneqq N/P$.  Using~\eqref{lipp}, we may replace $\delta_N$ by $\delta_{N'}$ with acceptable error.  As the factors of $W$ can be absorbed into the $\tilde o(\cdot)$ error, it suffices to show that
$$
\sum_{p_i \in {\mathcal P}_{i_k} \forall k\in I} \frac{1}{\prod_{j \in I} p_j} 
 \left| \frac{1}{N'}\sum_{N' < n' \leq 2N'} (g_1(P n')-\delta_{N'}) g_2(P (n' +\sigma h \prod_{j \in I} p_j)) 1_{n' \equiv b' \pmod W} \right|
\ll \tilde o\left(V^{|I|}\right)
$$
for any given $b' \in [W]$.

As $P$ has ${J'}-|I|$ prime factors, all of which lie in $(H_0,H)$, it is a routine matter to check by the triangle inequality that the contribution of the cases where $n'$ shares a common factor with $P$ is negligible, and similarly if $n' + \sigma h \prod_{j \in I} p_j$ shares a common factor with $P$.  In such cases, since either~\eqref{technical-1} holds or $\delta_{N'}=0$, we have
$$ g_1(Pn') - \delta_{N'} = g_1(P) (g_1(n')-\delta_{N'}); \quad g_2(P (n' +\sigma h \prod_{j \in I} p_j)) = g_2(P) g_2(n' +\sigma h \prod_{j \in I} p_j)  $$
and we reduce to showing that
$$
\sum_{p_i \in {\mathcal P}_{i_k} \forall k\in I} \frac{1}{\prod_{j \in I} p_j} 
 \left| \frac{1}{N'}\sum_{N' < n' \leq 2N'} (g_1(n')-\delta_{N'}) g_2(n' +\sigma h \prod_{j \in I} p_j) 1_{n' \equiv b' \pmod W} \right|
\ll \tilde o\left(V^{|I|}\right).
$$
The left-hand side is essentially an expression of the form considered in~\cite[Proposition C.1]{Pilatte} (see also~\cite[Lemma 3.6]{Tao}); the functions $(g_1(n')-\delta_{N'}) 1_{n' \equiv b' \pmod W}$ and $g_2$ appearing here are replaced with the Liouville function, and the dilation $\sigma h$ is not present, but the proof of this proposition (based on the dyadic decomposition and the circle method) is not sensitive to these changes (In particular, the proof of~\cite[Lemma C.2]{Pilatte} goes through verbatim). By repeating the  arguments in the proof of that proposition in~\cite[Appendix C]{Pilatte}, we reduce to showing that
$$ \frac{1}{N'} \sum_{N' < n \leq 2N'} \left|\frac{1}{M} \sum_{n \leq n' \leq n+M} (g_1(n')-\delta_{N'})  1_{n' \equiv b' \pmod W} e(\alpha n)\right| \ll \tilde o(1)$$
for all $H_0\leq M \leq H$ and $\alpha \in \R$.  We can Fourier expand $1_{n' \equiv b' \pmod W}$ to absorb it in the $e(n\alpha)$ term (absorbing factors of $W$ into the $\tilde o(\cdot)$ notation), thus we reduce to showing that
\begin{equation}\label{short}
\frac{1}{N'} \sum_{N' < n \leq 2N'} \left|\frac{1}{M} \sum_{n \leq n' \leq n+M} (g_1(n')-\delta_{N'}) e(\alpha n)\right| \ll \tilde o(1).
\end{equation}

\subsection{Proof of short exponential sum estimate}

It remains to prove~\eqref{short}. In the non-pretentious case (ii), this follows from~\cite[Theorem 1.7]{MRT}.  The equidistributed case was previously considered in~\cite[Lemma 3.6, Lemma 3.7]{teravainen-binary}, but with slightly different quantitative conclusions than are needed here.  We will therefore adapt the arguments from these papers as follows.
We introduce the scales $1 \ll Q_0 \ll Q \ll M$ by the formulae
\begin{align*}
    Q_0 &\coloneqq {\mathcal L}^{1/10} \\
    Q &\coloneqq M {\mathcal L}^{-1/10}
\end{align*}

We split into two cases:
\begin{itemize}
    \item (Major arc case) One has $\alpha \in \left[\frac{a}{q}-\frac{1}{qQ},\frac{a}{q}+\frac{1}{qQ}\right]$ for some $1 \leq q \leq Q_0$ and some $1 \leq a \leq q$ with $(a,q)=1$.
    \item (Minor arc case) One has $\alpha \notin \left[\frac{a}{q}-\frac{1}{qQ},\frac{a}{q}+\frac{1}{qQ}\right]$ whenever $1 \leq q \leq Q_0$ and $1 \leq a \leq q$.
\end{itemize}

\subsection{The minor arc case}

Suppose first that we are in the minor arc case.  Then by the geometric series formula we have
\begin{align*}
 \left|\sum_{x<n\leq x+M}e(\alpha n)\right|\ll Q
\end{align*}
for all $x\geq 1$, so by the triangle inequality and the choice of $Q$ we reduce to showing that
$$ \int_{N}^{2N}\left|\sum_{x<n\leq x+M} g(n) e(\alpha n)\right| \d x = \tilde o(MN).$$
We introduce two further scales $Q_0 \ll P_1 \ll P_2 \ll Q$ by the formulae
\begin{align*}
    P_1 &\coloneqq \log^{100} M \\
    P_2 &\coloneqq M^{1/2}.
\end{align*}
We repeat the arguments in~\cite[\S 3]{MRT}.  We recall the Ramar\'e identity
$$ 1 = \sum_{P_1 \leq p \leq P_2} \sum_{pm=n} \frac{1}{\omega_{[P_1,P_2]}(m) + 1_{p \nmid m}}$$
whenever $n$ is divisible by at least one prime between $P_1$ and $P_2$, and $\omega_{[P_1,P_2]}(m) \coloneqq \sum_{P_1 \leq p \leq p_2} 1_{p|m}$.  A standard sieve shows that for any $x$, at most $O(M\frac{\log P_1}{\log P_2})$ elements $n$ of $[x,x+M]$ are not divisible by such a prime, and hence
we can use the Ramar\'e identity to write
\begin{align*}
\sum_{x<n\leq x+M}g(n)e(\alpha n)=\sum_{P_1\leq p\leq P_2}\, \sum_{x/p<m\leq x/p+M/p}\frac{g(pm)e(\alpha pm)}{\omega_{[P_1,P_2]}(m)+1_{p \nmid m}} +O\left(M\frac{\log P_1}{\log P_2}\right).  
\end{align*}
We have $g(pm) = g(p) g(m)$ unless $p$ divides $m$; estimating this contribution by the triangle inequality and standard bounds, we conclude that
\begin{align*}
\sum_{x<n\leq x+M}g(n)e(\alpha n)=\sum_{P_1\leq p\leq P_2}\, \sum_{x/p<m\leq x/p+M/p}\frac{g(p) g(m)e(\alpha pm)}{\omega_{[P_1,P_2]}(m)+1} +O\left(M\frac{\log P_1}{\log P_2} + \frac{M}{P_1}\right).  
\end{align*}
The contribution of the error terms can easily be seen to be acceptable, so we reduce to showing that
$$ \int_{N}^{2N}\left|\sum_{P_1\leq p\leq P_2}\, \sum_{x/p<m\leq x/p+M/p}\frac{g(p) g(m)e(\alpha pm)}{\omega_{[P_1,P_2]}(m)+1}\right| \d x = \tilde o(MN).$$
We write the left-hand side as
$$ \int_{N}^{2N} \theta(x) \sum_{P_1\leq p\leq P_2}\, \sum_{x/p<m\leq x/p+M/p}\frac{g(p) g(m)e(\alpha pm)}{\omega_{[P_1,P_2]}(m)+1} \d x$$
for some $1$-bounded $\theta\colon \mathbb{R}\to \mathbb{C}$ supported on $[N,2N]$.  By partitioning the range $[P_1,P_2]$ in the sum over $p$ into intervals of the form $[P,P']$ with $3P/2\leq P'\leq 2P$, it suffices by the triangle inequality to show that
\begin{align*}
\sum_{1\leq m\leq 3N/P} \frac{g(m)}{\omega_{[P_1,P_2]}(m)+1} \sum_{P\leq p\leq P'} g(p) e(\alpha pm)\int_{\mathbb{R}}\theta(x)1_{x<mp\leq x+M}\d x = \tilde o(MN \log^{-1} P).
\end{align*}
By H\"older's inequality, the left-hand side can be bounded by 
\begin{align*}
\ll \left(\frac{N}{P}\right)^{3/4} S^{1/4},     
\end{align*}
where
\begin{align*}
S &\coloneqq \sum_{1\leq m\leq 3N/P}\left|\sum_{P\leq p\leq P'} g(p) e(\alpha pm)\int_{\mathbb{R}}\theta(x)1_{x<mp\leq x+M}\d x\right|^4 \\
&=\int_{(x_1,x_2,x_3,x_4)\in \mathbb{R}^4}\theta(x_1)\theta(x_2)\overline{\theta(x_3)}\overline{\theta(x_4)}\sum_{P\leq p_1,p_2,p_3,p_4\leq P'}g(p_1)g(p_2)\overline{g(p_3)}\overline{g(p_4)}\\
&\quad \quad\sum_{\substack{1\leq m\leq 3N/P\\ x/p_j<m\leq (x+M)/p_j\,\, \forall j\in [4]}}e(\alpha(p_1+p_2-p_3-p_4)m)\d x_1 \d x_2 \d x_3\d x_4.  
\end{align*}
By the geometric sum formula, the inner sum is $\ll \min\{M/P,\frac{1}{\|(p_1+p_2-p_3-p_4)\|}\}$. By Selberg's sieve, any integer $n=O(P)$ has $\ll P^3 \log^{-4} P$ representations as $n=p_1+p_2-p_3-p_4$. Hence we reduce to showing that
\begin{align*}
\sum_{|n|\leq 4P}\min\left\{\frac{M}{P},\frac{1}{\|\alpha n\|}\right\} = \tilde o(M).
\end{align*}
But recalling the minor arc hypothesis, this follows from~\cite[Page 346]{ik}.

\subsection{The major arc case}

Now suppose that $|\alpha-a/q|\leq \frac{1}{qQ}$ for some coprime $a,q$ with $1\leq a\leq q\leq Q_0$.  Our arguments here loosely follow those in~\cite{teravainen-binary}. We can subdivide $[x,x+M]$ (up to negligible error) into intervals of the slightly shorter length
$$ \tilde M \coloneqq M {\mathcal L}^{-1/5}$$
and also subdivide modulo $q$ to reduce to establishing the bound
$$ \int_{N}^{2N}\left|\sum_{x<n\leq x+\tilde M} (g(n)-\delta) e(\alpha n) 1_{n\equiv b\pmod q}\right| \d x \ll \tilde o\left(\frac{\tilde MN}{q}\right)$$
for all $b \in [q]$.  By hypothesis, $e(\alpha n)$ only varies by $O(\mathcal{L}^{-c/5} H)$ in the arithmetic progression $\{x < n \leq x+\tilde M: n \equiv b \pmod q\}$, so it suffices to show that
$$ \int_{N}^{2N}\left|\sum_{x<n\leq x+\tilde M} (g(n)-\delta) 1_{n\equiv b\pmod q}\right| \d x\ll \tilde o\left( \frac{\tilde MN}{q}\right).$$
Let $d=(b,q)$ and write $q=dq'$, $b=db'$, $N = dN'$, $M = dM'$. Then $(b',q')=1$ and we can rewrite the previous bound as
$$ \int_{N'}^{2N'}\left|\frac{1}{M'} \sum_{x<n\leq x+M'} (g(dn)-\delta) 1_{n\equiv b'\pmod q'}\right| \d x\ll \tilde o\left( \frac{N'}{q'}\right).$$
By Dirichlet character expansion and the triangle inequality, it suffices to show that
$$ \int_{N'}^{2N'}\left|\frac{1}{M'} \sum_{x<n\leq x+M'} g(dn) \chi(n) - \frac{\varphi(q')}{q'} \delta 1_{\chi=\chi_0}\right| \d x\ll N' {\mathcal L}^{-\eta}$$
for some absolute constant $\eta$ and for all characters $\chi$ of period $q$, with $\chi_0$ denoting the principal character.  Using the assumption~\eqref{equi} and the triangle inequality, we can rewrite this as
$$ \int_{N'}^{2N'}\left|\frac{1}{M'} \sum_{x<n\leq x+M'} g(dn) \chi(n) - \frac{1}{N'} \sum_{N' < n \leq 2N'} g(dn) \chi(n) \right| \d x\ll N' {\mathcal L}^{-\eta}.$$
We can write $g(dn) \chi(n) = g(dr) \chi(dr) g(n_1) \chi(n_1)$ where $r$ is only divisible by primes dividing $d$, and $n_1=n/r$ is coprime to $d$.  By the triangle inequality, we can bound the left-hand side by
$$ \sum_{r \mid d^\infty} \int_{N'/r}^{2N'/r}\left|\frac{r}{M'} \sum_{x_1<n_1\leq x_1+M'/r} g(n_1) \chi'(n_1) - \frac{r}{N'} \sum_{N'/r < n \leq 2N'/r} g(n_1) \chi'(n_1) \right| \d x_1$$
where $\chi'(n_1) \coloneqq \chi(n_1) 1_{(n_1,d)=1}$ and $r\mid d^{\infty}$ means $r\mid d^j$ for some natural number $j$.  The integral can be crudely bounded by $O(N'/t)$ when $t \leq 2N'$, and vanishes otherwise.  Since
$$ \sum_{r\mid d^\infty} \frac{N'}{r^{1/2}} \ll N' \exp\left(O\left(\sum_{p\mid d} \frac{1}{p^{1/2}}\right)\right) \ll N' \tau(d)^{O(1)}$$
we see from the divisor bound that the contribution of those $r$ with $r \geq {\mathcal L}^{1/2}$ (say) is acceptable.  For the remaining $r$, it will suffice to show that
$$ \int_{N'/r}^{2N'/r}\left|\frac{t}{M'} \sum_{x_1<n_1\leq x_1+M'/r} g(n_1) \chi'(n_1) - \frac{r}{N'} \sum_{N'/r < n \leq 2N'/r} g(n_1) \chi'(n_1) \right| \d x_1 \ll \frac{N'}{r} {\mathcal L}^{-2\eta}$$
(say).
If $\chi'$ is real, then since $g_1$ is also real-valued, this follows from the Matom\"aki--Radziwi\l{}\l{} theorem~\cite[Theorem 1]{mr-annals}.  If instead $\chi'$ is complex, we may use a complex-valued variant of the Matom\"aki--Radziwi\l{}\l{} theorem from~\cite[Theorem A.1]{MRT}, assuming that we can establish the non-pretentious estimate
\begin{align*}
 \mathbb{D}(g \chi',n \mapsto n^{it};N)^2 \geq \frac{1}{100}\log_2 H,   
\end{align*}
say, whenever $|t| \leq N$.  But this follows from~\cite[Lemma C.1]{MRT}\footnote{That lemma assumes that $\chi$ is fixed, but the proof clearly works for all $\chi$ of modulus at most $\log^{100} H$, say.}.

\subsection{Application to smooth numbers} 

We conclude this section with a quick deduction of Theorem~\ref{thm:smooth} from Theorem~\ref{thm:correlation}.

\begin{proof}[Proof of Theorem~\ref{thm:smooth}]
By a standard estimate on smooth numbers (see~\eqref{nsmooth}), we may assume that $u,v\in [1,\log_2 x]$.
By the union bound, if we have the claim with $\mathcal{X}$ depending on $u,v$, then we also have it (with a smaller value of $c$) for all $u,v$ that are multiples of $\log_2^{-4c} x$. But since $\rho$ is a Lipschitz function and the number of integers $n\in [1,x]$ that have a prime factor in the range $[x^{\alpha},x^{\alpha+\log^{-4c}x}]$  is  $\ll x\alpha^{-1}\log^{-4c}x$ for $1\leq \alpha\leq \log_2 x$ by Mertens' theorem, we can use the same exceptional set for all $u,v$.

Now, it suffices to show that for any given $u,v\in [1,c\log_2 x]$ we have~\eqref{eq:smoothcorr} for all $x\in \mathcal{X}$, where $\mathcal{X}$ satisfies~\eqref{eq:exceptional} for some constant $c>0$. By Theorem~\ref{thm:correlation}, we only need to check that for all $N\in [x^{0.4},x]$, $a,q\in \mathbb{N}$ and $\alpha\in [1,c\log_2 x]$ we have
\begin{align*}
\sum_{\substack{N<n\leq 2N\\ n\equiv a\pmod q}} 1_{n \textnormal{ is }x^{1/\alpha}\textnormal{ smooth}}=\frac{N}{q}\rho\left(\frac{\alpha\log N}{\log x}\right)+O(N\log^{-1/2} x),    
\end{align*}
say. But after using multiplicativity to reduce to the case $(a,q)=1$ (with the range of $N$ relaxed to $[x^{0.3},x]$) this follows from standard estimates for smooth numbers (see, e.g.,~\cite[Theorem 6.1 and (6.6)]{hildebrand-tenenbaum} or~\cite{fouvry-tenenbaum}) and the Lipschitz property of $\rho$.
\end{proof}

\section{Application to consecutive values of arithmetic functions}\label{hildebrand-sec}

In this section we prove \Cref{consecutive-omega-thm}.  We begin with the claim for $\Omega$, which already contains the main ideas. Let $X$ be a large asymptotic parameter going to infinity.  In this section, we use the notation that
\begin{align}\label{eq:tildedef3}
Y=\tilde o(Z)\,\, \textnormal{ means }\,\, Y=O(Z\log^{-c}_2 X)\,\, \textnormal{ for some absolute constant } c>0    
\end{align}
(note that this differs from our usage in Sections~\ref{consecutive-sec}--\ref{corr-sec}). It will suffice to show that
\begin{equation}\label{omega-control}
 |\{x \leq n < 2x\colon \Omega(n)=\Omega(n+1)\}|=\left(\frac{1}{2\sqrt{\pi}}+\tilde o(1)\right) x \log^{-1/2}_2 x
 \end{equation}
for all $x \in [\sqrt{X},X]$ outside of an exceptional set ${\mathcal E}$ with
\begin{equation}\label{excep}
 \int_{\mathcal E} \frac{\d x}{x} = \tilde o(\log X).
 \end{equation}
 By Markov's inequality, it suffices to show that
$$ \int_{\sqrt{X}}^X \left|\frac{1}{x}|\{x \leq n < 2x\colon \Omega(n)=\Omega(n+1)\}| - \frac{1}{2 \sqrt{\pi} \log^{1/2}_2 x} \right| \frac{\d x}{x} = \tilde o\left(\frac{\log X}{\log^{1/2}_2 X}\right).$$
In order to handle the $\tau$ case later, we will establish the following technical generalization: if $1 \leq w \leq \log_2 x$ and $d_0, d_1 \in \N_{\leq w}$ are coprime with
\begin{equation}\label{wbound}
d_0,d_1 \leq \exp(2\log^{1/4}_2 x),
\end{equation}
and $l = O(\log_3 x)$ is an integer, then
\begin{equation}\label{generalize}
 \int_{\sqrt{X}}^X \left|\frac{1}{x}|\{x \leq n < 2x\colon \Omega_{>w}(n)=\Omega_{>w}(n+1)-l;\, n_{\leq w} = d_0;\, (n+1)_{\leq w} = d_1\}| - \kappa \right| \frac{\d x}{x} = \tilde o\left(\kappa \log X\right)
 \end{equation}
 where
 $$
 \kappa \coloneqq \frac{1}{2 \sqrt{\pi} \log^{1/2}_2 X} \frac{\varphi(d_0 d_1)}{d_0^2 d_1^2} \prod_{\substack{p \leq w\\ p \nmid d_0 d_1}} \left( 1 - \frac{2}{p} \right).$$
Readers who are only interested in the $\Omega$ and $\omega$ cases and not the $\tau$ case may safely set $w = d_0 = d_1 = 1$ and $l=0$ in the arguments that follow, so that $\kappa = \frac{1}{2\sqrt{\pi} \log^{1/2}_2 X}$.

We use the circle method.  The left-hand side can be written as
$$ \int_{\sqrt{X}}^X \left| \frac{1}{x}\int_{-1/2}^{1/2} e(\alpha l) \sum_{x \leq n < 2x} e(\alpha \Omega_{>w}(n)) e(-\alpha \Omega_{>w}(n+1)) 1_{n_{\leq w} = d_0;\, (n+1)_{\leq w} = d_1}\ \d \alpha - \kappa \right| \frac{\d x}{x}.$$
We introduce the scale
$$ R \coloneqq \log^{1/4}_2 X$$
and let $\eta \colon \R \to \R$ be a smooth function supported in $[-2,2]$ that equals one on $[-1,1]$; then $\eta(R\cdot)$ is supported in $[-2/R,2/R]$.  It will suffice to establish the low frequency estimate
\begin{equation}\label{lowfreq}
\int_\R \eta(R\alpha) e(\alpha l) \sum_{x \leq n < 2x} e(\alpha \Omega_{>w}(n)) e(-\alpha \Omega_{>w}(n+1)) 1_{n_{\leq w} = d_0;\, (n+1)_{\leq w} = d_1}\d\alpha = (1+\tilde o(1)) \kappa x
\end{equation}
 for all $\sqrt{X} \leq x \leq X$, as well as the high frequency estimate
\begin{equation}\label{high-frequency}
\int_{\sqrt{X}}^X \left| \frac{1}{x} \int_{-1/2}^{1/2} (1-\eta(R\alpha)) e(\alpha l) \sum_{x \leq n < 2x} e(\alpha \Omega_{>w}(n)) e(-\alpha \Omega(n+1))1_{n_{\leq w} = d_0;\, (n+1)_{\leq w} = d_1}\ \d \alpha\right| \frac{\d x}{x} = \tilde o\left(\kappa \log X\right).
\end{equation}

\subsection{The high frequency estimate}

We begin with the high frequency estimate~\eqref{high-frequency}.  
Since $1-\eta(R\alpha)$ is supported on the region $1/R \leq |\alpha| \leq 1/2$, it suffices by the triangle inequality to show that
\begin{equation}\label{cor}
\int_{\sqrt{X}}^X \frac{1}{x} \left| \sum_{x \leq n < 2x} e(\alpha \Omega_{>w}(n)) e(-\alpha \Omega_{>w}(n+1)) 1_{n_{\leq w} = d_0;\, (n+1)_{\leq w} = d_1} \right|\frac{\d x}{x} = \tilde o\left(\kappa \log X\right).
\end{equation}
for all such $\alpha$.  

The function $n\mapsto e(\alpha \Omega_{>w}(n))$ is multiplicative\footnote{In fact it is completely multiplicative, but we will not use this fact as it is not true for the $\omega$ counterpart.} and of modulus $1$.  We will establish the lower bound
\begin{equation}\label{expy}
 M(n \mapsto e(\alpha \Omega_{>w}(n)); X^2, \log^{1/125} X) ) \gg \log_2 X.
 \end{equation}
Note that, writing $W=\prod_{p\leq w}p$ and $v_p(n)$ for the $p$-adic valuation of $n$, we can decompose
\begin{align*}
1_{n_{\leq w} = d_0;\, (n+1)_{\leq w} = d_1}=\prod_{j=0}^1 1_{n\equiv j\pmod{d_j}}1_{(n+j,W/d_j)=1}\cdot \prod_{p\mid d_0d_1}(1-1_{p^{v_p(d_0)+1}\mid n})(1-1_{p^{v_p(d_1)+1}\mid n+1}),          \end{align*}
and by the Chinese remainder theorem this can be expressed as an alternating sum of $\ll 4^{\omega(d_0d_1)}\ll d_0d_1$ terms of the form
\begin{align*}
1_{n \equiv b \pmod U}\cdot  g_0(n) g_1(n+1),    
\end{align*}
where $b\in \mathbb{N}$, $U\mid d_0^2d_1^2$, and $g_j$ is the completely multiplicative function with $g(p) = 1_{p \nmid W/d_j}$ for all primes $p$. Hence,~\eqref{expy} implies the desired bound~\eqref{cor} by the non-pretentious case of  \Cref{thm:correlation} with ${\mathcal L} = \log^{c} X$ for some absolute constant $c>0$, since we can use~\eqref{wbound} to absorb losses of $d_0^{O(1)} d_1^{O(1)}$ into the ${\mathcal L}^{-c}$ type factors (and also noting for similar reasons that ${\mathcal L}^{-c} \ll \kappa$).  

It remains to show~\eqref{expy}.  By the definition~\eqref{mdef2} of $M(\cdot)$, it suffices to show that
$$  \sum_{w < p \leq X^2} \frac{1 - \mathrm{Re}( e(\alpha) \chi(p) p^{it} )}{p} \gg \log^{1/2}_2 X$$
for all $|t| \leq X^2$ and all Dirichlet characters $\chi$ of period at most $\log^{1/125} X$. 
Suppose first that $\chi$ is non-principal, or that $\chi$ is principal with $|t| \geq 2 \log_2 X$.  In either case,
then $L(s,\chi)$ has the Vinogradov--Korobov zero-free (and pole-free) region of the form 
$$ \left\{ s\in \mathbb{C}\colon |\textnormal{Im}(s)-t| \leq \log_2 X; \textnormal{Re}(s) \geq 1 - \frac{c}{\log^{2/3} X \log_2^{1/3} X} \right\};$$
for some constant $c>0$; see, e.g.,~\cite[Theorem 2]{coleman}.
From this and a standard Perron formula argument we can show that
$$ \sum_{\exp(\log^{3/4} X) \leq p \leq X} \frac{\chi(p) p^{it}}{p} = O(1)$$
(say), and hence on discarding all primes in $[w,\exp(\log^{3/4} X)]$ we have
$$ \sum_{w<p \leq X^2} \frac{1 - \mathrm{Re}( e(\alpha) \chi(p) p^{it} )}{p} \gg \log_2 X$$
in this case from Mertens' theorem.  It remains to consider the case when $\chi$ is principal and $|t| < 2 \log_2 X$.  If $\log^{-1/2} X \leq |t| < 2 \log_2 X$, then from the prime number theorem with classical error term we can write
$$ \sum_{\exp(\log^{1/2}X)<p \leq X^2} \frac{\chi(p) p^{it}}{p} =\int_{\exp(\log^{1/2}X)}^{X^2}\frac{y^{it}}{y\log y}\d y\ll \frac{1}{|t|\log^{1/2}X}\ll 1,$$
so in this case from Mertens' theorem one has
$$  \sum_{w<p \leq X^2} \frac{1 - \mathrm{Re}( e(\alpha) \chi(p) p^{it} )}{p} \geq \frac{1}{2}\log_2 X - O(1),$$
which is acceptable. Finally, if $\chi$ is principal and $|t| < \log^{-1/2} X$, then for $\log^{1/125} X<p \leq \exp(\log^{1/4} X)$ one has
$$ 1 - \mathrm{Re}( e(\alpha) \chi(p) p^{it})=1-\textnormal{Re}(e(\alpha))+O(|t|\log^{1/4}X) \gg R^{-2},$$
and hence
$$ \sum_{w<p \leq X^2} \frac{1 - \mathrm{Re}( e(\alpha) \chi(p) p^{it} )}{p} \gg R^{-2} \log_2 X=\log^{1/2}X.$$
  Thus we obtain~\eqref{expy} in all cases.

\subsection{The low frequency estimate}

Now we prove the low frequency estimate~\eqref{lowfreq}.  We can rewrite the left-hand side as
$$ R \sum_{x \leq n < 2x} \hat \eta\left(2\pi \frac{\Omega_{>w}(n+1) - \Omega_{>w}(n) + l}{R} \right) 1_{n_{\leq w} = d_0; (n+1)_{\leq w} = d_1}$$
where 
$$ \hat \eta(u) \coloneqq \int_\R \eta(t) e^{-itu}\d t$$
is the Fourier transform of $\eta$.  We introduce the scale
$$ z \coloneqq x^{1/R^{1/2}}$$ 
and partition the primes into three classes:
\begin{itemize}
    \item \emph{tiny primes} $p \leq w$;
    \item \emph{small primes} $w < p \leq z$; and
    \item \emph{large primes} $p > z$.
\end{itemize}
A number of size $O(x)$ has at most $O(\log x/\log z) = O(R^{1/2})$ large prime factors (counting multiplicity), thus
$$ \Omega_{>z}(n) \ll R^{1/2}$$
and thus
$$ \Omega_{>w}(n) = \Omega_{w < \cdot \leq z}(n) + O(R^{1/2}).$$
Similarly for $\Omega_{>w}(n+1)$.  We can thus approximate
$$ \frac{\Omega_{>w}(n+1) - \Omega_{>w}(n) + l}{R} = \frac{\Omega_{w < \cdot \leq z}(n+1) - \Omega_{w < \cdot \leq z}(n)}{R}  + O(R^{-1/2}).$$
Since $\hat \eta$ is a Schwartz function, its derivative decays superpolynomially, so we may use the mean-value theorem to write the left-hand side of~\eqref{lowfreq} as the sum of a main term
\begin{equation}\label{hmain}
R \sum_{x \leq n < 2x} \hat \eta\left(2\pi \frac{\Omega_{w < \cdot \leq z}(n+1) - \Omega_{w < \cdot \leq z}(n)}{R}  \right) 1_{n_{\leq w} = d_0;\, (n+1)_{\leq w} = d_1}
\end{equation}
and an error term
\begin{equation}\label{herror}
O\left( R^{1/2} \sum_{x \leq n < 2x} \left(1 + \frac{\Omega_{w < \cdot \leq z}(n+1) - \Omega_{w < \cdot \leq z}(n)}{R} \right)^{-10} \right) 1_{n_{\leq w} = d_0;\, (n+1)_{\leq w} = d_1}
\end{equation}
(say). Let us first consider the main term~\eqref{hmain}.  By~\eqref{nsmooth} the proportion of $x \leq n < 2x$ with $n_{\leq z} > x^{1/R^{1/4}}$ is $O(\exp(-cR^{1/4}))$ for some absolute constant $c>0$.  Hence the contribution of such $n$ is negligible.  Similarly if $(n+1)_{\leq z} > x^{1/R^{1/4}}$.  Therefore, we can restrict our attention to natural numbers $n$ for which both $n_{\leq z}$ and $(n+1)_{\leq z}$ are bounded above by $x^{1/R^{1/4}}$, and hence the small prime factors $d'_0 \coloneqq n_{w < \cdot \leq z}$ and $d'_1 \coloneqq (n+1)_{w < \cdot \leq z}$ are also bounded by $x^{1/R^{1/4}}$.  
Thus this portion of the main term can be rewritten as
$$R  \sum_{(d'_0, d'_1) \in S'} \hat \eta\left(2\pi \frac{\Omega(d'_1)-\Omega(d'_0)}{R}  \right)
\sum_{x \leq n < 2x} 1_{n_{\leq w} = d_0;\, (n+1)_{\leq w} = d_1;\, n_{w < \cdot \leq z} = d'_0;\, (n+1)_{w < \cdot \leq z} = d'_1},$$
where $S'$ is the set of pairs of coprime $d'_0, d'_1 \in \N_{w < \cdot \leq z}$ with $d'_0,d'_1 \leq x^{1/R^{1/4}}$.  
By the fundamental lemma of sieve theory (\cite[Theorem 6.12]{
cribro} applied to the set $\mathcal{A}=\{n(n+1)/(d_0d_0'd_1d_1')\colon n\in [x,2x],\, d_0d_0'\mid n,\, d_1d_1'\mid n+1\}$), the inner sum is equal to
$$ (1 + O( \exp( - c R^{1/4} ))) \frac{1}{d_0 d_1 d'_0 d'_1}\prod_{p \mid d_0 d_1 d'_0 d'_1} \left(1-\frac{1}{p}\right)  \prod_{\substack{p \leq z\\ p \nmid d_0 d_1 d'_0 d'_1}} \left(1-\frac{2}{p}\right) x$$
for some $c>0$.  We can factor this as
$$ (2 \sqrt{\pi} + O( \exp( - c R^{1/4} ))) \kappa x(\log_2^{1/2} X) \frac{\varphi(d'_0 d'_1)}{(d'_0)^2 (d'_1)^2} \prod_{\substack{w < p \leq z\\ p \nmid d'_0 d'_1}} \left( 1 - \frac{2}{p} \right).$$
Thus, to control the main term, it suffices to show that
\begin{equation}\label{dax}
R  \sum_{(d'_0, d'_1) \in S'} \frac{\varphi(d'_0d'_1)}{(d'_0)^2 (d'_1)^2} \hat \eta\left(2\pi \frac{\Omega(d'_1)-\Omega(d'_0)}{R}  \right) \prod_{\substack{w < p \leq z\\ p \nmid d'_0 d'_1}} \left(1-\frac{2}{p}\right) = 
\left(\frac{1}{2\sqrt{\pi}}+\tilde o(1)\right) \log^{-1/2}_2 x
\end{equation}
and
\begin{equation}\label{rax}
R \exp(-cR^{1/4})  \sum_{(d'_0, d'_1) \in S'} \frac{\varphi(d'_0d'_1)}{(d'_0)^2 (d'_1)^2} \prod_{\substack{w < p \leq z\\ p \nmid d'_0 d'_1}} \left(1-\frac{2}{p}\right) \ll \log^{-1/2}_2 x
\end{equation}
The left-hand side of~\eqref{rax} may be bounded by
$$ \ll R \exp(-cR^{1/4})  \prod_{w < p \leq z} \left(1 + O\left(\frac{1}{p^2}\right) \right) $$
which is acceptable, so we now turn to~\eqref{dax}.
We can remove the constraints $d'_0,d'_1 \leq x^{1/R^{1/4}}$ with negligible error.  Indeed, the cost of removing this condition can be crudely bounded using Rankin's trick by
\begin{align*}
&\ll R x^{-1/R^{1/4}} \sum_{d'_0, d'_1 \in \N_{w < \cdot \leq z}} \frac{1}{(d'_0)^{1+1/\log z} (d'_1)^{1+1/\log z}} \prod_{\substack{w < p \leq z\\ p \nmid d'_0 d'_1}} \left(1-\frac{2}{p}\right)\\
&\ll  Rx^{-1/R^{1/4}} \prod_{w < p \leq z} \left( 1 + O\left(\frac{\log p}{p \log z} + \frac{1}{p^2} \right)\right) \\
&\ll R(\log z)^{O(1)} x^{-1/R^{1/4}}
\end{align*} 
by Mertens' theorem, which is acceptable.  Thus we now reduce to showing that
\begin{equation}\label{dax-2}
R  \sum_{(d'_0, d'_1) \in S''} \frac{\varphi(d'_0d'_1)}{(d'_0)^2 (d'_1)^2} \hat \eta\left(2\pi \frac{\Omega(d'_1)-\Omega(d'_0)}{R}  \right) \prod_{\substack{w < p \leq z\\ p \nmid d'_0 d'_1}} \left(1-\frac{2}{p}\right) = 
\left(\frac{1}{2\sqrt{\pi}}+\tilde o(1)\right) \log^{-1/2}_2 x
\end{equation}
where $S''$ is the set of pairs of coprime $d'_0, d'_1 \in \N_{w < \cdot \leq z}$.
We introduce the random variables $\mathbf{d}_0, \mathbf{d}_1 \in \N_{w < \cdot \leq z}$ by the formulae
$$ \mathbf{d}_i = \prod_{w < p \leq z} p^{\mathbf{a}_{i,p}} $$
for $i=0,1$, and the random pairs  $(\mathbf{a}_{0,p}, \mathbf{a}_{1,p})$ were defined in \Cref{ctau-def}.  Observe that $(\mathbf{d}_0, \mathbf{d}_1)$ takes values in $S''$, with
$$ \Prob( (\mathbf{d}_0, \mathbf{d}_1) = (d'_0, d'_1)) = \frac{\varphi(d'_0d'_1)}{(d'_0)^2 (d'_1)^2} \prod_{\substack{w < p \leq z\\ p \nmid d'_0 d'_1}} \left(1-\frac{2}{p}\right)$$
for all $(d'_0,d'_1) \in S''$. Indeed, from \Cref{ctau-def} one can check that
$$ \Prob((\mathbf{a}_{0,p}, \mathbf{a}_{1,p}) = (a_0,a_1)) = \frac{\varphi(p^{a_0+a_1})}{p^{2a_0+2a_1}} \left(1-\frac{2}{p} \right)^{1_{a_0=a_1=0}}$$
whenever $a_0, a_1$ are non-negative integers with at most one of $a_0,a_1$ non-zero; multiplying over all $w < p \leq z$ gives the claim.

The estimate~\eqref{dax-2} can now be expressed probabilistically as
$$ \E  R \hat \eta\left(2\pi \frac{\Omega(\mathbf{d}_1)-\Omega(\mathbf{d}_0)}{R}  \right) = \left(\frac{1}{2\sqrt{\pi}}+\tilde o(1)\right) \log^{-1/2}_2 x$$
where the random variables $\mathbf{d}_0, \mathbf{d}_1 \in \N_{w < \cdot \leq z}$ take the form
$$ \mathbf{d}_i = \prod_{w < p \leq z} p^{\mathbf{a}_{i,p}} $$
for $i=0,1$, and the random pairs  $(\mathbf{a}_{0,p}, \mathbf{a}_{1,p})$ were defined in \Cref{ctau-def}.  If we set 
$$\mathbf{b}_p \coloneqq \mathbf{a}_{1,p} - \mathbf{a}_{0,p},$$
then we have
$$ \Omega(\mathbf{d}_1) - \Omega(\mathbf{d}_0) = \sum_{w < p \leq z} \mathbf{b}_p$$
and $\mathbf{b}_p$ for small primes $p$ are independent random variables taking values in $\Z$ with distribution
\begin{align*}
    \Prob( \mathbf{b}_p = 0 ) &= 1 - \frac{2}{p} \\
    \Prob( \mathbf{b}_p = \pm j ) &= \left(1-\frac{1}{p}\right)\frac{1}{p^j}
\end{align*}
for $j \in \N$.  Our task is now to show that
\begin{equation}\label{pbp-targ}
\E  R \hat \eta\left(2\pi \frac{\sum_{w < p \leq z} \mathbf{b}_p}{R}  \right) = \frac{1+\tilde o(1)}{2\sqrt{\pi \log_2 x}}
\end{equation}
Each $\mathbf{b}_p$ has mean zero and variance $\frac{2}{p} + O(\frac{1}{p^2})$, hence by Mertens' theorem and independence, the sum $\sum_{w < p \leq z} \mathbf{b}_p$ has mean zero and variance $(2+\tilde o(1))\log_2 x$.  

\subsection{A local limit theorem}

At this point we invoke a local limit theorem.

\begin{lemma}[Local limit theorem]\label{llt} We have
\begin{equation}\label{pbp}
 \Prob\left( \sum_{w < p \leq z} \mathbf{b}_p = m \right) = \frac{e^{-m^2 /(4 \log_2 x)} + \tilde o(1)}{2\sqrt{\pi \log_2 x}} 
 \end{equation}
uniformly for integers $m$.  
\end{lemma}

\begin{proof}
    If the $\mathbf{b}_p$ were uniformly bounded, and we were satisfied with a qualitative $o(1)$ decay term, this would follow from Prokhorov's local limit theorem (see, e.g.,~\cite[Theorem 1.25]{szewczak}).  Because our variables are unbounded and we have a more quantitative error term, we perform a more explicit computation here.  By the circle method and independence, we can expand the left-hand side of~\eqref{pbp} as
$$ \int_{-1/2}^{1/2} e(-\alpha m) \prod_{w<p \leq z} \E e(\alpha \mathbf{b}_p)\d\alpha,$$
which after a change of variables becomes
$$ \frac{1}{\log_2^{1/2} x} \int_{-\log_2^{1/2} x/2}^{\log_2^{1/2} x/2} e(-\alpha m / \log_2^{1/2} x ) \prod_{w < p \leq z} \E e(\alpha \mathbf{b}_p / \log_2^{1/2} x )\d\alpha.$$
If we can show the approximation
\begin{equation}\label{al}
\prod_{w < p \leq z} \E e(\alpha \mathbf{b}_p / \log_2^{1/2} x ) = 
\exp( - 4\pi^2 \alpha^2 ) + \tilde o( \exp( - c \alpha^2 ) )
\end{equation}
for $|\alpha| \leq \log_2^{1/4} x$ and
\begin{equation}\label{bl}
\prod_{w < p \leq z} \E e(\alpha \mathbf{b}_p / \log_2^{1/2} x ) \ll \exp( - c \log_2^{1/2} x )
\end{equation}
for $\log_2^{1/4} x < |\alpha| \leq \log_2^{1/2} x/2$ and some $c>0$, then the claim follows from the standard Gaussian integral identity 
\begin{align*}
\int_{\mathbb{R}}\exp(-4\pi^2\alpha^2-2\pi iy\alpha)\d \alpha=\frac{1}{2\sqrt{\pi}}\exp\left(-\frac{y^2}{4}\right)    
\end{align*}
for  $y\in \mathbb{R}$ and the triangle inequality.

First suppose that $|\alpha| \leq \log_2^{1/4} x$.  By the explicit form of the distribution of $\mathbf{b}_p$, one can write
\begin{align*}
    \E e(\alpha \mathbf{b}_p / \log_2^{1/2} x ) &= 1 - \E (1 - \cos(2\pi \alpha \mathbf{b}_p / \log_2^{1/2} x )) \\
    &= 1 - 2\sum_{j=1}^\infty (1-\cos(2\pi \alpha j / \log_2^{1/2} x))\left(1-\frac{1}{p}\right)\frac{1}{p^j} \\
    &= 1 - \frac{4\pi^2\alpha^2}{p \log_2 x} + O\left(\frac{\alpha^2}{p^2 \log_2 x} + \frac{\alpha^4}{p \log_2^2 x} \right) \\
    &= \exp\left( - \frac{4\pi^2\alpha^2}{p \log_2 x} + O\left(\frac{\alpha^2}{p^2 \log_2 x} + \frac{\alpha^4}{p \log_2^2 x} \right) \right)\\
    &= \exp\left( - \frac{4\pi^2\alpha^2}{p \log_2 x} + O\left(\frac{\alpha^2}{p^2 \log_2 x} + \frac{\alpha^2}{p \log_2^{3/2} x} \right) \right). 
\end{align*}
Taking products and using Mertens' theorem in the form 
\begin{equation}\label{mert2}
\sum_{w<p\leq z}\frac{1}{p}=(1+\tilde{o}(1))\log_2 x,
\end{equation}
we can bound the left-hand side of~\eqref{al} by
$$ \exp( -  4\pi^2 \alpha^2 + \tilde o(\alpha^2) )$$
and the claim~\eqref{al} follows from the mean value theorem.  For $\log_2^{1/4} x < |\alpha| \leq \log_2^{1/2} x/2$, we use the monotone decreasing nature of cosine on $[0,\pi]$ and the estimates $1-\cos t\asymp t^2$ and $e^t=1+O(t)$ for $|t|\leq 1$ to instead bound
\begin{align*}
    \E e(\alpha \mathbf{b}_p / \log_2^{1/2} x ) &= 1 - 2\sum_{j=1}^\infty (1-\cos(2\pi \alpha j / \log_2^{1/2} x))\left(1-\frac{1}{p}\right)\frac{1}{p^j} \\
    &\leq 1 - 2(1-\cos(2\pi \alpha / \log_2^{1/2} x)) \left(1-\frac{1}{p}\right)\frac{1}{p}+ O\left(\frac{1}{p^2} \right) \\
    &\leq 1 - 2(1-\cos(2\pi / \log_2^{1/4} x)) \left(1-\frac{1}{p}\right)\frac{1}{p}  + O\left(\frac{1}{p^2} \right) \\
    &\leq \exp\left( - \frac{c}{p \log^{1/2}_2 x} + O\left(\frac{1}{p^2} \right) \right) 
\end{align*}
for some $c>0$; multiplying in $p$ and using~\eqref{mert2} again, we obtain~\eqref{bl} in this case.
\end{proof}

\subsection{Conclusion of the argument}

We apply Lemma~\ref{llt} to~\eqref{pbp-targ}. The contribution of the $\tilde o(1)$ error can be seen to be negligible due to the decay of $\hat \eta$, so it suffices to show that
$$ \sum_{m=-\infty}^\infty R \hat \eta\left(2\pi\frac{m}{R}\right) e^{-m^2 / 4 \log_2 x} = 1 + \tilde o(1).$$
But this is clear since  by the Poisson summation formula we have $$\sum_{m=-\infty}^{\infty}R \hat \eta\left(2\pi\frac{m}{R}\right)=\sum_{k=-\infty}^{\infty}\eta(Rk)=1$$
and  $e^{-m^2 / 4 \log_2 x}=1 + \tilde o(1)$ for $|m| \leq R^{3/2}$ (say), and $R \hat \eta(2\pi\frac{m}{R})\ll R^{-10}$ for $|m|> R^{3/2}$.  This concludes the proof of~\eqref{hmain}.  A similar analysis (but without the need for such fine control on main terms) lets one bound the left-hand side of~\eqref{herror} by
$$ \ll R^{1/2} \sum_{m=-\infty}^\infty \left(1 + \frac{m}{R} \right)^{-10} \frac{e^{-m^2 / 4 \log_2 x} +\tilde o(1)}{2\sqrt{\pi \log_2 x}}$$
which one can easily verify to be $O(R^{-1/2}) = \tilde o(1)$.  This concludes the proof of \Cref{consecutive-omega-thm} for $\Omega$.

For $\omega$, the argument is similar, but the random variable $\mathbf{b}_p$ now has to be defined as
$$\mathbf{b}_p \coloneqq \min(\mathbf{a}_{1,p},1) - \min(\mathbf{a}_{0,p},1),$$
which now takes values in $\{-1,0,1\}$ with distribution
\begin{align*}
    \Prob( \mathbf{b}_p = 0 ) &= 1 - \frac{2}{p} \\
    \Prob( \mathbf{b}_p = \pm 1 ) &= \frac{1}{p}
\end{align*}
However, the mean and variance of $\sum_{p \leq z} \mathbf{b}_p$ obey the same asymptotics as before.  The analogue of~\eqref{llt} can then be achieved by similar calculations.

\subsection{The analogue for the divisor function}

Now we turn to the analogous problem for $\tau$, which is a little trickier.  With $X$ as before, it suffices by Markov's inequality to show that
\begin{equation}\label{nx}
\int_{\sqrt{X}}^X \left| \frac{1}{x} |\{x \leq n < 2x\colon \tau(n)=\tau(n+1)\}| - \frac{c_\tau}{2\sqrt{\pi \log_2 X}} \right| \frac{\d x}{x}
= \tilde o\left(\frac{\log X}{\log^{1/2}_2 X} \right).
 \end{equation}

We introduce a cutoff
$$ w \coloneqq \log_2 X.$$
The contribution of those $n$ which are divisible by $p^2$ for some $p > w$ to the integrand is $O(1/w)$ which is acceptable, and similarly for $n+1$.  Thus we may restrict attention to those $n$ for which $n_{>w}$ and $(n+1)_{>w}$ are both squarefree.  In that case we have
$$ \tau(n) = \tau(n_{\leq w}) 2^{\Omega(n_{>w})}$$
and similarly for $n+1$.  Next, by~\eqref{nsmooth} the contribution to~\eqref{nx} from those $n$ for which $n_{\leq w} > \exp(\log^{1/4}_2 X)$ is by Rankin's trick
\begin{align*}
    &\ll \sum_{\substack{d \in \N_{\leq w}\\ d > \exp(\log^{1/4}_2 X)}} \int_{\sqrt{X}}^X \frac{1}{d} \frac{\d x}{x} \\
    &\ll (\log X)(\log w) \exp( - (\log^{1/4}_2 X)/(\log w)) \\
    &= \tilde o\left(\frac{\log X}{\log^{1/2}_2 X} \right)
\end{align*}
for some constant $c>0$, thus giving a negligible contribution; similarly for $(n+1)_{\leq w}$.  Finally,
we perform the standard first moment calculation
$$ \frac{1}{x} \sum_{x \leq n < 2x} \tau(n_{\leq w}) \ll \sum_{d \in \N_{\leq w}} \frac{1}{d} \ll \exp\left( \sum_{p \leq w} \frac{1}{p} \right) \ll \log w = \log_3 X,$$
 so up to acceptable errors, we can restrict to the case where $\tau(n_{\leq w}) \leq \log_2 X$, and similarly for $\tau((n+1)_{\leq w})$.  Thus we can replace the condition $\tau(n+1)=\tau(n)$ by the three conditions
\begin{itemize}
\item[(i)] $n_{\leq w}, (n+1)_{\leq w} \leq \exp(\log^{1/4}_2 X)$
\item[(ii)] $\tau(n_{\leq w}), \tau((n+1)_{\leq w}) \leq \log_2 X$ and $\tau((n+1)_{\leq w}) = 2^l \tau(n_{\leq w})$ for some integer $l$ (which is then necessarily in $[-2\log_3 X, 2\log_3 X]$).
\item[(iii)] $\Omega((n+1)_{>w}) = \Omega(n_{>w}) - l$.
\end{itemize}
Also we must have $n_{\leq w}, (n+1)_{\leq w}$ coprime.
We can thus replace the left-hand side of~\eqref{nx} by
$$ \int_{\sqrt{X}}^X \left| \sum_{|l|\leq 2\log_3 X} \sum_{(d_0,d_1) \in S_l}
 \frac{1}{x} \sum_{x \leq n < 2x} 1_{n_{\leq w} = d_0;\, (n+1)_{\leq w} = d_1} 1_{\Omega_{>w}(n+1) = \Omega_{>w}(n) - l} - \frac{c_\tau}{2\sqrt{\pi \log_2 X}}\right| \frac{\d x}{x},$$
where $S_l$ is the set of coprime pairs $d_0,d_1 \in \N_{\leq w}$ such that $d_0,d_1 \leq \exp(\log^{1/4}_2 X)$, $\tau(d_0), \tau(d_1) \leq \log_2 X$, and $\tau(d_1) = 2^l \tau(d_0)$.  Since any factor of
$$ \sum_{|l|\leq 2\log_3 X} \sum_{(d_0,d_1) \in S_l} \frac{1}{d_0 d_1} \ll \log_3 X \prod_{p \leq w} \left(1 + \frac{2}{p} + O\left(\frac{1}{p^2}\right)\right) \ll \log_3^3 X$$
can be absorbed into the $\tilde o(\cdot)$ type error terms, we may safely invoke~\eqref{generalize} and replace this expression by
$$ \int_{\sqrt{X}}^X \left| \sum_{|l|\leq 2\log_3 X} \sum_{(d_0,d_1) \in S_l} \frac{\varphi(d_0 d_1)}{d_0^2 d_1^2} \frac{1}{2\sqrt{\pi \log_2 X}}
\prod_{\substack{p \leq w\\ p \nmid d_0 d_1}} \left(1-\frac{2}{p}\right)
 - \frac{c_\tau}{2\sqrt{\pi \log_2 X}}\right| \frac{\d x}{x}.$$
Thus it will suffice to show that
\begin{equation}\label{ltau}
  \sum_{|l|\leq 2\log_3 X} \sum_{(d_0,d_1) \in S_l} \frac{\varphi(d_0 d_1)}{d_0^2 d_1^2} 
\prod_{\substack{p \leq w\\ p \nmid d_0 d_1}} \left(1-\frac{2}{p}\right) = c_\tau + \tilde o(1).
\end{equation}
We may replace the double sum with a single sum $\sum_{(d_0,d_1) \in S}$, where $S$ is the set of coprime pairs $d_0,d_1 \in \N_{\leq w}$ such that $d_0,d_1 \leq \exp(\log^{1/4}_2 X)$, $\tau(d_0), \tau(d_1) \leq \log_2 X$, and $\frac{\log \tau(d_1)-\log \tau(d_0)}{\log 2}$ is an integer.  The left-hand side can then be expressed probabilistically as
$$ \Prob( (\mathbf{d}_0,\mathbf{d}_1) \in S )$$
where the random variables $\mathbf{d}_0, \mathbf{d}_1$ are constructed by the formula
$$ \mathbf{d}_i \coloneqq \prod_{p \leq w} p^{\mathbf{a}_{i,p}} $$
for $i=0,1$, where the random pairs $\mathbf{a}_{0,p}, \mathbf{a}_{1,p}$ are defined as in \Cref{ctau-def}.

An easy application of the first moment method shows that $\log \mathbf{d}_0 \leq \log^{1/4}_2 X$ and $\tau(\mathbf{d}_0) \leq \log_2 X$ with probability $1-\tilde o(1)$, and similarly for $\mathbf{d}_1$.  Thus the above probability is equal to
$$\Prob \left( \frac{\log \tau(\mathbf{d}_1)-\log \tau(\mathbf{d}_0)}{\log 2} \equiv 0 \pmod 1 \right) - \tilde o(1).$$
We can write
$$ \frac{\log \tau(\mathbf{d}_1)-\log \tau(\mathbf{d}_0)}{\log 2} \equiv \sum_{p \leq w} \mathbf{c}_p \pmod 1 $$
where the $\mathbf{c}_p$ are defined in~\eqref{eq:cpdef}.  By~\eqref{cp0} and the union bound we have
$$\sum_{p \leq w} \mathbf{c}_p = \sum_{p} \mathbf{c}_p$$
with probability $1 - O(1/w) = 1 - \tilde o(1)$.  The claim follows.

\begin{remarks} \begin{itemize}
    \item The above analysis in fact gives a ``local limit theorem'' 
$$ \frac{1}{x} |\{n\leq x\colon \omega(n+1)-\omega(n)=m\}|=\frac{e^{-m^2 / 4 \log_2 x}+O(\log^{-c}_2 x)}{2\sqrt{\pi\log_2 x}}$$
uniformly in $m$ for all $x$ outside of an exceptional set~\eqref{excep}, and similarly for $\Omega$.  Similarly, one has
$$ \frac{1}{x} \left|\left\{n\leq x\colon \frac{\tau(n+1)}{\tau(n)}=2^m \frac{a}{b}\right\}\right|= \frac{c_{\tau,a/b} e^{-m^2 / 4 \log_2 x}+O(\log^{-c}_2 x)}{2\sqrt{\pi\log_2 x}}$$
uniformly for integers $m$ for a fixed rational $a/b > 0$ with odd numerator and denominator, where
$$ c_{\tau,a/b} \coloneqq \Prob\left( \sum_p \mathbf{c}_p = \frac{\log a - \log b}{\log 2} \pmod 1 \right)$$
is a generalization of $c_\tau = c_{\tau,1}$: we leave these generalizations of \Cref{consecutive-omega-thm} to the interested reader.  Among other things, this can be used to recover the recent result of 
 Eberhard~\cite{eberhard} that in fact every element of $(0,\infty)$ is a limit point of the sequence $(\tau(n+1)/\tau(n))_{n=1}^\infty$, which solved a problem of Erd\H{o}s~\cite{erd-boca},~\cite[Problem \#864]{bloom}.

\item With more effort, we expect it to be possible to use these methods to obtain a two-dimensional local limit theorem of the form
$$ \frac{1}{x} |\{n\leq x\colon \omega(n)=m_0, \omega(n+1)=m_1\}|=\frac{e^{-((m_0-\log_2x)^2+(m_1-\log_2 x)^2) / (4 \log_2 x)}+O(\log^{-c}_2 x)}{2\pi\log_2 x}$$
(outside of an exceptional set of $x$, of course) and similarly for $\Omega$ and $\tau$, which would in particular provide an alternative proof of the results in~\cite{charamaras} (and improve the upper bound result in~\cite[Th\'eor\`eme 3]{goudout}). However, this would require several technical modifications to the argument, in particular replacing the single frequency $\alpha$ with a pair $\alpha_0, \alpha_1$, and also establishing a two-dimensional analogue of \Cref{llt}, and we will not pursue these modifications here.

\item The results for $\omega(n) = \omega(n+1)$ can be adapted to other pairwise equations such as $\omega(n)=\omega(n+h)$ for fixed $h\neq 0$ without significant modifications.  However, to handle three-point equations such as $\omega(n)=\omega(n+1)=\omega(n+2)$ one would require, in addition to the previously mentioned extension of the methods to two dimensions, a version of \Cref{thm:correlation} for triple correlations, which does not appear to be within current\footnote{Such methods may however be able to control partially relaxed versions of such equations, such as $\omega(n)=\omega(n+1)=\omega(n+2)+O(1)$, in which the contribution of large prime factors of $n+2$ can be ignored.} technology.  For similar reasons we are currently unable to remove the exceptional set in  \Cref{consecutive-omega-thm}.
\end{itemize}
\end{remarks}

\section{Application to an irrationality problem of Erd\H{o}s}\label{irrational-sec}

In this section we establish \Cref{thm:irrational}.  As with~\cite{pratt}, our arguments will proceed by understanding the distribution of consecutive blocks $\omega(n+1),\dots,\omega(n+L)$ of the prime counting function $\omega$ for some relatively small $L$; but without the prime tuples conjecture, the control on this distribution will be quite crude.

Suppose for contradiction that the above sum is equal to $a/q$ for some integers $a, q \geq 1$, thus
$$ q \sum_{h=1}^{\infty}\frac{\omega(h)}{2^h} = 0 \pmod 1.$$
We view $q$ as fixed, and allow all implicit constants to depend on $q$.

\subsection{Shifting and dilating}

Reindexing $h$ as $n+h$ for some natural number $n$ and then multiplying by $2^n$, we conclude that
\begin{equation}\label{q1}
 q \sum_{h=1}^{\infty}\frac{\omega(n+h)}{2^h} = 0 \pmod 1
\end{equation}
for any $n \in \mathbb{N}$.  

To make use of these congruences, we exploit the additivity of $\omega$ to dilate the shift $h$.  If $p$ is a prime, we have the identity
$$ \omega(pn) = \omega(n) + 1 - 1_{p\mid n}$$
for any natural number $n$, and hence if $n$ is divisible by $p$,
$$ \omega(n+ph) = \omega\left(\frac{n}{p}+h\right) + 1 - 1_{p^2\mid n+ph}.$$
Since
$$ q \sum_{h=1}^\infty \frac{1}{2^h} = q = 0 \pmod 1$$
we conclude (after applying~\eqref{q1} with $n$ replaced with $n/p$) the identity
\begin{equation}\label{q2}
 q \sum_{h=1}^{\infty}\frac{\omega(n+ph)}{2^h} = -q\delta_p(n) \pmod 1
\end{equation}
whenever $p\mid n$, where $\delta_p(n)$ is the quantity
\begin{equation}\label{deltap-def}
 \delta_p(n) \coloneqq \sum_{h=1}^\infty \frac{1_{p^2\mid n+ph}}{2^h}.
\end{equation}
In practice, $\delta_p(n)$ will be negligible for typical $n$ and primes $p$ that are not too small.

\subsection{Taking an alternating sum to cancel terms}

Inspired by the theory of the Gowers uniformity norms, we now manipulate~\eqref{q2} to eliminate the role of the first few terms in $h$, at the cost of increasing the number of remaining terms.  Let $K$ be a natural number to be chosen later, and suppose that we can find parameters $p_0 \in \mathbb{N}$ and $v_1,\dots,v_K \in \Z$ such that the quantities
\begin{equation}\label{peps}
 p_\epsilon \coloneqq p_0 + \epsilon_1 v_1 + \dots + \epsilon_K v_K
 \end{equation}
are distinct primes for all $\epsilon = (\epsilon_1,\dots,\epsilon_K) \in \{0,1\}^K$.  Let $n$ be a natural number obeying the congruences
\begin{equation}\label{pcong}
 n \equiv  \sum_{k=1}^K k \epsilon_k v_k\pmod{p_{\epsilon}}\,\, \textnormal{ for all }\,\, \epsilon = (\epsilon_1,\dots,\epsilon_K) \in \{0,1\}^K;
 \end{equation}
we also assume the lower bound
\begin{equation}\label{nlower}
n > \sum_{k=1}^K k |v_k|
\end{equation}
so that the quantities $n - \sum_{k=1}^K k \epsilon_k v_k$ are all positive.  From~\eqref{q2} with $p=p_{\epsilon}$ and $n-\sum_{k=1}^K k \epsilon_k v_k$ in place of $n$, we have
$$ q \sum_{h=1}^{\infty}\frac{\omega(n+r_{\epsilon,h})}{2^h} = -q \delta_{p_\epsilon}\left(n - \sum_{k=1}^K k \epsilon_k v_k\right) \pmod 1$$
for each $\epsilon \in \{0,1\}^K$,
where the integer shifts $r_{\epsilon,h}$ are defined by the formula
\begin{align}\label{eq:rhk}
r_{\epsilon,h} \coloneqq p_\epsilon h - \sum_{k=1}^K k \epsilon_k v_k  = p_0h + \sum_{k=1}^K (h-k) \epsilon_k v_k.
\end{align}
Taking an alternating sum over $\epsilon$, we conclude that
$$ q \sum_{h=1}^{\infty}\frac{\sum_{\epsilon \in \{0,1\}^K} (-1)^{|\epsilon|} \omega(n+r_{\epsilon,h})}{2^h} =-q \sum_{\epsilon \in \{0,1\}^K} (-1)^{|\epsilon|} \delta_{p_\epsilon}\left(n - \sum_{k=1}^K k \epsilon_k v_k\right) \pmod 1$$
where $|\epsilon| \coloneqq \epsilon_1 + \dots + \epsilon_K$.  The point of doing this is that for $1 \leq h \leq K$, the shift $r_{\epsilon,h}$ is independent of $\epsilon_h$ thanks to~\eqref{eq:rhk}.  As a consequence, we have 
$$\sum_{\epsilon_h \in\{0,+1\}} (-1)^{\epsilon_h}\omega(n+r_{\epsilon,h})=0,$$
so the summand on the left-hand side vanishes for such $h$.  Eliminating these terms and reindexing, we obtain the identity
\begin{equation}\label{qho}
q \sum_{h=1}^{\infty}\frac{\sum_{\epsilon \in \{0,1\}^K} (-1)^{|\epsilon|} \omega(n+r_{\epsilon,h+K})}{2^{h+K}} = -q\sum_{\epsilon \in \{0,1\}^K} (-1)^{|\epsilon|} \delta_{p_\epsilon}\left(n - \sum_{k=1}^K k \epsilon_k v_k\right) \pmod 1.
\end{equation}
Applying the Lipschitz function $t\mapsto e(t)$, we conclude that
$$
e\left(q \sum_{h=1}^{\infty}\frac{\sum_{\epsilon \in \{0,1\}^K} (-1)^{|\epsilon|} \omega(n+r_{\epsilon,h+K})}{2^{h+K}}\right) = 1 + O\left(\sum_{\epsilon \in \{0,1\}^K} \delta_{p_\epsilon}\left(n - \sum_{k=1}^K k \epsilon_k v_k\right)\right).
$$

\subsection{Taking expectations and truncating}

If we let $\mathbf{n}$ be a random natural number (whose distribution is to be determined) obeying~\eqref{pcong},~\eqref{nlower}, and take expectations, we conclude
\begin{equation}\label{expectation}
\E e\left(q \sum_{h=1}^{\infty}\frac{\sum_{\epsilon \in \{0,1\}^K} (-1)^{|\epsilon|} \omega( \mathbf{n}+r_{\epsilon,h+K})}{2^{h+K}}\right) = 1 + O(\kappa_1)
\end{equation}
where $\kappa_1$ is the error term
\begin{equation}\label{kappa1-def}
\kappa_1 \coloneqq \sum_{\epsilon \in \{0,1\}^K} \E \delta_{p_\epsilon}\left(\mathbf{n} - \sum_{k=1}^K k \epsilon_k v_k\right)
= \sum_{\epsilon \in \{0,1\}^K} \sum_{h=1}^\infty \frac{\Prob( p_\epsilon^2 \mid  \mathbf{n} + r_{\epsilon,h})}{2^h}.
\end{equation}
In practice, it will be easy to show that $\kappa_1$ is negligible.  We split $\sum_{h=1}^{\infty}$ as $\sum_{h \leq H} + \sum_{h>H}$ for some cutoff $H$ to be chosen later, and use Taylor expansion to rewrite~\eqref{expectation} as
$$
\E e\left(q \sum_{1 \leq h \leq H} \frac{\sum_{\epsilon \in \{0,1\}^K} (-1)^{|\epsilon|} \omega(\mathbf{n}+r_{\epsilon,h+K})}{2^{h+K}}\right) = 1 + O(\kappa_1) + O(\kappa_2)$$
where
\begin{equation}\label{kappa2-def}
\kappa_2 \coloneqq \sum_{h > H} \frac{\sum_{\epsilon \in \{0,1\}^K} \E \omega(\mathbf{n}+r_{\epsilon,h+K})}{2^{h+K}}.
\end{equation}
We will be able to make the tail error $\kappa_2$ negligible by setting the parameter $H$ to be suitably large.

Using~\eqref{prime-count}, we can rewrite this as
\begin{equation}\label{ekk}
 \E e\left( \sum_p \mathbf{X}_p\right) = 1 + O(\kappa_1) + O(\kappa_2)
 \end{equation}
where for each prime $p$, $\mathbf{X}_p$ denotes the random variable
\begin{equation}\label{Xp-def}
\mathbf{X}_p \coloneqq q \sum_{1 \leq h \leq H} \sum_{\epsilon \in \{0,1\}^K} (-1)^{|\epsilon|} \frac{1_{p\mid \mathbf{n}+r_{\epsilon,h+K}}}{2^{h+K}}.
\end{equation}
In practice, each of the $\mathbf{X}_p$ will be individually small, but a certain sum $\sum_{p \in S_1} \mathbf{X}_p$ will be moderately large (in variance), which we will be able to use to ultimately contradict~\eqref{ekk}.

Suppose that we have some partition $S_0 \uplus S_1 \uplus S_2$ of the set of primes $p$ into a core set $S_1$ and exceptional sets $S_0$, $S_2$, with the hypothesis that the $\mathbf{X}_p$ are deterministic constants $X_p$ for $p \in S_0$ (in practice, this will be accomplished by enforcing a congruence condition on $\mathbf{n}$ for each prime in $S_0$).  By Taylor expansion and Cauchy--Schwarz one has
\begin{align*}
    \E e\left( \sum_p \mathbf{X}_p \right) &= e\left( \sum_{p \in S_0} X_p \right) \E e\left( \sum_{p \in S_1} \mathbf{X}_p \right) + O\left( \E \left|\sum_{p \in S_2} \mathbf{X}_p\right| \right) \\
    &= e\left( \sum_{p \in S_0} X_p \right) \E e\left( \sum_{p \in S_1} \mathbf{X}_p \right) + O(\kappa_3)
\end{align*}
where
\begin{equation}\label{kappa3-def}
\kappa_3 \coloneqq \left( \E \left|\sum_{p \in S_2} \mathbf{X}_p\right|^2\right)^{1/2}.
\end{equation}
So now we have
\begin{equation}\label{epsi}
 \left|\E e\left( \sum_{p \in S_1} \mathbf{X}_p \right)\right| = 1 + O(\kappa_1) + O(\kappa_2) + O(\kappa_3).
 \end{equation}
Obtaining an adequate control on $\kappa_3$ will be somewhat delicate, particularly for very large primes where the correlation estimate of Theorem~\ref{thm:correlation} is needed.

\subsection{Extracting a variance bound}

Intuitively, the estimate~\eqref{epsi} suggests that $\sum_{p \in S_1} \mathbf{X}_p$ has small variance, but to make this precise we will need to assume some approximate independence properties between the $\mathbf{X}_p$.
We let $\mathbf{X}'_p$ be a copy of $\mathbf{X}_p$, but with the $\mathbf{X}'_p$ jointly independent in $p$.  We consider the difference
$$ \E e\left( \sum_{p \in S_1} \mathbf{X}_p \right) - \E e\left( \sum_{p \in S_1} \mathbf{X}'_p \right).$$
By Taylor expansion with remainder to some even index $M$, we can write this as
\begin{align*}
    &\sum_{m=0}^{M-1} \frac{(2\pi i)^m}{m!} \E \left( \left(\sum_{p \in S_1} \mathbf{X}_p\right)^m - \E \left(\sum_{p \in S_1} \mathbf{X}'_p\right)^m\right) \\
& \quad + O\left( \frac{(2\pi)^M}{M!} \left(\E \left(\sum_{p \in S_1} \mathbf{X}_p\right)^M + \E \left(\sum_{p \in S_1} \mathbf{X}'_p\right)^M\right) \right)
\end{align*} 
and hence
$$ \E e\left( \sum_{p \in S_1} \mathbf{X}_p \right) - \E e\left( \sum_{p \in S_1} \mathbf{X}'_p \right) = O(\kappa_4) + O\left(\frac{(2\pi)^M}{M!} \E \left(\sum_{p \in S_1} \mathbf{X}'_p\right)^M\right)$$
where
\begin{equation}\label{kappa4-def}
\kappa_4 \coloneqq \sum_{m=0}^M \frac{(2\pi)^m}{m!} \left|\E \left(\sum_{p \in S_1} \mathbf{X}_p\right)^m - \E \left(\sum_{p \in S_1} \mathbf{X}'_p\right)^m\right|
\end{equation}
We will choose the primes in $S_1$ to be small enough that $\kappa_4$ is easy to make negligible.   We have the power series inequality
$$ e^t + e^{-t} \geq e^{|t|}=\sum_{j=0}^{\infty}\frac{|t|^j}{j!}\geq t^M/M!$$
for any real $t$.  Applying this with $t = 2\pi e \sum_{p \in S_1} \mathbf{X}_p$, and then taking expectations, using independence and rearranging, we conclude that
\begin{align*} 
\frac{(2\pi)^M}{M!} \E \left(\sum_{p \in S_1} \mathbf{X}'_p\right)^M &\leq \E \left(\exp\left( 2\pi e \sum_{p \in S_1} \mathbf{X}'_p\right) + \exp\left( -2\pi e \sum_{p \in S_1} \mathbf{X}'_p\right)\right) \bigg/ e^M \\
&= \frac{\prod_{p \in S_1} \E \exp(2\pi e \mathbf{X}_p) + \prod_{p \in S_1} \E \exp(-2\pi e \mathbf{X}_p)}{e^M}.
\end{align*}
As $M$ is even, the left-hand side is non-negative, hence by~\eqref{kappa4-def} and the triangle inequality we have
$$
\frac{(2\pi)^M}{M!} \left(\E \left(\sum_{p \in S_1} \mathbf{X}_p\right)^M + \E \left(\sum_{p \in S_1} \mathbf{X}'_p\right)^M\right) \ll \kappa_4 + \kappa_5$$
where
\begin{equation}\label{kappa5-def}
\kappa_5 \coloneqq \frac{\prod_{p \in S_1} \E \exp(2\pi e \mathbf{X}_p) + \prod_{p \in S_1} \E \exp(-2\pi e \mathbf{X}_p)}{e^M}.
\end{equation}
We will be able to make $\kappa_5$ small by choosing the degree $M$ to be somewhat large.

By independence,~\eqref{epsi}, and the triangle inequality, we conclude that
\begin{align*}
\left|\prod_{p \in S_1} \E e(\mathbf{X}_p)\right| &= \left| \E e\left( \sum_{p \in S_1} \mathbf{X}'_p \right)\right|\\
&= \left| \E e\left( \sum_{p \in S_1} \mathbf{X}_p \right) \right| + O(\kappa_4+\kappa_5) \\
&= 1 + \sum_{j=1}^5 O(\kappa_j).   
\end{align*} 
Assuming that 
\begin{equation}\label{kappas}
\kappa_1,\dots,\kappa_5 = o(1)
\end{equation}
as some asymptotic parameter $x$ goes to infinity along some sequence, we can then take logarithms to conclude (for sufficiently large $x$) that
$$ \sum_{p \in S_1} \log \frac{1}{|\E e(\mathbf{X}_p)|} \ll \sum_{j=1}^5 \kappa_j.$$
If we have 
\begin{equation}\label{xp}
\mathbf{X}_p = o(1)
\end{equation}
for all $p \in S_1$, then by Taylor expansion of the complex exponential there exists some $\theta\in [-1/2,1/2]$ for which we can have
$$ |\E e(\mathbf{X}_p)| = \mathrm{Re} \E e(\mathbf{X}_p - \theta)=1 - \E (1-\cos(2\pi(\mathbf{X}_p-\theta))),$$
 and then using $1-\cos(t) \asymp t^2$ for $|t|\leq (2-1/10)\pi$ and $1-t\leq e^{-t}$ for $t\in \mathbb{R}$, we have
\begin{align*}
|\E e(\mathbf{X}_p)| &\leq 1 - c \E (\mathbf{X}_p - \theta)^2 \\
&\leq 1 - c \Var \mathbf{X}_p\\
&\leq \exp( - c \Var \mathbf{X}_p )
\end{align*}
for some constant $c> 0$, and hence
\begin{equation}\label{varb}
\sum_{p \in S_1} \Var \mathbf{X}_p \ll \sum_{j=1}^5 \kappa_j.
\end{equation}

To obtain the desired contradiction, it now suffices to show

\begin{theorem}[Technical reduction]\label{thm:tech}  There exists a sequence of $x$ going to infinity for which the following statement holds: there exist natural numbers $K, H, M$ with $M$ even, integers $p_0\geq 1$ and $v_1,\dots,v_K$ with $p_\epsilon\coloneqq p_0+\epsilon_1v_1+\cdots+ \epsilon_Kv_K$ distinct primes for all $\epsilon \in \{0,1\}^K$, a partition of the primes into disjoint sets $S_0,S_1,S_2$, and a random natural number $\mathbf{n}$ obeying~\eqref{pcong},~\eqref{nlower}--- all of which can depend on $x$ --- such that (recalling~\eqref{eq:rhk}) the quantities $\kappa_1,\dots,\kappa_5$ defined in~\eqref{kappa1-def},~\eqref{kappa2-def}, \eqref{kappa3-def},~\eqref{kappa4-def},~\eqref{kappa5-def} obey~\eqref{kappas}, but the random variables $\mathbf{X}_p$ defined in~\eqref{Xp-def} are deterministic for $p \in S_0$ and also obey the condition
\begin{equation}\label{var2}
\sum_{p \in S_1} \Var \mathbf{X}_p \gg 1.
\end{equation}
\end{theorem}

Indeed, it is clear that~\eqref{kappas},~\eqref{varb}, and~\eqref{var2} will contradict each other for $x$ large enough.

\subsection{Selecting parameters}

It remains to establish \Cref{thm:tech}. In the rest of the section we assume that $x$ is sufficiently large.  There is some flexibility in the choice of parameters, but for the sake of concreteness we shall make some (slightly arbitrary) choices.  We define
\begin{align}
K &\coloneqq \left\lfloor \frac{2 \log_4 x}{\log 2} \right\rfloor \label{K-def} \\
H &\coloneqq \lfloor \log_3^2 x\rfloor \label{H-def} \\
M &\coloneqq 2 \lfloor \log_2 x \rfloor \label{M-def}\\
H_+ &\coloneqq \lfloor \log_2^2 x \rfloor \label{Hplus-def} \\
P &\coloneqq \exp( \log^3_3 x ) \label{P-def} \\
L &\coloneqq \log^{1/\log_3 x} x\label{L-def} \\
R &\coloneqq x^{1/\log_2^2 x}\label{Rplus-def} \\
R_{+} &\coloneqq x^{1/100} \label{Rplusplus-def}
\end{align}

Recalling the notations $\lll$ and $\llll$  from Subsection~\ref{notation-sec}, for future reference we observe the hierarchy of scales
\begin{equation}\label{hierarchy}
\begin{split}
1 \llll K \llll  \log_3 x \asymp &(\log_2 R_{+} - \log_2 R) \lll H \asymp 2^K  \\
\lll \log P \asymp \log(H_+ K P) \llll M \asymp\log_2 R &\asymp \log_2 x \lll H_+ \llll 2^H \llll P \llll P^{2^K}  \\
\llll  (HKP)^{2^{2K} H^2} \lll L \lll \log x \llll P^M \llll 2^{H_+}&  \llll R \llll R^M \llll R_{+} \lll x.
\end{split}
\end{equation}
Thus for instance
\begin{equation}\label{P-decay}
\frac{2^{2K} H^2}{P} = \frac{1}{P^{1-o(1)}} \ll \left(o\left(\frac{1}{\log P}\right)\right)^{2^K}. 
\end{equation}

Similarly to \Cref{consecutive-sec}, we divide the shifts $h \geq 1$ into three classes:
\begin{itemize}
    \item \emph{near shifts} $1 \leq h \leq H$;
    \item \emph{far shifts} $H < h \leq H_+$; and
    \item \emph{very far shifts} $h > H_+$.
\end{itemize}

Now we select the parameters $p_0, v_1, \dots, v_K$.

\begin{lemma}  There exist integers $p_0,v_1,\dots,v_K$ with $p_0\geq 1$ such that the integers $p_0 + \epsilon_1 v_1 + \dots + \epsilon_K v_K$ for $\epsilon \in \{0,1\}^K$ are distinct primes in $[P/2,P]$, and furthermore for all $\epsilon \in \{0,1\}^K$ and near shifts $1 \leq h \leq H$, the quantities $r_{\epsilon,h+K}$ (defined in~\eqref{eq:rhk})
are all distinct.
\end{lemma}

\begin{proof} This is an application of the ``Hilbert cube lemma'', which we will utilize via the modern language of Gowers uniformity norms.
There are $\asymp \frac{P}{\log P}$ primes in $[P/2,P]$. Embedding this interval in $\Z/(5P)\Z$, and letting $f \colon \Z/(5P)\Z \to \{0,1\}$ be the indicator function of the embedded primes (that is, $f(x)=\sum_{P/2 \leq p \leq P} 1_{x \equiv p \pmod{5P}}$), we thus have
$$ \|f\|_{U^1(\Z/((5P)\Z)} \gg \frac{1}{\log P},$$
where the Gowers norms are defined in~\eqref{def:gowers}. Thus by monotonicity of the Gowers norms (see, e.g.,~\cite[(11.7)]{tao-vu})
$$ \|f\|_{U^K(\Z/((5P)\Z)} \gg \frac{1}{\log P}.$$
We conclude that if we select $\tilde p_0, \tilde v_1,\dots,\tilde v_K$ uniformly at random from $\Z/(5P)\Z$, the probability that $f(\tilde p_0 + \epsilon_1 \tilde v_1 + \dots + \epsilon_K \tilde v_K)=1$ for each $\epsilon \in \{0,1\}^K$ is
$$ \geq \left(\frac{c}{\log P}\right)^{2^K}$$
for some constant $c>0$.  On the other hand, the probability that two of the $\tilde p_0 + \epsilon_1 \tilde v_1 + \dots + \epsilon_K \tilde v_K$ coincide is
$$ O( 2^{2K} / P ).$$
Similarly, the probability that $\tilde r_{\epsilon,h+K}$ and $\tilde r_{\epsilon',h'+K}$ coincide for some $\epsilon,\epsilon' \in \{0,1\}^K$ and $h,h' \in [H]$ with $(\epsilon,h) \neq (\epsilon',h')$ is
$$ O( 2^{2K} H^2 / P ),$$
where we set
$$\tilde r_{\epsilon,h} \coloneqq \tilde p_\epsilon h - \sum_{k=1}^K k \epsilon_k \tilde v_k.$$

By~\eqref{P-decay}, we conclude that for $x$ large enough, we can select $\tilde p_0,\tilde v_1,\dots,\tilde v_K\in \mathbb{Z}/(5P)\mathbb{Z}$ such that
\begin{itemize}
    \item[(i)] $f(\tilde p_0 + \epsilon_1 \tilde v_1 + \dots + \epsilon_K \tilde v_K)=1$ for each $\epsilon \in \{0,1\}^K$;
    \item[(ii)] No two of the $\tilde p_0 + \epsilon_1 \tilde v_1 + \dots + \epsilon_K \tilde v_K$ for $\epsilon \in \{0,1\}^K$ are equal; and
    \item[(iii)] No two of the $\tilde r_{\epsilon,h+K}$ for $\epsilon \in \{0,1\}^K$ and $h \in [H]$ are equal.
\end{itemize}
We now lift this situation from $\Z/((5P)\Z)$ back to $\Z$.  By definition of $f$ and (i), each of the $\tilde p_0 + \epsilon_1 \tilde v_1 + \dots + \epsilon_K \tilde v_K$ is equal to a prime in $[P/2,P]$.  In particular, $\tilde p_0 \equiv p_0 \pmod{5P}$ for some prime $p_0 \in [P/2,P]$, and on comparing two adjacent elements of the cube $\{0,1\}^K$ we see that $\tilde v_k \equiv v_k \pmod{5P}$ for some integer $-P/2 \leq v_k \leq P/2$.  By traversing the cube one edge at a time, we also see that $p_0 + \epsilon_1 v_1 + \dots + \epsilon_K v_K$ is equal (not just equal mod $5P$) to a prime $p_\epsilon$ in $[P/2,P]$ for every $\epsilon \in \{0,1\}^K$ (the modulus $5P$ is large enough to prevent ``wraparound'' from occurring).  By (ii), the $p_\epsilon$ are all distinct.  By construction we have $\tilde r_{\epsilon,h} \equiv r_{\epsilon,h}$ for all $\epsilon \in \{0,1\}^K$ and $h \in [H]$, and so by (iii) the $r_{\epsilon,h+K}$ are distinct as well, giving the claim.
\end{proof}

Let $p_0,v_1,\dots,v_K$ be as above.  
From the triangle inequality we note the trivial bounds
\begin{equation}\label{rhk}
r_{\epsilon, h+K} \ll (h+K) K P 
\end{equation}
for any $\epsilon \in \{0,1\}^K$ and shifts $h \geq 1$.

We introduce the discriminant
$$W \coloneqq \prod_{\epsilon,\epsilon' \in \{0,1\}^K}\prod_{\substack{1 \leq h,h' \leq H\\ (\epsilon,h) \neq (\epsilon',h')}} |r_{\epsilon,h+K} - r_{\epsilon',h'+K}|.$$
Then
$$ W \leq (O(HKP))^{2^{2K} H^2} \ll \exp( \log^{O(1)}_3 x) \llll L \llll x.$$
Again similarly to \Cref{consecutive-sec}, we divide the set of primes into four classes:
\begin{itemize}
    \item \emph{small primes} in which $p \mid W$ or $p = p_\epsilon$ for some $\epsilon \in \{0,1\}^K$;
    \item \emph{medium primes}\footnote{The terminology here is not completely ordered, as it is possible for medium primes to occasionally be less than some of the small primes; but this will not affect the arguments.} $p \leq R$ not dividing $W$ and not equal to any of the $p_\epsilon$;
    \item \emph{large primes} $R \leq p < R_{+}$; and
    \item \emph{very large primes} $p > R_{+}$;
\end{itemize}
and set $S_0$ to be the set of small primes, $S_1$ to be the set of medium primes, and $S_2$ to be the set of large and very large primes.  Note from~\eqref{hierarchy} that this is indeed a partition of the primes (all the small primes are bounded by $\max(W,P)$ and hence by $R$).

\subsection{Constructing the random variable}

We select the random variable $\mathbf{n}$ uniformly at random among the integers in $[x/2,x]$ that obey the congruences~\eqref{pcong}, as well as the additional congruence condition
\begin{equation}\label{tiny}
\mathbf{n}\equiv 0\pmod W^{\flat}, 
\end{equation}
where
\begin{align*}
W^{\flat}\coloneqq \frac{W}{\prod_{\epsilon\in\{0,1\}^K} p_\epsilon^{v_{p_\epsilon}(W)}}.    
\end{align*}
In particular, $p_{\epsilon}\nmid W^{\flat}$ for all $\epsilon$, so \eqref{tiny} is automatically compatible with the congruences \eqref{pcong}.

Thus, by the Chinese remainder theorem,~\eqref{hierarchy}, and the prime number theorem, $\mathbf{n}$ is restricted to an arithmetic progression $\mathcal{Q}$ of $[x/2,x]$ with spacing $Q$ satisfying
\begin{align}\label{eq:Qbound}
Q \ll P^{2^K} W^{\flat} \ll  \exp( \log^{O(1)}_3 x) \llll L \llll x.    
\end{align}

From this we obtain the relation
\begin{equation}\label{chop}
\E f(\mathbf{n}) = \bar{f} + O\left( \frac{u \|f\|_\infty}{x^{1-o(1)}} \right)
\end{equation}
for any $u$-periodic function $f$ with mean $\bar{f}$ and supremum norm $\|f\|_\infty$, where $u$ is coprime to $Q$.

The condition~\eqref{nlower} is now automatic since
$$ \sum_{k=1}^K k |v_k| \ll K^2 P $$
is much smaller than $x/2$ by~\eqref{hierarchy}.  The congruence condition~\eqref{tiny} also ensures that $\mathbf{X}_p$ is deterministic for all small primes $p \in S_0$.
For non-small primes $p \in S_1 \cup S_2$, we have the bound
\begin{equation}\label{xp-upper}
\mathbf{X}_p = O(1/2^K)
\end{equation}
since all the shifts $r_{\epsilon,h+K}$ reduce to distinct residue classes mod $p$ for $\epsilon \in \{0,1\}^K$ and near shifts $1 \leq h \leq H$, by definition of $W$, and the fact that $p$ does not divide $W^{\flat}$.  In particular,~\eqref{xp} is obeyed.

\subsection{Estimation of \texorpdfstring{$\kappa_1$}{kappa1}}

We first estimate the quantity $\kappa_1$ defined in~\eqref{kappa1-def}. For each $\epsilon \in \{0,1\}^K$ and near shift $1 \leq h \leq H$, the condition $p_\epsilon^2 \mid  \mathbf{n} - \sum_{k=1}^K k \epsilon_k v_k + p_\epsilon h$ restricts $\mathbf{n}$ to lie in a subprogression of $\mathcal{Q}$ of $p_\epsilon$ times the spacing, and thus
$$ \Prob( p_\epsilon^2 \mid  \mathbf{n} + r_{\epsilon,h}) \ll \frac{1}{p_\epsilon} \ll \frac{1}{P}.$$
For far or very far shifts $h > H$, we trivially have
$$ \Prob( p_\epsilon^2 \mid  \mathbf{n} + r_{\epsilon,h}) \leq 1.$$
Inserting these bounds into~\eqref{kappa1-def}, we conclude that
$$ \kappa_1 \ll 2^K \left( \frac{1}{P} + \frac{1}{2^H}\right).$$
By~\eqref{hierarchy}, we conclude that the contribution of $\kappa_1$ to~\eqref{kappas} is acceptable.

\subsection{Estimation of \texorpdfstring{$\kappa_2$}{kappa2}}

Now we estimate $\kappa_2$, given by~\eqref{kappa2-def}.  For the very far shifts $h > H_+$ we use the very crude bound $\omega(n) \ll \log n$, so that by~\eqref{rhk} we have
$$ \omega(\mathbf{n}+r_{\epsilon,h+K}) \ll \log x+\log(r_{\epsilon,h+K})\ll \log x+\log(hKP)$$
(noting that $h+K \asymp h$ for $h > H$) which implies that the contribution of very far shifts $h > H_+$ to~\eqref{kappa2-def} is 
$$\ll \frac{\log x + \log( P H_+ )}{2^{H_+}}.$$
For the far shifts $H < h \leq H_+$, we observe that as $\mathbf{n} + r_{\epsilon,h+K} = O(x)$ can have at most $O(1)$ prime factors larger than $R_+$, we have
$$ \omega_{>R_+}(\mathbf{n} + r_{\epsilon,h+K}) \ll 1,$$
and from the congruence condition~\eqref{pcong} we have
$$\sum_{p\text{ small}} 1_{p \mid \mathbf{n} + r_{\epsilon,h+K}} = \omega(r_{\epsilon,h+K}) \ll \log r_{\epsilon,h+K} \ll \log(H_+ K P)$$
thanks to~\eqref{rhk}.  Thus by~\eqref{fadd} we have
$$ \omega(\mathbf{n} + r_{\epsilon,h+K}) \ll \sum_{p \text{ medium or large}} 1_{p\mid \mathbf{n}+r_{\epsilon,h+K}} + \log(H_+ K P).$$ 
For medium or large primes $p$ (so in particular $p \leq R_+$ and $p \nmid W$) we see from~\eqref{chop} (and~\eqref{hierarchy}) that
$$ \E 1_{p\mid \mathbf{n}+r_{\epsilon,h+K}} \ll \frac{1}{p} + \frac{p}{x^{1-o(1)}} \ll \frac{1}{p}.$$
  Thus by Mertens' theorem we conclude that
$$ \E \omega(\mathbf{n}+r_{\epsilon,h+K}) \ll \log_2 R_+ + \log(H_+ K P)$$
and hence
$$ \kappa_2 \ll \frac{\log_2 R_+ + \log(H_+ K P)}{2^H} + \frac{\log x + \log(H_+ K P)}{2^{H_+}}.$$
By~\eqref{hierarchy}, we conclude that the contribution of $\kappa_2$ to~\eqref{kappas} is acceptable.

\subsection{Estimation of \texorpdfstring{$\kappa_4$}{kappa4}}

We postpone the treatment of $\kappa_3$ for now and move on to $\kappa_4$ defined in~\eqref{kappa4-def}.
For $0 \leq m \leq M$ and medium primes $p_1,\dots,p_m \in S_1$, the product $p_1 \cdots p_m$ is at most $R^M$.  From~\eqref{chop},~\eqref{xp-upper} we conclude that
$$ \E \mathbf{X}_{p_1} \cdots \mathbf{X}_{p_m} - \E \mathbf{X}'_{p_1} \cdots \mathbf{X}'_{p_m}\ll \frac{R^M}{2^K x^{1-o(1)}} \ll \frac{1}{x^{1-o(1)}}$$
for all such $p_1,\dots,p_m$, thanks to~\eqref{hierarchy}.  Inserting this into~\eqref{kappa4-def}, we conclude that
$$ \kappa_4 \ll M (2\pi)^M P^M \frac{1}{x^{1-o(1)}}.$$
By~\eqref{hierarchy}, we conclude that the contribution of $\kappa_4$ to~\eqref{kappas} is acceptable.

\subsection{Estimation of \texorpdfstring{$\kappa_5$}{kappa5}}

Now we estimate $\kappa_5$ defined in~\eqref{kappa5-def}.  For any medium prime $p \in S_1$, we see from~\eqref{xp-upper} (or~\eqref{xp}) that
$$ \exp( \pm 2\pi e \mathbf{X}_p ) = 1 \pm 2\pi e \mathbf{X}_p + O( |\mathbf{X}_p|^2)$$
while from~\eqref{chop} (and~\eqref{hierarchy}) one has
\begin{equation}\label{xp-1}
 \E \mathbf{X}_p \ll \frac{p}{x^{1-o(1)}} = \frac{1}{x^{1-o(1)}}
 \end{equation}
and
\begin{equation}\label{xp-2}
\E |\mathbf{X}_p|^2 \asymp \frac{1}{2^K p} + O\left(\frac{p}{x^{1-o(1)}}\right) \asymp \frac{1}{2^K p},
\end{equation}
so that
$$ \E \exp( \pm 2\pi e\mathbf{X}_p ) = \exp\left( O\left( \frac{1}{2^K p} \right) \right),$$
and hence by Mertens' theorem and~\eqref{kappa5-def} we have
$$ \kappa_5 \ll \exp\left( O\left( \frac{\log_2 R}{2^K} \right) - M \right).$$
By~\eqref{hierarchy}, we conclude that the contribution of $\kappa_4$ to~\eqref{kappas} is acceptable.

\subsection{Variance bound}
From~\eqref{xp-1},~\eqref{xp-2} (and~\eqref{hierarchy}) we have
$$ \Var \mathbf{X}_p \asymp \frac{1}{2^K p},$$
so by Mertens' theorem
$$ \sum_{p \in S_1} \Var \mathbf{X}_p \gg \frac{\log_2 R - \log_2 W}{2^K} $$
and the claim~\eqref{var2} follows from~\eqref{hierarchy}.

\subsection{Estimation of \texorpdfstring{$\kappa_3$}{kappa3}}

The last remaining task required for \Cref{thm:tech} is to establish $\kappa_3 = o(1)$ given by~\eqref{kappa3-def}.  By~\eqref{kappa3-def} and the triangle inequality, we may split this into two subtasks, corresponding to the large and very large primes respectively:
\begin{align}
 \left|\E \sum_{R < p \leq R_{+}} \mathbf{X}_p\right|^2 &= o(1) \label{med-est} \\
 \left|\E \sum_{R_{+} < p \leq 2x} \mathbf{X}_p\right|^2 &= o(1) \label{hi-est}
\end{align}
(note that $\mathbf{X}_p$ will necessarily vanish for $p> 2x$).

\subsection{Contribution of large primes to \texorpdfstring{$\kappa_3$}{kappa3}}  
We now prove~\eqref{med-est}.  For distinct large primes $R < p,p' \leq R_+$, 
the product $\mathbf{X}_{p} \mathbf{X}_{p'}$ is periodic in $\mathbf{n}$ with period $p p' \leq R_+^2$,  bounded by~\eqref{xp-upper}, and mean zero thanks to the alternating factor $(-1)^{|\epsilon|}$ in~\eqref{Xp-def}, so from~\eqref{chop} we have
$$
\E \mathbf{X}_{p} \mathbf{X}_{p'} = O\left( \frac{R_+^2}{x^{1-o(1)}} \right).$$
Hence, on expanding out the square, we have
$$  \E \left|\sum_{R < p \leq R_{+}} \mathbf{X}_p\right|^2 =  \sum_{R < p \leq R_+} \E |\mathbf{X}_p|^2 + O\left( \frac{R_{+}^4}{x^{1-o(1)}} \right).$$
There are no collisions between the various congruence classes defining $\mathbf{X}_p$ in~\eqref{Xp-def}, and so by~\eqref{chop} one has
$$ \E |\mathbf{X}_p|^2 \ll \frac{1}{p2^K} + \frac{R_{+}}{x^{1-o(1)}}.$$
Hence by Mertens' theorem, we can bound the left-hand side of~\eqref{med-est} by
$$ \ll \frac{R_{+}^4}{x^{1-o(1)}} + \frac{\log_2 R_{+} - \log_2 R}{2^K}.$$
By~\eqref{hierarchy}, we conclude~\eqref{med-est}.

\subsection{Contribution of very large primes to \texorpdfstring{$\kappa_3$}{kappa3}}  We now prove~\eqref{hi-est}.  By inserting an alternating sum of $\frac{1}{p}$ and applying the Cauchy--Schwarz inequality, it will suffice to show that
$$ \E \left|\sum_{\epsilon \in \{0,1\}^K} (-1)^{|\epsilon|} \sum_{R_{+} < p \leq 2x} \left(1_{p\mid \mathbf{n}+r_{\epsilon,h+K}}-\frac{1}{p}\right)\right|^2 = o( 2^{2K})$$
for all near shifts $1 \leq h \leq H$. 

Fix $h$.  On squaring, it suffices to show that
\begin{equation}\label{easy}
 \E \left|\sum_{R_{+} < p \leq 2x} \left(1_{p\mid \mathbf{n}+r_{\epsilon,h+K}}-\frac{1}{p}\right)\right|^2 \ll 1
 \end{equation}
for $\epsilon \in \{0,1\}^K$ and
\begin{align}\label{eq:largeprimes_estimate}
\E \left(\sum_{R_{+} < p \leq 2x} \left(1_{p\mid \mathbf{n}+r_{\epsilon,h+K}} - \frac{1}{p}\right)\right) \left(\sum_{R_{+} < p' < 2x} \left(1_{p'\mid \mathbf{n}+r_{\epsilon',h+K}}-\frac{1}{p'}\right)\right) = o(1)
\end{align}
for distinct $\epsilon, \epsilon' \in \{0,1\}^K$.  The former claim~\eqref{easy} is easy since the expression in parentheses is bounded, so we turn to the latter claim~\eqref{eq:largeprimes_estimate}.

Partition the interval $(R_+,2x]$ into $L$ subintervals $J_\ell \coloneqq (y_{\ell-1}, y_{\ell}]$ for $\ell \in [L]$, where
$$ R_+ = y_0 < y_1 < \dots < y_L = 2x$$
and $\log_2 y_\ell - \log_2 y_{\ell-1} \asymp \frac{\log_2 x}{L}$ for all $\ell \in [L]$; for instance one could define $y_\ell$ by the formula
$$ \log_2 y_\ell \coloneqq \log_2 y_0 + \frac{\ell}{L} (\log_2(2x) - \log_2 y_0).$$
For any $n = O(x)$, we can then write
\begin{equation}\label{mean}
\sum_{R_+<p < 2x}1_{p\mid n} = \sum_{\ell \in [L]} \left( 1 - g_\ell(n) + O(1_{n\in E_\ell})\right)
\end{equation}
where, for each $\ell \in [L]$, $g_\ell$ is the completely multiplicative function with $g_\ell(p) \coloneqq 1_{p \notin J_\ell}$ for all primes $p$, and $E_\ell$ is the set of $n\in [x/2,x]\cap \mathbb{N}$ for which there exist two primes 
$$ y_{\ell-1} < p \leq p' \leq y_{\ell}$$
which both divide $n$.  

The sets $E_\ell$ are quite sparse: a standard application of Mertens' theorem shows that
\begin{equation}\label{standard}
 \frac{1}{x/2} \sum_{x/2 \leq n \leq x} 1_{n\in E_\ell} \ll \left(\frac{\log_2 x}{L}\right)^2,
 \end{equation}
so if we set
$$ \delta_{\ell,x} \coloneqq \frac{1}{x/2} \sum_{x/2 \leq n \leq x} g_\ell(n),$$
then on taking means in~\eqref{mean}
$$ \sum_{R_+<p < 2x} \frac{1}{p} = \sum_{\ell \in [L]} \left( 1 - \delta_{\ell,x} + O\left(\left(\frac{\log_2 x}{L}\right)^2\right)\right).$$
Subtracting, we conclude that
\begin{equation}\label{pm}
 \sum_{R_+<p < 2x} \left( 1_{p\mid n} - \frac{1}{p} \right) = -\sum_{\ell \in [L]} \left( g_\ell(n) - \delta_{\ell,x} + O(1_{n\in E_\ell})\right) + o(1).
 \end{equation}
In particular, we apply~\eqref{pm} with $n$ set to $\mathbf{n}+r_{\epsilon,h+K}$ and insert it into~\eqref{eq:largeprimes_estimate}.  The contribution of a given $E_\ell$ term is then bounded by
$$ \ll \E 1_{\mathbf{n}+r_{\epsilon,h+K}\in E_{\ell}} \ll \frac{Q(\log_2 x)^2}{L^2}$$
where $Q$ is the period of the residue class defining $\mathbf{n}$, and where we have used~\eqref{standard} and the fact that sums over $p,p'$ in~\eqref{eq:largeprimes_estimate} are bounded.  Thus the net effect of all of the $E_\ell$ terms is $O(Q(\log_2 x)^2/L) = o(1)$ by~\eqref{hierarchy}, so these contributions can be discarded.  Similarly if we apply~\eqref{pm} with $n$ set to
$\mathbf{n}+r_{\epsilon',h+K}$.  We may thus replace the left-hand side of~\eqref{eq:largeprimes_estimate} by
$$
\sum_{\ell,\ell' \in [L]} \E (g_\ell(\mathbf{n}+r_{\epsilon,h+K}) - \delta_{\ell,x}) (g_{\ell'}(\mathbf{n}+r_{\epsilon',h+K}) - \delta_{\ell',x}).$$
By the triangle inequality and~\eqref{hierarchy}, it thus suffices to show that
$$ \E (g_\ell(\mathbf{n}+r_{\epsilon,h+K}) - \delta_{\ell,x}) (g_{\ell'}(\mathbf{n}+r_{\epsilon',h+K}) - \delta_{\ell',x}) \ll \log^{-c} x$$
for some absolute constant $c>0$ and all $\ell,\ell' \in [L]$.

Applying \Cref{thm:correlation} (in the equidistributed case with ${\mathcal L}=\log^c x$) with $g_1$ set to $g_\ell$ and $g_2$ set to either $g_{\ell'}$ or $1$, and noting that importantly the spacing $Q$ of the progression to which $\mathbf{n}$ is restricted is $\ll \log^{o(1)} x$ by~\eqref{eq:Qbound}, we see that it will suffice (after restricting $x$ to a subsequence) to show that
$$
\sum_{\substack{N<n\leq 2N\\ n\equiv a\pmod r}}g_\ell(n) = \frac{N}{r}\delta_{\ell,N} + O(N \log^{-c} x)
$$
for some absolute constant $c>0$ and all $x^{0.4} \leq N \leq x$ and $a,r\in \mathbb{N}$, and similarly for $g_{\ell'}$.  Of course it suffices to show the $g_\ell$ bound.  By dividing out common factors of $a,r$ (and using~\eqref{lipo}) we may assume without loss of generality that $(a,r)=1$ (after worsening the lower bound on $N$ to say $x^{0.3} \leq N$).
We may assume that $r = O(\log^c x)$, as the claim is trivial otherwise.

By inclusion-exclusion and the fact that $J_\ell$ lies in $(R_+,2x]$, we have
$$ \frac{1}{N} \sum_{N < n \leq 2N} g_\ell(n) = \sum_{i=0}^{101} (-1)^i \sum_{\substack{p_1 < \dots < p_i\\ p_1,\dots,p_i \in J_\ell\\ p_1 \cdots p_i \leq 2N}} \left( \frac{1}{p_1 \cdots p_i} + O(\frac{1}{N}) \right) $$
and similarly
$$ \frac{r}{N} \sum_{\substack{N < n \leq 2N\\ n\equiv a\pmod r}} g_\ell(n) = \sum_{i=0}^{101} (-1)^i \sum_{\substack{p_1 < \dots < p_i\\ p_1,\dots,p_i \in J_\ell\\ p_1 \cdots p_i \leq 2N}} \left( \frac{1}{p_1 \cdots p_i} + O\left(\frac{r}{N}\right) \right).$$
Thus it will suffice to show that for all $0 \leq i \leq 101$ that
$$ \sum_{\substack{p_1 < \dots < p_i\\ p_1,\dots,p_i \in J_\ell\\ p_1 \cdots p_i \leq 2N}} \frac{r}{N} \ll \log^{-c} N.$$
But by a simple application of Selberg's sieve there are only $O(N \log^{-1} N)$ numbers $n \leq 2N$ that are not divisible by any prime less than or equal to $R_+$, so the number of summands here is $N \log^{-1} N$, and the claim follows since $r = O(\log^c x)$.

\appendix

\section{Empirical distribution of consecutive values of prime counting functions}\label{appendix:a}

We conclude this paper by presenting some figures on the convergence of the probabilities $\Prob(f(\mathbf{n})=f(\mathbf{n}+1))$ for $f=\omega, \Omega, \tau$ to their conjectured asymptotics with lower order terms, relating to Theorem~\ref{consecutive-omega-thm}. The initial codes for these plots were provided by ChatGPT.

 \begin{figure}[H]
    \centering
\centerline{\includegraphics[width=\textwidth]{ctau.png}}
    \caption{The probability of the event that $\tau(\mathbf{n}+1)/\tau(\mathbf{n})$ is a power of $2$ numerically converges (fairly rapidly) to about $0.4888$.  For comparison we also include the larger events $\nu_p(\tau(\mathbf{n}+1)/\tau(\mathbf{n}))=0$ for $p=3$, $p=3,5$, and $p=3,5,7$, which are easier to calculate asymptotically (they correspond to return probabilities of certain one-dimensional, two-dimensional, and three-dimensional random walks respectively), as well as the two smaller events that $\tau(\mathbf{n}), \tau(\mathbf{n}+1)$ are both powers of two, or either both powers of two or both three times powers of two.  The parameter $x$ ranges between $10$ and $10^7$ (and is restricted to multiples of $10$ or $10^2$ for $x \geq 10^3$ and $x \geq 10^4$ respectively) and is plotted logarithmically.  }
    \label{ctau-fig}
\end{figure}

\bibliography{refs}
\bibliographystyle{plain}

\end{document}